\documentclass{article}
\usepackage{amssymb}
\usepackage{amsfonts}
\usepackage{hyperref}
\usepackage{amsmath}

\newtheorem{theorem}{Theorem}

\newtheorem{corollary}[theorem]{Corollary}

\newtheorem{definition}[theorem]{Definition}

\newtheorem{lemma}[theorem]{Lemma}

\newtheorem{proposition}[theorem]{Proposition}
\newtheorem{remark}[theorem]{Remark}

\numberwithin{theorem}{section}
\input{tcilatex}
\begin{document}

\title{Cohomological Weight Shiftings for Automorphic Forms on Definite
Quaternion Algebras}
\author{Davide A. Reduzzi \\
%EndAName
University of California at Los Angeles\\
{\scriptsize devredu83@math.ucla.edu}}
\maketitle

\begin{abstract}
Let $F/%
%TCIMACRO{\U{211a} }%
%BeginExpansion
\mathbb{Q}
%EndExpansion
$ be a totally real field extension of degree $g$ and let $D$ be a definite
quaternion algebra with center $F$. Fix an odd prime $p$ which is unramified
in $F$ and $D$. We produce weight shiftings between $(\func{mod}p)$
automorphic forms on $D^{\times }$ of a fixed level $U$. When the starting
weight does not contain any $(2,...,2)$-block, we obtain these shiftings via
maps induced in cohomology by intertwining operators acting on $\mathbb{\bar{%
F}}_{p}$-representations of $GL_{2}(\mathcal{O}_{F}/p\mathcal{O}_{F})$. We
construct two families of such operators, each of cardinality $g^{2}$, and
we produce between others weight shiftings by cyclic permutations of the
blocks $(p^{r},0,...,0,\pm 1,0,...,0)$ where the number of zeros between $%
p^{r}$ and $\pm 1$ depends upon the value of the integer $r$. In particular,
we produce shiftings by $(p,\pm 1,0,...,0)$.

\noindent Shiftings by $(p-1,...,p-1)$ for weights containing $(2,...,2)$%
-blocks are obtained following the methods of \cite{EK}, as half of our
intertwining operators becomes trivial in this case.
\end{abstract}

\tableofcontents

\section{Introduction}

Let $F$ be a totally real number field and $p$ be an odd prime which is
unramified in $F$. Let $D$ be a definite quaternion algebra with center $F$
and assume that $D$ is split at all primes of $F$ above $p$. In this paper
we produce congruences modulo $p$ between automorphic forms on $D^{\times }$
having fixed level and varying weights. Our interest in this matter is
partially motivated by the study of the weight part of\ Serre's modularity
conjecture over a totally real field, as formulated in \cite{BDJ} and proven
in many cases in \cite{Gee22} and \cite{Gee333}. We hope that our
constructions could be further generalizable to algebraic groups other than $%
D^{\times }.$

Weight shiftings for modulo $p$ elliptic modular forms can be obtained via
the classical theory of Hasse invariants and theta operators (\cite{Serre333}%
), and have been studied via cohomological methods by Ash-Stevens (\cite{AS}%
) and Edixhoven-Khare (\cite{EK}). In \cite{Re}, the author studied
cohomological weight shiftings of Hasse-type, adopting the viewpoint
conceived by C. Khare on this matter. \noindent In \cite{EK}, Edixhoven and
Khare produce parallel weight shiftings by $p-1$ for automorphic forms on $%
D^{\times }$ of parallel weight two, under the assumption that $p$ is inert
in $F$. In \cite{Gee22}, Gee uses a construction of Kisin (\cite{Kisin3})
relying on the classification of the irreducible admissible $\mathbb{\bar{F}}%
_{p}$-representations of $GL_{2}(%
%TCIMACRO{\U{211a} }%
%BeginExpansion
\mathbb{Q}
%EndExpansion
_{p})$ to produce some non-parallel weight shiftings by $p-1$ on forms on $%
D^{\times }$, assuming that $p$ is totally split in $F$: this is a crucial
step in his proof of the weight conjecture of \cite{BDJ} in the totally
split case. Weight shiftings for geometric Hilbert modular forms over $F$
can be obtained via operators constructed by Andreatta and Goren in \cite{AG}%
.

As these results show, there are two possible approaches to the study of
congruences between automorphic forms of different weights: geometric and
cohomological. Assume that $p>3$ and that $N\geq 5$ is coprime to $p$. The
theta operator and the Hasse invariant are geometrically defined operators
acting on spaces of $(\func{mod}p)$ elliptic modular forms of level $N$. The
first operator induces a Hecke equivariant injection increasing weights by $%
p-1$; the latter is a Hecke twist-equivariant map that shifts weights by $%
p+1 $. The geometric approach to weight shiftings passes through the
generalization of these operators to spaces of geometric Hilbert modular
forms over the totally real field $F$. This is carried over by Goren in \cite%
{Go}, where $[F:%
%TCIMACRO{\U{211a} }%
%BeginExpansion
\mathbb{Q}
%EndExpansion
]$ partial Hasse invariants are constructed assuming $p$ is unramified in $F$%
, and by Andreatta-Goren in \cite{AG}, where generalized theta operators are
considered and the unramifiedness assumption of $p$ in $F$ is dropped.

The Eichler-Shimura isomorphism translates the study of Hecke eigensystems
of $(\func{mod}p)$ elliptic modular forms of weight $k\geq 2$ and level $N$
into the study of the Hecke action on the cohomology group $H^{1}(\Gamma
_{1}(N),\limfunc{Sym}^{k-2}\mathbb{\bar{F}}_{p}^{2})$. In \cite{AS}, Ash and
Stevens identify a cohomological analogue of the theta operator in the map
induced in cohomology by the Dickson polynomial $\Theta
_{p}=XY^{p}-X^{p}Y\in \mathbb{F}_{p}\left[ X,Y\right] .$ For the
cohomological counterpart of the Hasse invariant, we must restrict ourselves
to work with $p$-small weights. Additionally, the case of weight two$\ $must
be treated \textit{per se }(this dichotomy between weight two and weight
larger than two will appear, \textit{mutatis mutandis}, also in this paper).
In \cite{EK}, a cohomological analogue of the Hasse invariant acting upon
weight two forms is constructed by studying a degeneracy map:%
\begin{equation*}
H^{1}\left( \Gamma _{1}(N),\mathbb{\bar{F}}_{p}\right) ^{2}\longrightarrow
H^{1}\left( \Gamma _{1}(N)\cap \Gamma _{0}\left( p\right) ,\limfunc{Sym}%
\nolimits^{p-1}\mathbb{\bar{F}}_{p}^{2}\right) .
\end{equation*}%
A $GL_{2}(\mathbb{F}_{p})$-equivariant derivation $D$ of $\mathbb{F}%
_{p}[X,Y] $ defined by Serre by:%
\begin{equation*}
\text{ }Df=X^{p}\partial _{X}f+Y^{p}\partial _{Y}f
\end{equation*}

\noindent is used in \cite{Re} \ to produce weight shiftings by $p-1$
starting from forms of weight $2<k\leq p+1$. The cokernels of the operators $%
\Theta _{p}$ and $D$ are related to the characteristic zero theory of
representations of $GL_{2}(\mathbb{F}_{p})$ (cf. \cite{AS}, \cite{Re}).

Fred Diamond suggested to look for a generalization of the results of \cite%
{Re} to other contexts. Following his suggestion, in this paper we construct
weight shiftings for Hilbert modular forms in cohomological settings. The
geometric picture had a motivational role in our study. 

Let us mention first a few advantages of working with cohomology groups
rather than with geometrically defined modular forms. First, by the
Jacquet-Langlands correspondence, we are led to the more general study of
weight shiftings for adelic automorphic forms on definite quaternion $F$%
-algebras, where $F$ is a field as above. These adelic spaces (cf. \ref%
{adelic HMF}) seem to be well suited for computations. Furthermore, their
formation is compatible with base change (Proposition \ref{onto}), while the
formation of spaces of geometric Hilbert modular forms is not, in general:
none of the geometrically constructed partial Hasse invariants lift to
characteristic zero if $F\neq 
%TCIMACRO{\U{211a} }%
%BeginExpansion
\mathbb{Q}
%EndExpansion
$. Finally, our methods produce a larger variety of weight shiftings than
the ones arising from the geometric setup (cf. Remark\ \ref{more}).

\bigskip

The paper is divided into two parts: in the first part, consisting of
sections 2 and 3, we study weight shiftings for Serre's weights, i.e., for
(irreducible) $\overline{\mathbb{F}}_{p}$-linear representations of $GL_{2}(%
\mathcal{O}_{F}/p\mathcal{O}_{F})$; in the second part, consisting of
sections 4 and 5, we address the problem of weight shiftings for automorphic
forms associated to the definite quaternion algebra $D$ with center $F$.

Set $[F:%
%TCIMACRO{\U{211a} }%
%BeginExpansion
\mathbb{Q}
%EndExpansion
]=g_{\ast }$.\ The main novelties of the paper consist in: (1) the
introduction of $g_{\ast }^{2}$ generalized Dickson operators and $g_{\ast
}^{2}$ generalized $D$-operators acting on $\overline{\mathbb{F}}_{p}[GL_{2}(%
\mathcal{O}_{F}/p\mathcal{O}_{F})]$-modules: these maps will induce
cohomological generalizations of the theta operators and of the partial
Hasse invariants, respectively; (2) the determination of many non-parallel
weight shiftings for automorphic forms on $D^{\times }$ of a fixed level. In
particular, we will answer a question of Diamond, as for any prime $%
\mathfrak{P}$ of $F$ above $p$ we will produce weight shiftings that
increase those entries of the weight parameter $\vec{k}$ associated to the
embeddings $F_{\mathfrak{P}}\hookrightarrow \bar{%
%TCIMACRO{\U{211a}}%
%BeginExpansion
\mathbb{Q}%
%EndExpansion
}_{p}$ by $(p,-1,0,...,0)\in 
%TCIMACRO{\U{2124} }%
%BeginExpansion
\mathbb{Z}
%EndExpansion
^{f(\mathfrak{P}/p)}$ - or by any cyclic permutation of this tuple, cf. \ref%
{GDO}.

\noindent Along the way, we will also obtain new identities between virtual
modular representations of $GL_{2}(\mathcal{O}_{F}/p\mathcal{O}_{F})$ that
will allow us to give an algorithm to compute the Jordan-H\"{o}lder
constituents of any product of symmetric power representations of this group.

Let us now summarize the content of each section of the paper. Fix a
positive integer $g$ and set $q=p^{g}$ and $G=GL_{2}\left( \mathbb{F}%
_{q}\right) $. For any non-negative integer $k$, define the $\mathbb{F}%
_{q}[G]$-module $M_{k}=\limfunc{Sym}^{k}\mathbb{F}_{q}^{2}$. In \cite{Se01},
after extending the definition of the $M_{k}$'s for $k<0$ in a suitable way,
Serre proves the following identity, valid in the Grothendieck ring of
finitely generated $\mathbb{F}_{q}[G]$-modules for any integer $k$:

\begin{equation*}
M_{k}-\det \cdot M_{k-(q+1)}=M_{k-\left( q-1\right) }-\det \cdot M_{k-2q}.
\end{equation*}

\noindent The weight shiftings by $q-1$ and by $q+1$ appearing in the above
formula are induced by the Dickson invariant $\Theta _{q}$ and the
derivation map $D$ mentioned above. In section 2 we recall some
constructions associated to these operators and some weight shiftings
results for elliptic modular forms (cf. \cite{Re}).

In section 3 we derive the new identity (Corollary \ref{inttt}):%
\begin{equation}
M_{k}^{[i]}M_{h}^{[i+1]}-\det\nolimits^{p^{i+1}}\cdot
M_{k-p}^{[i]}M_{h-1}^{[i+1]}=M_{k-p}^{[i]}M_{h+1}^{[i+1]}-\det%
\nolimits^{p^{i+1}}\cdot M_{k-2p}^{[i]}M_{h}^{[i+1]}  \tag{1}
\end{equation}

\noindent valid for any $h,k,i\in 
%TCIMACRO{\U{2124} }%
%BeginExpansion
\mathbb{Z}
%EndExpansion
$. Here the superscript $[i]$ indicates that the $G$-action is twisted by
the $i$th power of the absolute Frobenius morphism of $\mathbb{F}_{q}$. For $%
g>1$, this identity allows us to explicitly compute the Jordan-H\"{o}lder
factors of any virtual representations of the form $\tprod%
\nolimits_{i=0}^{g-1}M_{k_{i}}^{[i]}$ (Theorem \ref{final}; cf. also \cite%
{Re3}).

\noindent As for the case $g=1$, also for $g>1$ the periods appearing in $%
(1) $, i.e., the cyclic permutations of the $g$-tuples $(p,1,0,...,0)$ and $%
(p,-1,0,...,0)$, correspond to weight shiftings arising from $G$-equivariant
operators. In \ref{GDI} and \ref{GDO}, we define two families of such
operators, each containing $g^{2}$ maps. For any integers $\alpha ,\beta $
subject to the constraints $0\leq \alpha \leq g-1$ and $1\leq \beta \leq g-1$%
, we construct generalized Dickson operators\ $\Theta _{\beta }^{[\alpha ]}$
and generalized $D$-operators $D_{\beta }^{[\alpha ]}$ giving rise, for any
set of non-negative integers $k_{0},...,k_{g-1}$, to the $G$-modules
monomorphism:%
\begin{equation*}
\Theta _{\beta }^{[\alpha ]}:\det\nolimits^{p^{\alpha }}\otimes
\dbigotimes\nolimits_{i}M_{k_{i}}^{[i]}\longrightarrow \left(
\dbigotimes\nolimits_{i\neq \alpha ,\alpha +\beta }M_{k_{i}}^{[i]}\right)
\otimes M_{k_{\alpha }+1}^{[\alpha ]}\otimes M_{k_{\alpha +\beta
}+p^{g-\beta }}^{[\alpha +\beta ]}
\end{equation*}

\noindent and to the $G$-morphism:%
\begin{equation*}
D_{\beta }^{[\alpha
]}:\dbigotimes\nolimits_{i}M_{k_{i}}^{[i]}\longrightarrow \left(
\dbigotimes\nolimits_{i\neq \alpha ,\alpha +\beta }M_{k_{i}}^{[i]}\right)
\otimes M_{k_{\alpha }-1}^{[\alpha ]}\otimes M_{k_{\alpha +\beta
}+p^{g-\beta }}^{[\alpha +\beta ]}.
\end{equation*}

\noindent We study some properties of these and other $G$-operators at the
end of section 3.

In section 4 we use the above results to obtain weight shiftings for modulo $%
p$ automorphic forms on $D^{\times }$ having fixed level. We start by\
treating the case in which the tensor factors - corresponding to the prime
decomposition of $p$ in $F$ - of the weight that we want to shift are all of
dimension greater than one: this is what we call a weight not containing a $%
(2,...,2)$-block.

More precisely, write $p\mathcal{O}_{F}=\mathfrak{P}_{1}...\mathfrak{P}_{r}$
and denote by $f_{j}$ the residual degree of $\mathfrak{P}_{j}$ over $p$.
Let $\mathcal{O}$ be the ring of integers of the smallest unramified
extension of $%
%TCIMACRO{\U{211a} }%
%BeginExpansion
\mathbb{Q}
%EndExpansion
_{p}$ inside $\bar{%
%TCIMACRO{\U{211a}}%
%BeginExpansion
\mathbb{Q}%
%EndExpansion
}_{p}$\ containing the image of $F$ under all the embeddings $F\rightarrow 
\bar{%
%TCIMACRO{\U{211a}}%
%BeginExpansion
\mathbb{Q}%
%EndExpansion
}_{p}$; let $\mathbb{F}$ be the residue field of $\mathcal{O}$. Let $A$ be a
topological $%
%TCIMACRO{\U{2124} }%
%BeginExpansion
\mathbb{Z}
%EndExpansion
_{p}$-algebra and denote by $S_{\tau ,\psi }(U,A)$ the space of $A$-valued
adelic automorphic forms on $D$ having level $U\subset \left( D\otimes _{F}%
\mathbb{A}_{F}^{\infty }\right) ^{\times }$, weight $\tau :U\rightarrow 
\limfunc{Aut}(W_{\tau })$ and Hecke character $\psi :\left( \mathbb{A}%
_{F}^{\infty }\right) ^{\times }/F^{\times }\rightarrow A^{\times }$. For
any set $S$ of primes\ of $F$\ containing the ramification set of $D$, the
primes above $p$ and the primes $v$ for which $U_{v}$ is not a maximal
compact subgroup of $D_{v}^{\times }$, the universal Hecke algebra $\mathbb{T%
}_{S,A}^{univ}=A[T_{v},S_{v}:v\notin S]$ acts upon this space. We assume
that $U$ is small enough.

We use the generalized Dickson and $D$-operators from section $3$, together
with some classical results of Ash-Stevens and Deligne-Serre, to produce
congruences modulo $p$ between Hecke eigenforms arising from the spaces $%
S_{\tau ,\psi }(U,\mathcal{\bar{%
%TCIMACRO{\U{2124}}%
%BeginExpansion
\mathbb{Z}%
%EndExpansion
}}_{p})$ for fixed $U$ and varying $\tau $. Some technical difficulties
arise, as we are mainly interested in forms having holomorphic weight in the
sense of \ref{conv_on_emb}, but, in general, our intertwining operators do
not preserve holomorphicity. Furthermore, we need to make sure that when we
transfer forms from weight $\tau $ to weight $\tau ^{\prime }$, we can lift
the reduction of $\psi $ modulo the maximal ideal of $\mathcal{\bar{%
%TCIMACRO{\U{2124}}%
%BeginExpansion
\mathbb{Z}%
%EndExpansion
}}_{p}$ to a compatible $\mathcal{\bar{%
%TCIMACRO{\U{2124}}%
%BeginExpansion
\mathbb{Z}%
%EndExpansion
}}_{p}^{\times }$-valued Hecke character for $\tau ^{\prime }$.

One of the weight shiftings result we can prove is the following (Theorem %
\ref{general}). Assume that $\tau $ is the $\mathcal{O}$-linear weight with
holomorphic parameters $(\vec{k},w)\in 
%TCIMACRO{\U{2124} }%
%BeginExpansion
\mathbb{Z}
%EndExpansion
_{\geq 2}^{g}\times \left( 2%
%TCIMACRO{\U{2124} }%
%BeginExpansion
\mathbb{Z}
%EndExpansion
+1\right) $ and that $\psi $ is a Hecke character compatible with $\tau $.
Let $f$ be the minimum of the residual degrees of the primes $\mathfrak{P}%
_{j}$ and fix an integer $\beta $ such that $1\leq \beta \leq f$. For any
integers $i$ and $j$ with $1\leq j\leq r$ and $0\leq i\leq f_{j-1}\ $choose $%
a_{i}^{(j)}\in \{p^{\beta }-1,p^{\beta }+1\}$ and set $\vec{a}=(\vec{a}%
^{(1)},...,\vec{a}^{(r)})$ with $\vec{a}%
^{(j)}=(a_{0}^{(j)},...,a_{f_{j}-1}^{(j)})$, and $w^{\prime }=w+(p^{\beta
}-1).\noindent $ Then:

\bigskip

\textbf{Theorem }\textit{Suppose that the weight }$(\vec{k},w)$\textit{\ is }%
$p$\textit{-small and generic, i.e., }$2<k_{i}^{(j)}\leq p+1$\textit{\ for
all }$i,j$\textit{. Then, if }$\Omega $\textit{\ is a Hecke eigensystem
occurring in the space }$S_{\tau ,\psi }(U,\mathcal{O})$\textit{, there is a
finite local extension of discrete valuation rings }$\mathcal{O}^{\prime }/%
\mathcal{O}$\textit{\ and an }$\mathcal{O}^{\prime }$\textit{-valued Hecke
eigensystem }$\Omega ^{\prime }$\textit{\ occurring in holomorphic weight }$(%
\vec{k}+\vec{a},w^{\prime })$\textit{\ and with associated Hecke character }$%
\psi ^{\prime }$\textit{\ such that }$\Omega ^{\prime }(\func{mod}\mathfrak{M%
}_{\mathcal{O}^{\prime }})=\Omega (\func{mod}\mathfrak{M}_{\mathcal{O}}).$ 
\textit{The character }$\psi ^{\prime }$\textit{\ is compatible with the
weight }$(\vec{k}+\vec{a},w^{\prime })$\textit{\ and it can be chosen so
that }$\bar{\psi}^{\prime }=\bar{\psi}$\textit{.}

\bigskip

More weight shiftings results are proved in \ref{more3} under the assumption 
$f_{j}<3$ for all $j$. The combinatoric involved in describing all the
holomorphic weight shiftings arising from the generalized Dickson and $D$
operators becomes more complicated as the $%
%TCIMACRO{\U{2124} }%
%BeginExpansion
\mathbb{Z}
%EndExpansion
_{p}$-rank of $\mathcal{O}$ grows.

The techniques of sections 3 and 4 cannot be successfully applied to obtain
weight shiftings by $p-1$ when starting from weights that contain at least a 
$(2,...,2)$-block (for example parallel weight two). In section 5 we
therefore re-present a result due to Edixhoven and Khare (\cite{EK}):

\bigskip

\textbf{Theorem }\textit{Assume that }$\tau $\textit{\ is an irreducible
(non necessarily holomorphic) }$\mathbb{F}$\textit{-linear weight with
parameters }$(\vec{k},\vec{w})\in 
%TCIMACRO{\U{2124} }%
%BeginExpansion
\mathbb{Z}
%EndExpansion
_{\geq 2}^{g}\times 
%TCIMACRO{\U{2124} }%
%BeginExpansion
\mathbb{Z}
%EndExpansion
^{g}$\textit{\ such that }$\vec{k}^{(j)}=\vec{2}$\textit{\ for some }$1\leq
j\leq r$\textit{. Let }$\tau ^{\prime }$\textit{\ be the }$\mathbb{F}$%
\textit{-linear weight associated to the parameters }$\vec{k}^{\prime }=(%
\vec{k}^{(1)},...,\vec{k}^{(j)}+\overrightarrow{p-1},...,\vec{k}^{(r)})$%
\textit{\ and }$\vec{w}^{\prime }=\vec{w}$\textit{. For any non-Eisenstein
maximal ideal }$\mathfrak{M}$\textit{\ of }$T_{S,\mathbb{F}}^{univ}$\textit{%
, there is an injective Hecke-equivariant }$\mathbb{F}$\textit{-morphism:}%
\begin{equation*}
S_{\tau }(U,\mathbb{F})_{\mathfrak{M}}\hookrightarrow S_{\tau ^{\prime }}(U,%
\mathbb{F})_{\mathfrak{M}}.
\end{equation*}

\bigskip

\noindent The proof of this result of Edixhoven and Khare relies on the
determination of the $\mathbb{T}_{S,\mathbb{F}}^{univ}$-support of the
kernel of a degeneracy map $S_{\tau ,\psi }(U,\mathbb{F})^{2}\rightarrow
S_{\tau ,\psi }(U_{0},\mathbb{F})$. We remark that in the above theorem the
weight $\tau $ is not assumed parallel. The weight shiftings produced by
repeatedly applying this theorem are not parallel, but parallel in blocks.
We do not know if, starting from weight two, weight shiftings by $p-1$ which
are not of this type are possible or if they can be obtained via the above
methods.

\bigskip

\textit{Conventions} Unless otherwise stated, in this paper all rings are
assumed to have an identity element and are commutative. All the group
representations are assumed to be left representations on a module of finite
length over a fixed coefficient ring. The letter $p$ always denotes a
positive rational prime.

\bigskip

\textbf{Acknowledgements}

The problem of constructing cohomological weight shiftings in the context of
the Hasse invariant was conceived and suggested to me by Chandrashekhar
Khare, to whom I am deeply indebted. His viewpoint on the problem is the
core of \cite{Re}, in which the elliptic case is considered, and it is
therefore also the base of our study in the present paper. The idea of the
existence of a possible connection between cohomological weight shifting
operators and integral models of some irreducible characteristic zero
representations of $GL_{2}(\mathbb{F}_{q})$, is entirely his. \noindent I
would like to further thank C. Khare for his constant, patient and generous
advices during the years of my Ph.D. in UCLA, when this paper was written.

I would like to thank Jean-Pierre Serre, for defining the differential
operator $D$ acting on the polynomial algebra $\mathbb{F}_{q}[X,Y]$, which
plays a crucial role in \cite{Re} and in the present paper.

I am very grateful to Fred Diamond, who suggested to look into
generalizations of the results of \cite{Re}, and in particular into the
existence of cohomological weight shiftings "by $(p,-1,0,...,0)$". This
induced me to construct and study the operators $D_{\beta }^{[\alpha ]}$
appearing in this paper.

I would like to express my gratitude to Claus Sorensen, for sharing with me
his unpublished results on cohomological weight shiftings for modular forms
having trivial weight. I mostly do not address this case in the present
paper, and the results of Edixhoven-Khare (\cite{EK}) and Sorensen complete
in this sense the picture of weight shiftings considered here.

I\ would like to thank Don Blasius, Haruzo Hida, Gordan Savin, and Jacques
Tilouine for their precious comments and questions on the topics studied in
this paper.

\newpage

\part{Weight shiftings for $GL_{2}(\mathbb{F}_{q})$-modules}

\section{Untwisted $GL_{2}(\mathbb{F}_{q})$-modules\label{sec2}}

Fix a rational prime $p$, a positive integer $g$, and set $q=p^{g}$. Denote
by $\mathbb{F}_{q}$ a finite field with $q$ elements and fix an algebraic
closure $\overline{\mathbb{F}}_{q}$ of $\mathbb{F}_{q}$; denote by $\sigma
\in \limfunc{Gal}\left( \mathbb{F}_{q}/\mathbb{F}_{p}\right) $ the
arithmetic Frobenius element. Let $G=GL_{2}\left( \mathbb{F}_{q}\right) $
and let $M$ be a representation of $G$ over $\mathbb{F}_{q}$; for any $n\in 
%TCIMACRO{\U{2124} }%
%BeginExpansion
\mathbb{Z}
%EndExpansion
$, the Frobenius element $\sigma ^{n}$ induces a map $G\rightarrow G$
obtained by applying $\sigma ^{n}$ to each entry of the matrices in $G$:
composing this map with the action of $G$ on $M$,\ we give to the latter a
new structure of $G$-module, that is denoted $M^{[n]}$ and called the $n$th
Frobenius twist of $M.$ If $f:M\rightarrow N$ is a $G$-homomorphism and $%
n\in 
%TCIMACRO{\U{2124} }%
%BeginExpansion
\mathbb{Z}
%EndExpansion
$, denote by $f^{[n]}:M^{[n]}\mathbb{\rightarrow }N^{[n]}$ the map defined
by $f^{[n]}(x)=f(x)$ for all $x\in M^{[n]}$: $f^{[n]}$ is a $G$-homomorphism.

Let $M_{1}$ denote the standard representation of $G$ on $\mathbb{F}_{q}^{2}$
and, for any positive integer $k$, define $M_{k}=\limfunc{Sym}^{k}M_{1}$ to
be the $k$th symmetric power of $M_{1}$. We identify $M_{k}$ with the $%
\mathbb{F}_{q}$-vector space of homogeneous polynomials over $\mathbb{F}_{q}$
in two variables and of degree $k$, endowed with the action of $G$ induced
by:

\begin{center}
\begin{equation*}
\left( 
\begin{array}{cc}
a & b \\ 
c & d%
\end{array}%
\right) \cdot X=aX+cY,\ \left( 
\begin{array}{cc}
a & b \\ 
c & d%
\end{array}%
\right) \cdot Y=bX+dY.
\end{equation*}
\end{center}

\noindent We set $M_{0}$ to be the trivial representation of $G$.\ \noindent
Denote by $\det :G\rightarrow \mathbb{F}_{q}^{\times }$ the determinant
character of $G$, so that $\det^{[n]}=\det^{p^{n}}$.

Recall (cf. \cite{St}, \cite{St2} \S 13) that the irreducible
representations of $G$ over $\mathbb{F}_{q}$\ are all and only of the form:%
\begin{equation*}
\det\nolimits^{m}\otimes _{\mathbb{F}_{q}}\dbigotimes%
\nolimits_{i=0}^{g-1}M_{k_{i}}^{[i]},
\end{equation*}

\noindent where $k_{0},...,k_{g-1}$ and $m$ are integers such that $0\leq
k_{i}\leq p-1$ for $i=0,...,g-1$, $0\leq m<q-1$, and all the tensor products
are over $\mathbb{F}_{q}.$ The above representations are pairwise
non-isomorphic.

We denote by $K_{0}(G)$ the Grothendieck group of finitely generated $%
\mathbb{F}_{q}[G]$-modules: it can be identified with the free abelian group
generated by the isomorphism classes of irreducible representations of $G$
over $\mathbb{F}_{q}$ (\cite{Serre}). If $M$ is an $\mathbb{F}_{q}[G]$%
-module, we denote by $[M]$ its class in $K_{0}(G)$ and set $e=[\det ]$; if
no confusion arises we also write $M$ to denote $[M]$. Tensor product over $%
\mathbb{F}_{q}$ induces on $K_{0}(G)$ a structure of commutative ring with
identity; we denote the product in $K_{0}(G)$ by $\cdot $ or by
juxtaposition.

\subsection{Identities in $K_{0}(G)$ (I)}

We present some identities between virtual representations in $K_{0}(G)$
that we will need later.

\paragraph{Negative weights}

We extend the definition of $M_{k}\in K_{0}\left( G\right) $ for $k<0$ in a
way that is coherent with Brauer character computations, as suggested by
Serre in \cite{Se01}. We briefly explain this: a more detailed account of
what follows is contained in \cite{Re} 2.1.

Let $\mathbf{G}=GL_{2}$ as an algebraic group over $\mathbb{F}_{q}$, and let 
$\mathbf{T\subset G\ }$be the maximal split torus of diagonal matrices.
Identify the character group $X(\mathbf{T})$ of $\mathbf{T}$ with $%
%TCIMACRO{\U{2124} }%
%BeginExpansion
\mathbb{Z}
%EndExpansion
^{2}$ in the usual way, so that the roots associated to $\left( \mathbf{G},%
\mathbf{T}\right) $ are $(1,-1)$ and $(-1,1)$;\ fix a choice of positive
root $\alpha =(1,-1)$. The corresponding Borel subgroup $\mathbf{B}$ is the
group of upper triangular matrices in $\mathbf{G}$; we denote by $\mathbf{B}%
^{-}$ the opposite Borel subgroup. For a fixed $\lambda \in $ $X(\mathbf{T})$%
, let $\mathbf{M}_{\lambda }$ be the one dimensional left $\mathbf{B}^{-}$%
-module on which $\mathbf{B}^{-}$ acts (through $\mathbf{T}$) via the
character $\lambda $ . Denote by $\limfunc{ind}\nolimits_{\mathbf{B}^{-}}^{%
\mathbf{G}}\mathbf{M}_{\lambda }$ the left $\mathbf{G}$-module given by
algebraic induction from $\mathbf{B}^{-}$ to $\mathbf{G}$ of $\mathbf{M}%
_{\lambda }.$ Define the following generalization of the dual Weyl module
for $\lambda $ (cf. \cite{Jan}, II.5):%
\begin{equation*}
W\left( \lambda \right) =\sum\nolimits_{i\geq 0}\left( -1\right) ^{i}\cdot
R^{i}\limfunc{ind}\nolimits_{\mathbf{B}^{-}}^{\mathbf{G}}\left( \mathbf{M}%
_{\lambda }\right) \text{,}
\end{equation*}

\noindent where $R^{i}\limfunc{ind}\nolimits_{\mathbf{B}^{-}}^{\mathbf{G}%
}\left( \cdot \right) $ denote the $i$th right derived functor of $\limfunc{%
ind}\nolimits_{\mathbf{B}^{-}}^{\mathbf{G}}\left( \cdot \right) $. $W\left(
\lambda \right) $ is an element of the Grothendieck group $K_{0}(\mathbf{G})$
of $\mathbf{G}$, because each $R^{i}\limfunc{ind}\nolimits_{\mathbf{B}^{-}}^{%
\mathbf{G}}\left( \mathbf{M}_{\lambda }\right) $ is a finite dimensional $%
\mathbf{G}$-module, and $R^{i}\limfunc{ind}\nolimits_{\mathbf{B}^{-}}^{%
\mathbf{G}}\left( \mathbf{M}_{\lambda }\right) $ is zero for $i>1$ (\cite%
{Jan}, II.4.2). For $\lambda _{k}=(k,0)\in X(\mathbf{T})$ with $k$ any
integer we have: 
\begin{equation*}
R^{i}\limfunc{ind}\nolimits_{\mathbf{B}^{-}}^{\mathbf{G}}\left( \mathbf{M}%
_{\lambda _{k}}\right) \simeq H^{i}(\mathbb{P}_{\mathbb{F}_{q}}^{1},\mathcal{%
O}\left( k\right) ).
\end{equation*}

If $k\geq 0$, $H^{1}(\mathbb{P}_{\mathbb{F}_{q}}^{1},\mathcal{O}\left(
k\right) )=0$ so that $W\left( \lambda _{k}\right) =H^{0}(\mathbb{P}_{%
\mathbb{F}_{q}}^{1},\mathcal{O}\left( k\right) )=\limfunc{Sym}\nolimits^{k}%
\mathbb{F}_{q}^{2}$; if $k<0$ we have $H^{0}(\mathbb{P}_{\mathbb{F}_{q}}^{1},%
\mathcal{O}\left( k\right) )=0$\ and $W\left( \lambda _{k}\right) =-H^{1}(%
\mathbb{P}_{\mathbb{F}_{q}}^{1},\mathcal{O}\left( k\right) )$; the canonical
perfect pairing of $\mathbf{G}$-modules: 
\begin{equation*}
H^{0}(\mathbb{P}_{\mathbb{F}_{q}}^{1},\mathcal{O}\left( -k-2\right) )\times
H^{1}(\mathbb{P}_{\mathbb{F}_{q}}^{1},\mathcal{O}\left( k\right)
)\rightarrow H^{1}(\mathbb{P}_{\mathbb{F}_{q}}^{1},\mathcal{O}\left(
-2\right) )\simeq \det\nolimits^{-1}\otimes \mathbb{G}_{a,\mathbb{F}_{q}},
\end{equation*}%
brings naturally to the following:

\begin{definition}
Let $k<0$ be an integer. Define the element $M_{k}$ of the Grothendieck
group $K_{0}\left( G\right) $ of $G$ over $\mathbb{F}_{q}$ by:%
\begin{equation*}
M_{k}=\left\{ 
\begin{array}{cc}
0 & \text{if }k=-1 \\ 
-e^{1+k}\cdot M_{-k-2} & \text{if }k\leq -2%
\end{array}%
.\right.
\end{equation*}
\end{definition}

\begin{lemma}
For any $k\in 
%TCIMACRO{\U{2124} }%
%BeginExpansion
\mathbb{Z}
%EndExpansion
$\ we have in $K_{0}\left( G\right) $ the identity:%
\begin{equation}
M_{k}+e^{1+k}\cdot M_{-k-2}=0.  \tag{$\Delta _{g,k}$}
\end{equation}
\end{lemma}

\paragraph{Weight shifting by $q\pm 1$}

Let us fix an embedding $\iota :\mathbb{F}_{q^{2}}\rightarrow M_{2}(\mathbb{F%
}_{q})$ corresponding to a choice of $\mathbb{F}_{q}$-basis for the degree $%
2 $ extension of $\mathbb{F}_{q}$ inside $\overline{\mathbb{F}}_{q}$. Let $%
\overline{%
%TCIMACRO{\U{211a} }%
%BeginExpansion
\mathbb{Q}
%EndExpansion
}_{p}$\ be a fixed algebraic closure of the $p$-adic field $%
%TCIMACRO{\U{211a} }%
%BeginExpansion
\mathbb{Q}
%EndExpansion
_{p}$ and let us fix an isomorphism between $\overline{\mathbb{F}}_{q}$ and
the residue field of the ring of integers $\overline{%
%TCIMACRO{\U{2124} }%
%BeginExpansion
\mathbb{Z}
%EndExpansion
}_{p}$ of $\overline{%
%TCIMACRO{\U{211a} }%
%BeginExpansion
\mathbb{Q}
%EndExpansion
}_{p}$; denoting by $\chi :\overline{\mathbb{F}}_{q}^{\times }\mathbb{%
\rightarrow }\overline{%
%TCIMACRO{\U{2124} }%
%BeginExpansion
\mathbb{Z}
%EndExpansion
}_{p}^{\times }$ the corresponding Teichm\"{u}ller character, the Brauer
character $G_{reg}\rightarrow \overline{%
%TCIMACRO{\U{211a} }%
%BeginExpansion
\mathbb{Q}
%EndExpansion
}_{p}$ of the representations $M_{k}$ ($k\geq 1$) is given as follows:%
\begin{eqnarray*}
\left( 
\begin{array}{cc}
a &  \\ 
& a%
\end{array}%
\right) &\mapsto &(k+1)\chi \left( a\right) ^{k},\text{ \ \ \ \ \ \ \ \ \ \
\ \ \ \ \ }a\in \mathbb{F}_{q}^{\times } \\
\left( 
\begin{array}{cc}
a &  \\ 
& b%
\end{array}%
\right) &\mapsto &\frac{\chi \left( a\right) ^{k+1}-\chi \left( b\right)
^{k+1}}{\chi \left( a\right) -\chi \left( b\right) },\text{ \ \ \ \ \ }%
a,b\in \mathbb{F}_{q}^{\times },a\neq b \\
\iota \left( c\right) &\mapsto &\frac{\chi \left( c\right) ^{q\left(
k+1\right) }-\chi \left( c\right) ^{k+1}}{\chi \left( c\right) ^{q}-\chi
\left( c\right) },\text{ \ \ \ }c\in \mathbb{F}_{q^{2}}^{\times }\backslash 
\mathbb{F}_{q}^{\times }.
\end{eqnarray*}

Using the above formulae, the following is proved in \cite{Se01}:

\begin{lemma}
For any $k\in 
%TCIMACRO{\U{2124} }%
%BeginExpansion
\mathbb{Z}
%EndExpansion
$\ we have in $K_{0}\left( G\right) $ the identity:%
\begin{equation}
M_{k}-e\cdot M_{k-(q+1)}=M_{k-\left( q-1\right) }-e\cdot M_{k-2q}. 
\tag{$\Sigma _{g,k}$}
\end{equation}
\end{lemma}

\paragraph{Product formula}

It is a result of Glover that for any positive integers $n,m$ there exists a
short exact sequence of $\overline{\mathbb{F}}_{q}[SL_{2}(\overline{\mathbb{F%
}}_{q})]$-modules of the form:%
\begin{equation*}
0\rightarrow M_{n-1}\otimes _{\overline{\mathbb{F}}_{q}}M_{m-1}\overset{j}{%
\rightarrow }M_{n}\otimes _{\overline{\mathbb{F}}_{q}}M_{m}\overset{\pi }{%
\rightarrow }M_{n+m}\rightarrow 0,
\end{equation*}

\noindent where $j$ is induced by the assignment $u\otimes v\mapsto
uX\otimes vY-uY\otimes vX$ and $\pi $ is induced by multiplication inside
the algebra $\overline{\mathbb{F}}_{q}[X,Y].$ \noindent The following is an
easy extension to $GL_{2}$ of Glover's result:

\begin{lemma}
For any $n,m\in 
%TCIMACRO{\U{2124} }%
%BeginExpansion
\mathbb{Z}
%EndExpansion
$\ we have in $K_{0}\left( G\right) $ the identity:%
\begin{equation}
M_{n}M_{m}=M_{n+m}+eM_{n-1}M_{m-1}.  \tag{$\Pi _{g,n,m}$}
\end{equation}
\end{lemma}

\textbf{Proof }Let $\tau $ be the Brauer character of the virtual
representation $M_{n}M_{m}-M_{n+m}-eM_{n-1}M_{m-1}$. Let $a,b\in \mathbb{F}%
_{q}^{\times }$ such that $a\neq b$; denote by $\tilde{x}$ the Teichm\"{u}%
ller lift of $x\in \overline{\mathbb{F}}_{q}^{\times }$ taken via $\chi $.
We have:

\begin{eqnarray*}
\tau \left( 
\begin{array}{cc}
a &  \\ 
& a%
\end{array}%
\right) &=&(n+1)(m+1)\tilde{a}^{n+m}-(n+m+1)\tilde{a}^{n+m}+ \\
&&-\tilde{a}^{2}\cdot nm\tilde{a}^{(n-1)+(m-1)}; \\
\tau \left( 
\begin{array}{cc}
a &  \\ 
& b%
\end{array}%
\right) &=&\frac{\tilde{a}^{n+1}-\tilde{b}^{n+1}}{\tilde{a}-\tilde{b}}\frac{%
\tilde{a}^{m+1}-\tilde{b}^{m+1}}{\tilde{a}-\tilde{b}}-\frac{\tilde{a}%
^{n+m+1}-\tilde{b}^{n+m+1}}{\tilde{a}-\tilde{b}}+ \\
&&-\tilde{a}\tilde{b}\frac{(\tilde{a}^{n}-\tilde{b}^{n})(\tilde{a}^{m}-%
\tilde{b}^{m})}{(\tilde{a}-\tilde{b})^{2}}.
\end{eqnarray*}

\noindent Both these expressions are trivially zero. If $c\in \mathbb{F}%
_{q^{2}}^{\times }\backslash \mathbb{F}_{q}^{\times }$ and $\iota :\mathbb{F}%
_{q^{2}}{\rightarrow }M_{2}\left( \mathbb{F}_{q}\right) $ is as above, then $%
\det \iota \left( c\right) =c^{1+q}$, so that:

\begin{eqnarray*}
\tau (\iota \left( c\right) ) &=&\frac{\tilde{c}^{q(n+1)}-\tilde{c}^{n+1}}{%
\tilde{c}^{q}-\tilde{c}}\frac{\tilde{c}^{q(m+1)}-\tilde{c}^{m+1}}{\tilde{c}%
^{q}-\tilde{c}}+ \\
&&-\frac{\tilde{c}^{q(n+m+1)}-\tilde{c}^{n+m+1}}{\tilde{c}^{q}-\tilde{c}}-%
\tilde{c}^{1+q}\frac{(\tilde{c}^{qn}-\tilde{c}^{n})(\tilde{c}^{qm}-\tilde{c}%
^{m})}{(\tilde{c}^{q}-\tilde{c})^{2}},
\end{eqnarray*}

\noindent and this is also zero. As $\tau $ is identically zero on $G^{reg}$%
, $M_{n}M_{m}-M_{n+m}-eM_{n-1}M_{m-1}$ is the zero element of $K_{0}\left(
G\right) $. $\blacksquare $

\bigskip

We summarize the three identities obtained so far:

\begin{proposition}
\label{3Id}Let $q=p^{g}$ ($g\geq 1$) and let $k,n,m\in 
%TCIMACRO{\U{2124} }%
%BeginExpansion
\mathbb{Z}
%EndExpansion
$. The following identities hold in $K_{0}\left( G\right) $:%
\begin{equation}
M_{k}=-e^{1+k}\cdot M_{-k-2}  \tag{$\Delta _{g,k}$}
\end{equation}%
\begin{equation}
M_{k}-e\cdot M_{k-(q+1)}=M_{k-\left( q-1\right) }-e\cdot M_{k-2q} 
\tag{$\Sigma _{g,k}$}
\end{equation}%
\begin{equation}
M_{n}M_{m}=M_{n+m}+eM_{n-1}M_{m-1}.  \tag{$\Pi _{g,k}$}
\end{equation}
\end{proposition}

\subsection{Intertwining operators for the periods $q-1$ and $q+1\label%
{g=1_op}$}

Recall that the irreducible complex representations of $G$ (of dimension
larger than one) that are not twists of the Steinberg representation are of
two types:\ the principal series representations, having dimension $q+1$ and
obtained by inducing to $G$ characters of the Borel subgroup of $G$, and the
cuspidal representations, having dimension $q-1$ and characterized by the
property that they do not occur as a factor of a principal series.

The two periods\ $q+1$ and $q-1$ appear in the identity $\left( \Sigma
_{g,k}\right) $ and suggest the existence of intertwining operators that
shift weights by $q+1$ and $q-1$ respectively; furthermore one expects these
operators to give a bridge between the modular representations of $G$ and
the above mentioned characteristic zero representations of $G$. We recall
below the known results on this matter, as the operators we introduce here
will be the starting point of the generalizations considered in the
following sections (cf. \ref{intertnow}).

\paragraph{The period $q+1\label{Dickson}$}

Let $k>q$ be an integer and let $\Theta _{q}=XY^{q}-X^{q}Y\in \mathbb{F}_{q}%
\left[ X,Y\right] $. (Dickson proved that this polynomial is one of the two
generators of the ring of $SL_{2}\left( \mathbb{F}_{q}\right) $-invariants
in the symmetric algebra $\limfunc{Sym}^{\ast }\mathbb{F}_{q}^{2}$, so we
will call it the Dickson invariant). Let us denote by $\Theta _{q}$ also the 
$G$-equivariant map $\det \otimes M_{k-\left( q+1\right) }\mathbb{%
\rightarrow }M_{k}$ given by multiplication by $\Theta _{q}$.

\begin{proposition}
\label{Dickson prop}For $k>q$, there is an exact sequence of $G$-modules:%
\begin{equation*}
0\rightarrow \det \otimes M_{k-\left( q+1\right) }\overset{\Theta _{q}}{%
\rightarrow }M_{k}\rightarrow \limfunc{Ind}\nolimits_{B}^{G}\left( \eta
^{k}\right) \rightarrow 0,
\end{equation*}

\noindent where $B$ is the subgroup of $G$ consisting of upper triangular
matrices, and $\eta $ is the character of $B$ defined extending the
character $\limfunc{diag}(a,b)\mapsto a$ of the standard maximal torus of $G$%
. Furthermore, for any integer $\lambda \geq 0$\ there are isomorphisms of $%
G $-modules:%
\begin{equation*}
\frac{M_{k}}{\det \otimes M_{k-(q+1)}}\simeq \frac{M_{k+\lambda (q-1)}}{\det
\otimes M_{k+\lambda (q-1)-(q+1)}},
\end{equation*}%
\noindent where the inclusion $\det \otimes M_{k+\lambda
(q-1)-(q+1)}\hookrightarrow M_{k+\lambda (q-1)}$ is induced by the
multiplication by $\Theta _{q}$.
\end{proposition}

\textbf{Proof }The above result is standard; cf. \cite{Re}, Proposition 2.7. 
$\blacksquare $

\paragraph{The period $q-1\label{D-map}$}

The period $q-1$ is studied in \cite{Re}: here we just recall the main
result proved there. The starting point is the $G$-equivariant derivation $D:%
\mathbb{F}_{q}[X,Y]\mathbb{\rightarrow }\mathbb{F}_{q}[X,Y]$ defined by
Serre as:%
\begin{equation*}
D:\text{ \ }f(X,Y)\longmapsto X^{q}\frac{\partial f}{\partial X}(X,Y)+Y^{q}%
\frac{\partial f}{\partial Y}(X,Y).
\end{equation*}

This map defines by restriction an intertwining operator $M_{k}\mathbb{%
\rightarrow }M_{k+(q-1)}$ for any $k\geq 0$, giving rise to a weight
shifting by $q-1$. The kernel of $D$ is often non trivial (\cite{Re},
Proposition 3.3), and $D$\ captures essential properties related to the
existence or non-existence of embeddings of $G$-modules of the form $%
M_{k}\rightarrow M_{k+(q-1)}$ (\cite{Re}, Proposition 3.5 and Proposition
3.6).

We now assume, for the rest of this paragraph, that $p$ is an odd prime. If
we restricted ourselves to weights $2\leq k\leq p-1$ we have the following
exact sequence:%
\begin{equation*}
0\rightarrow \det \otimes M_{k-2}\overset{\overline{\Theta }_{q}}{%
\rightarrow }\dfrac{M_{k+(q-1)}}{D(M_{k})}\rightarrow \limfunc{coker}%
\overline{\Theta }_{q}\rightarrow 0,
\end{equation*}

\noindent where $\overline{\Theta }_{q}=\Theta _{q}(\func{mod}D(M_{k}))$ is
induced by the Dickson invariant.

\begin{theorem}
\label{cris}\textit{Let }$q\neq 2$\textit{, }$2\leq k\leq p-1$\textit{\ with 
}$k\neq \frac{q+1}{2}$\textit{\ and let us denote by }$\Xi \left( \chi
^{k}\right) $\textit{\ the cuspidal }$\overline{%
%TCIMACRO{\U{211a} }%
%BeginExpansion
\mathbb{Q}
%EndExpansion
}_{p}$\textit{-representation of }$G$\textit{\ associated to the }$k$th%
\textit{-power of the Teichm\"{u}ller character }$\chi $\textit{. Let }$%
\mathcal{C}$ be the Deligne-Lusztig variety of $SL_{2/\mathbb{F}_{q}}$. 
\textit{There exists a canonical }$W(\mathbb{F}_{q})$-\textit{integral model 
}%
\begin{equation*}
\tilde{\Xi}\left( \chi ^{k}\right) :=H_{\limfunc{cris}}^{1}(\mathcal{C}_{/%
\mathbb{F}_{q}})_{-k}
\end{equation*}%
\textit{\ of }$\Xi \left( \chi ^{k}\right) $,\textit{\ arising from the }$%
\left( -k\right) $\textit{-eigenspace of the first crystalline cohomology
group of }$\mathcal{C}_{/\mathbb{F}_{q}}$,\textit{\ such that there is an
isomorphism of }$\mathbb{F}_{q}\left[ G\right] $\textit{-modules:}%
\begin{equation*}
\dfrac{M_{k+(q-1)}}{D\left( M_{k}\right) }\simeq \tilde{\Xi}\left( \chi
^{k}\right) \text{ }\otimes _{W(\mathbb{F}_{q})}\mathbb{F}_{q}.
\end{equation*}

\noindent \noindent (Here the $(-k)$-eigenspace of $H_{\limfunc{cris}}^{1}(%
\mathcal{C}_{/\mathbb{F}_{q}})$ is computed with respect to the natural
action of $\ker (\limfunc{Nm}_{\mathbb{F}_{q^{2}}^{\times }/\mathbb{F}%
_{q}^{\times }})$ on $H_{\limfunc{cris}}^{1}(\mathcal{C}_{/\mathbb{F}_{q}})$%
).
\end{theorem}

\textbf{Proof }\cite{Re}, Theorem 4.2. $\blacksquare $

\subsection{Determination of Jordan-H\"{o}lder constituents: the case $g=1$}

Assume $g$ is any positive integer. For convenience, we give the following
non standard definition:

\begin{definition}
Let $M\in K_{0}(G)$ be of the form $M=e^{m}\tprod%
\nolimits_{i=0}^{g-1}M_{k_{i}}^{[i]}$ where $m,k_{0},...,k_{g-1}\in 
%TCIMACRO{\U{2124} }%
%BeginExpansion
\mathbb{Z}
%EndExpansion
$. We say that the Jordan-H\"{o}lder factors of $M$\ can be computed in the
standard form (using $(\Delta _{g}),(\Sigma _{g})$\ and $(\Pi _{g})$)\ if,
by applying finitely many times the identities of Proposition \ref{3Id},
together with the identities $e^{q-1}=1$ and $\sigma ^{g}=1$,\ we can write $%
M$ as:%
\begin{equation*}
M=\sum\nolimits_{j\in J}n_{j}\left(
e^{m_{j}}\prod\nolimits_{i=0}^{g-1}M_{k_{i}^{(j)}}^{[i]}\right) ,
\end{equation*}

where $J$ is a finite set and for any $j\in J$ we have $%
n_{j},m_{j},k_{0}^{(j)},...,k_{g-1}^{(j)}\in 
%TCIMACRO{\U{2124} }%
%BeginExpansion
\mathbb{Z}
%EndExpansion
$ such that $n_{j}\neq 0$, $0\leq m_{j}<q-1,$ $0\leq
k_{0}^{(j)},...,k_{g-1}^{(j)}\leq p-1$ and, if $j,j^{\prime }\in J$ with $%
j\neq j^{\prime }$ then $(m_{j},k_{0}^{(j)},...,k_{g-1}^{(j)})\neq
(m_{j^{\prime }},k_{0}^{(j^{\prime })},...,k_{g-1}^{(j^{\prime })})$.
(Notice that the integers $n_{j},m_{j},k_{0}^{(j)},...,k_{g-1}^{(j)}$ are
uniquely determined by $M$).

Similarly one defines the notion of computability in standard form for an
element of $K_{0}(G)$ that is given as an algebraic sum of products of
elements of the form $M=e^{m}\tprod\nolimits_{i=0}^{g-1}M_{k_{i}}^{[i]}$.
Also, in an obvious way, one defines computability in standard form using
any subset or superset of the identities $(\Delta ),(\Sigma )$ and $(\Pi )$
(together with the identities $e^{q-1}=1$ and $\sigma ^{g}=1$).
\end{definition}

\begin{lemma}
\label{product}Let $g$ be any positive integer and let $n,m\in 
%TCIMACRO{\U{2124} }%
%BeginExpansion
\mathbb{Z}
%EndExpansion
$ such that $n,m\geq 0$. By applying $(\Pi _{g})$ we obtain the following
identity in $K_{0}(G):$%
\begin{equation*}
M_{n}M_{m}=\dsum\nolimits_{i=0}^{\min \{n,m\}}e^{i}M_{n+m-2i}.
\end{equation*}
\end{lemma}

\textbf{Proof }We induct on $n$. For $n=0$ the statement is true; for $n\geq
0$ we have, assuming $m>0$:%
\begin{eqnarray*}
M_{n+1}M_{m} &=&M_{n+m+1}+eM_{n}M_{m-1}= \\
&=&M_{n+m+1}+\dsum\nolimits_{i=0}^{\min \{n,m-1\}}e^{i+1}M_{n+m-1-2i}= \\
&=&M_{n+m+1}+\dsum\nolimits_{i=1}^{\min \{n+1,m\}}e^{i}M_{(n+1)+m-2i}= \\
&=&\dsum\nolimits_{i=0}^{\min \{n+1,m\}}e^{i}M_{(n+1)+m-2i}.\text{ \ \ }%
\blacksquare
\end{eqnarray*}

\bigskip

\begin{corollary}
\label{product-corollary}For any positive integer $t$ and any integers $%
n_{1},...,n_{t}\geq 0$ we have:%
\begin{equation*}
\prod\nolimits_{i=1}^{t}M_{n_{i}}=\sum\nolimits_{\alpha \in A}e^{s_{\alpha
}}M_{r_{\alpha }},
\end{equation*}

\noindent where $A$ is a finite set and $s_{\alpha },r_{\alpha }\geq 0$ for
any $\alpha \in A$.
\end{corollary}

\textbf{Proof }It follows from applying $(\Pi _{g})$ and inducting on $t$. $%
\blacksquare $

\bigskip

The following proposition guarantees that, if $g=1$, $(\Delta _{1})\ $and $%
(\Sigma _{1})$ are enough to compute explicitly the Jordan-H\"{o}lder
factors of any of the modules $M_{k}$ for $k\in 
%TCIMACRO{\U{2124} }%
%BeginExpansion
\mathbb{Z}
%EndExpansion
$.

\begin{proposition}
\label{JHg=1}Let $g=1$.\ For any $m,k\in 
%TCIMACRO{\U{2124} }%
%BeginExpansion
\mathbb{Z}
%EndExpansion
$, we can compute the Jordan-H\"{o}lder factors of $e^{m}M_{k}$ in the
standard form, using\textbf{\ }$(\Delta _{1})\ $and $(\Sigma _{1})$.
Furthermore, by using also $(\Pi _{1})$, we can compute the Jordan-H\"{o}%
lder factors in the standard form for any algebraic sum of products of $%
e^{m}M_{k}$'s.
\end{proposition}

\textbf{Proof }The second assertion in the statement of the proposition
follows from the first one, together with Lemma \ref{product}. To prove the
first assertion, we can assume $m=0$ and, using $(\Delta _{1})$, we also
suppose $k\geq 0$. Write $k=np+r$ where $n$ is a non-negative integer and $r$
is an integer such that $0\leq r\leq p-1$. We induct on $n$.

If $n=0$, there is nothing to prove. Assume $n\geq 1$ is fixed and that we
can compute the Jordan-H\"{o}lder factors of $M_{k}$ in the standard form,
using $(\Delta _{1}),(\Sigma _{1})$ and $(\Pi _{1})$, for any $k$ of the
form $k=n^{\prime }p+r^{\prime }$ where $0\leq n^{\prime }\leq n-1$ and $%
0\leq r^{\prime }\leq p-1$. If $0\leq r\leq p-1$ we have, applying $(\Sigma
_{1})$, that $M_{np+r}=M_{(n-1)p+(r+1)}+e(M_{(n-1)p+(r-1)}-M_{(n-2)p+r})$.
If $r\neq 0,p-1$ we are done by induction assumption. If $r=0$, then $%
M_{np}=M_{(n-1)p+1}+e(M_{(n-2)p+(p-1)}-M_{(n-2)p})$ and we are done. If $%
r=p-1$, just notice that $%
M_{(n-1)p+p}=M_{np}=M_{(n-1)p+1}+e(M_{(n-2)p+(p-1)}-M_{(n-2)p})$. (When $n=1$
one sometimes needs to apply $(\Delta _{1})$ to canonically compute the
constituents of the virtual representations appearing in these identities). $%
\blacksquare $

\subsection{Application to elliptic modular forms}

In this section we summarize some weight shifting results for elliptic
modular forms modulo $p$ in terms of cohomology of groups; the main
references are \cite{AS}, \cite{EK}, \cite{Re}. We assume $p>3$; by a
modular form mod $p$ we mean the reduction modulo $p$ of a form in
characteristic zero - as defined by Serre and Swinnerton-Dyer, unless
otherwise specified.

Let $N\geq 5$ be a positive integer not divisible by $p$ and denote by $%
M_{k}(N,\overline{\mathbb{F}}_{p})$ the $\overline{\mathbb{F}}_{p}$-vector
space of mod $p$ modular forms for the group $\Gamma _{1}(N)$ having weight $%
k\geq 2$ and with coefficients in $\overline{\mathbb{F}}_{p}$; the Hecke
algebra $\mathcal{H}_{N},$ generated over $\overline{\mathbb{F}}_{p}$ by the
operators $T_{l}$ for $p\nmid l$, acts on this space. The $q$-expansion
homomorphism is an injective map $M_{k}(N,\overline{\mathbb{F}}%
_{p})\hookrightarrow \overline{\mathbb{F}}_{p}[[q]]$.

\noindent The theta operator $\Theta :M_{k}(N,\overline{\mathbb{F}}_{p})%
\mathbb{\rightarrow }M_{k+(p+1)}(N,\overline{\mathbb{F}}_{p})$ is defined on 
$q$-expansion by the formula $\Theta
(\tsum\nolimits_{n}a_{n}q^{n})=\tsum\nolimits_{n}na_{n}q^{n}$; it satisfies $%
\Theta T_{l}=lT_{l}\Theta $ for any prime $l\neq p$ ($T_{l}\in \mathcal{H}%
_{N}$). Denote by $E_{p-1}$ the normalized form of the classical
characteristic zero Eisenstein series whose $q$-expansion is given by:%
\begin{equation*}
E_{p-1}=1-2(p-1)/B_{p-1}\tsum\nolimits_{n}\sigma _{p-2}(n)q^{n};
\end{equation*}

\noindent then $E_{p-1}\in M_{p-1}(1,%
%TCIMACRO{\U{2124} }%
%BeginExpansion
\mathbb{Z}
%EndExpansion
_{(p)})$ and $E_{p-1}\equiv 1(\func{mod}p%
%TCIMACRO{\U{2124} }%
%BeginExpansion
\mathbb{Z}
%EndExpansion
_{(p)}[[q]])$, as $2\zeta (2-p)^{-1}\equiv 0(\func{mod}p)$ by the
Clausen-von Staudt theorem. Multiplication by the reduction mod $p$ of $%
E_{p-1}$ gives rise to a Hecke-equivariant map $M_{k}(N,\overline{\mathbb{F}}%
_{p})\mathbb{\rightarrow }M_{k+(p-1)}(N,\overline{\mathbb{F}}_{p})$, that we
refer to as the Hasse invariant.

In view of the Eichler-Shimura isomorphism, the study of Hecke eigensystems
of $\func{mod}p$ modular forms of weight $k\geq 2$ and level $N$ leads to
the study of the eigenvalues of the Hecke algebra $\mathcal{H}_{N}$ acting
on the cohomology group $H^{1}(\Gamma _{1}(N),M_{k-2})$, where $\Gamma
_{1}(N)$ acts on $M_{k-2}$ via its reduction $\func{mod}p$, and the action
of $\mathcal{H}_{N}$ comes from the $G$-action on $M_{k-2}$ and it\ is
defined as in \cite{AS}. The weight shiftings realized on the spaces of
modular forms by the theta operator and the Hasse invariant have
cohomological counterparts. In \cite{AS}, Ash and Stevens identifies a
group-theoretical analogue of the $\Theta $-operator in the
Hecke-equivariant map induced in cohomology by the Dickson invariant (cf. %
\ref{Dickson}):%
\begin{equation*}
\Theta _{p,\ast }:H^{1}(\Gamma _{1}(N),\det \otimes M_{k-2})\longrightarrow
H^{1}(\Gamma _{1}(N),M_{k+p-1}).
\end{equation*}

\noindent Here the twisting by $\det $ on the left hand side is a
manifestation of the fact that the $\Theta $ operator on spaces of modular
forms is twist-Hecke-equivariant.

Edixhoven and Khare identifies in \cite{EK} a cohomological analogue of the
Hasse invariant in the case $k=2$ by studying the degeneracy map $%
H^{1}\left( \Gamma _{1}(N),M_{0}\right) ^{\oplus 2}\rightarrow H^{1}\left(
\Gamma _{1}(N)\cap \Gamma _{0}\left( p\right) ,M_{p-1}\right) $. In \cite{Re}%
, the $D$-derivation defined in \ref{D-map} is used to allow weight shifting
by $p-1$ for $3\leq k\leq p+1$:

\begin{theorem}
Let $\mathfrak{M}$ be a non-Eisenstein maximal ideal of the Hecke algebra $%
\mathcal{H}_{N}$. If $k\geq 0$ and $H^{1}(\Gamma _{1}(N),M_{k})_{\mathfrak{M}%
}\neq 0$, then also $H^{1}(\Gamma _{1}(N),M_{k+(p-1)})_{\mathfrak{M}}\neq 0$%
. If $0\leq k\leq p-1$, there is a Hecke-equivariant embedding $H^{1}(\Gamma
_{1}(N),M_{k})_{\mathfrak{M}}\hookrightarrow H^{1}(\Gamma
_{1}(N),M_{k+(p-1)})_{\mathfrak{M}}$ \noindent \noindent that is induced by
the derivation $D$ if $0<k\leq p-1$, and is the map defined in \cite{EK} if $%
k=0$.
\end{theorem}

\textbf{Proof }The first statement and the second statement for $k\neq 0$
are proved in \cite{Re}, Proposition 5.1; the second statement in the case $%
k=0$ is treated in \cite{EK}, 2. $\blacksquare $

\bigskip

Notice that the above theorem cannot be deduced only by the existence of the
map $D$, as for $k=0$ the virtual representation $M_{p-1}-M_{0}$ is not
positive in $K_{0}\left( G\right) $. A similar situation will occur later on
when we will consider the more general case of Hilbert modular forms (cf.
section \ref{22222}).

\newpage

\section{Twisted $GL_{2}(\mathbb{F}_{q})$-modules and intertwining operators
for $g>1$\label{sec3}}

We keep the notation of the previous section, so that $p$ is a prime number, 
$g$ a positive integer, and $q=p^{g}$; we denote by $\mathbb{F}_{q}$ a
finite field with $q$ elements and we fix an algebraic closure $\overline{%
\mathbb{F}}_{q}$ of $\mathbb{F}_{q}$; we let $\sigma \in \limfunc{Gal}\left( 
\mathbb{F}_{q}/\mathbb{F}_{p}\right) $ be the arithmetic Frobenius element
and $G=GL_{2}\left( \mathbb{F}_{q}\right) $. If $k\in 
%TCIMACRO{\U{2124} }%
%BeginExpansion
\mathbb{Z}
%EndExpansion
$, $M_{k}^{[i]}$ is its $i$th Frobenius twist of the virtual representation $%
M_{k}$, for any integer $i$.

\subsection{Identities in $K_{0}(G)$ (II)}

None of the identities in $K_{0}(G)$ appearing in Proposition \ref{3Id}
contains a Frobenius twist; this implies that, while $(\Delta _{g}),(\Sigma
_{g}),(\Pi _{g})$ are all we need to compute the Jordan-H\"{o}lder factors
of products of virtual representations of the form $M_{k}$ ($k\in 
%TCIMACRO{\U{2124} }%
%BeginExpansion
\mathbb{Z}
%EndExpansion
$) when $g=1$ (Proposition \ref{JHg=1}), these same three families of
identities are not enough to work out such a computation when $g>1$. For
example, when $g>1$, the Jordan-H\"{o}lder factors of $M_{p}$ are $%
\{M_{1}^{[1]},eM_{p-2}\}$ and they cannot be found using $(\Sigma _{g}).$

\begin{proposition}
\label{Frobtwist}Let $g\geq 1$. For any $k\in 
%TCIMACRO{\U{2124} }%
%BeginExpansion
\mathbb{Z}
%EndExpansion
$ we have in $K_{0}(G)$ the identity:%
\begin{equation}
M_{k}=M_{k-p}M_{1}^{[1]}-e^{p}M_{k-2p}.  \tag{$\Phi _{g,k}$}
\end{equation}
\end{proposition}

\textbf{Proof }Fix an embedding $\iota :\mathbb{F}_{q^{2}} \rightarrow
M_{2}\left( \mathbb{F}_{q}\right) $ and denote by $\tilde{x}\in \overline{%
%TCIMACRO{\U{2124} }%
%BeginExpansion
\mathbb{Z}
%EndExpansion
}_{p}$ the Teichm\"{u}ller lift of $x\in \overline{\mathbb{F}}_{q}^{\times }$
taken via the Teichm\"{u}ller character we previously fixed.\textbf{\ }Let $%
\tau $ be the Brauer character of the virtual representation $%
M_{k}-M_{k-p}M_{1}^{[1]}+e^{p}M_{k-2p}$. Let $a,b\in \mathbb{F}_{q}^{\times
} $ such that $a\neq b$. We have:

\begin{eqnarray*}
\tau \left( 
\begin{array}{cc}
a &  \\ 
& a%
\end{array}%
\right) &=&(k+1)\tilde{a}^{k}-(k-p+1)\tilde{a}^{k-p}\cdot 2\tilde{a}^{p}+ \\
&&+\tilde{a}^{2p}\cdot (k-2p+1)\tilde{a}^{k-2p}; \\
\tau \left( 
\begin{array}{cc}
a &  \\ 
& b%
\end{array}%
\right) &=&\frac{\tilde{a}^{k+1}-\tilde{b}^{k+1}}{\tilde{a}-\tilde{b}}-\frac{%
\tilde{a}^{k-p+1}-\tilde{b}^{k-p+1}}{\tilde{a}-\tilde{b}}\frac{\tilde{a}%
^{2p}-\tilde{b}^{2p}}{\tilde{a}^{p}-\tilde{b}^{p}}+ \\
&&+\tilde{a}^{p}\tilde{b}^{p}\frac{\tilde{a}^{k-2p+1}-\tilde{b}^{k-2p+1}}{%
\tilde{a}-\tilde{b}}.
\end{eqnarray*}

\noindent Both these expressions are zero. If $c\in \mathbb{F}%
_{q^{2}}^{\times }\backslash \mathbb{F}_{q}^{\times }$ then $\det \iota
\left( c\right) =c^{1+q}$; also notice that $\limfunc{tr}(\iota \left(
c\right) ;M_{1}^{[1]})=\limfunc{tr}(\iota \left( c\right) ^{\sigma };M_{1})=%
\limfunc{tr}\iota \left( c\right) ^{\sigma }=c^{p}+c^{pq}$, so that:

\begin{eqnarray*}
\tau (\iota \left( c\right) ) &=&\frac{\tilde{c}^{q(k+1)}-\tilde{c}^{k+1}}{%
\tilde{c}^{q}-\tilde{c}}-\frac{\tilde{c}^{q(k-p+1)}-\tilde{c}^{k-p+1}}{%
\tilde{c}^{q}-\tilde{c}}(c^{p}+c^{pq})+ \\
&&+\tilde{c}^{(1+q)p}\frac{\tilde{c}^{q(k-2p+1)}-\tilde{c}^{k-2p+1}}{\tilde{c%
}^{q}-\tilde{c}},
\end{eqnarray*}

\noindent and this is also zero. As $\tau $ is identically zero on $G^{reg}$%
, $M_{k}-M_{k-p}M_{1}^{[1]}+e^{p}M_{k-2p}$ is the zero element of $%
K_{0}\left( G\right) $. $\blacksquare $

\bigskip

\begin{corollary}
\label{Gen. Serre}Let $g\geq 1$; for any $k,h\in 
%TCIMACRO{\U{2124} }%
%BeginExpansion
\mathbb{Z}
%EndExpansion
$ the following identity holds in $K_{0}(G):$%
\begin{equation}
M_{k}M_{h}^{[1]}-e^{p}M_{k-p}M_{h-1}^{[1]}=M_{k-p}M_{h+1}^{[1]}-e^{p}M_{k-2p}M_{h}^{[1]}.
\tag{$\Phi _{g,k,h}^{\prime }$}
\end{equation}
\end{corollary}

\textbf{Proof }Multiplying $(\Phi _{g,k})$ by $M_{h}^{[1]}$ we obtain the
identity 
\begin{equation*}
M_{k}M_{h}^{[1]}=M_{k-p}\left( M_{1}M_{h}\right)
^{[1]}-e^{p}M_{k-2p}M_{h}^{[1]}.
\end{equation*}%
Applying $(\Pi _{g,1,h})$ and distributing the Frobenius action, we deduce
that the left hand side of this equation equals $M_{k-p}\left(
M_{h+1}^{[1]}+e^{p}M_{h-1}^{[1]}\right) -e^{p}M_{k-2p}M_{h}^{[1]}$. $%
\blacksquare $

\bigskip

\begin{corollary}
\label{inttt}Let $g\geq 1$. For any $k,h,i\in 
%TCIMACRO{\U{2124} }%
%BeginExpansion
\mathbb{Z}
%EndExpansion
$ we have in $K_{0}(G)$ the identity:%
\begin{equation*}
M_{k}^{[i]}M_{h}^{[i+1]}-e^{p^{i+1}}M_{k-p}^{[i]}M_{h-1}^{[i+1]}=M_{k-p}^{[i]}M_{h+1}^{[i+1]}-e^{p^{i+1}}M_{k-2p}^{[i]}M_{h}^{[i+1]}.
\end{equation*}
\end{corollary}

\textbf{Proof }Just apply $i$th Frobenius twist to $(\Phi _{g,k,h}^{\prime
}) $. $\blacksquare $

\begin{remark}
\begin{enumerate}
\item \label{comments1}By applying the product formula, one sees that $(\Phi
_{1})$ and $\left( \Sigma _{1}\right) $ are equivalent.

\item Equation $(\Phi _{g}^{\prime })$ ($g>1$) has a structure similar to
equation $\left( \Sigma _{1}\right) $: the weight shiftings appearing in the
latter are by $p+1$ and $p-1$ (corresponding respectively to the degree of
the Dickson invariant and of Serre's derivation map); in equation $(\Phi
_{g}^{\prime })$, the weight shiftings occurring are by $(p,1,0,...,0)$ and $%
(p,-1,0,...,0)$ - the commas separate the shifting constants for tensor
factors corresponding to Frobenius twistings by $\sigma ^{0},\sigma
^{1},...,\sigma ^{g-1}$. In this sense we can think of $(\Phi _{g}^{\prime
}) $ as a generalization of $\left( \Sigma _{1}\right) $ for $g>1$.

\item The reason for which only three (possibly) non-zero terms appear in $%
(\Phi _{g})$ instead of four - as one could have expected by looking at $%
\left( \Sigma _{1}\right) $, is that by applying weight-shifting of $%
(p,1,0,...,0)$ to $M_{k}$ we obtain $e^{p}M_{k-p}M_{-1}^{[1]}$ that is the
zero module: this phenomenon cannot happen when $g=1$.

\item The reason for which, when $g>1$, we were expecting an identity in $%
K_{0}(G)$ involving\ weight shiftings by $(p,\pm 1,0,...,0)$ (and cyclic
permutations of this) resides in the existence of the partial Hasse
invariants and theta operators acting on spaces of mod $p$ Hilbert modular
forms of genus $g$. Also, for good reasons we do not have weight shiftings
by $(\pm 1,p,0,...,0)$, as long as $g>2$: cf. \ref{motivations}.
\end{enumerate}
\end{remark}

\subsection{Determination of Jordan-H\"{o}lder constituents: the case $g>1$}

We know show that equations $(\Delta _{g}),(\Phi _{g}),(\Pi _{g})$ are
enough to compute the Jordan-H\"{o}lder constituents of products of virtual
representations of the form $e^{m}\tprod\nolimits_{i=0}^{g-1}M_{k_{i}}^{[i]}$
($m,k_{0},...,k_{g-1}\in 
%TCIMACRO{\U{2124} }%
%BeginExpansion
\mathbb{Z}
%EndExpansion
$).

\begin{lemma}
Let $g\geq 1$;\ for any $k\in 
%TCIMACRO{\U{2124} }%
%BeginExpansion
\mathbb{Z}
%EndExpansion
$, we can compute the Jordan-H\"{o}lder factors of $M_{k}$ in the standard
form, using\textbf{\ }$(\Delta _{g}),(\Phi _{g}),(\Pi _{g})$.
\end{lemma}

\textbf{Proof }By applying $(\Delta _{g})$ if necessary we can assume $k\geq
0$.\ If $g=1$, the lemma follows from the last remark and Proposition \ref%
{3Id}. For $g\geq 2$, write $k=np+r$ where $n,r\in 
%TCIMACRO{\U{2124} }%
%BeginExpansion
\mathbb{Z}
%EndExpansion
$ are such that $n\geq 0$ and $0\leq r\leq p-1$. We induct on $n$.

If $n=0$, there is nothing to prove. Assume $n\geq 1$ is fixed and that we
can compute the Jordan-H\"{o}lder factors of $M_{k}$ in the standard form,
using $(\Delta _{g}),(\Phi _{g})$ and $(\Pi _{g})$, for any $k$ of the form $%
k=n^{\prime }p+r^{\prime }$ where $0\leq n^{\prime }\leq n-1$ and $0\leq
r^{\prime }\leq p-1$. We have $%
M_{np+r}=M_{(n-1)p+r}M_{1}^{[1]}-e^{p}M_{(n-2)p+r}$ by $(\Phi _{g})$; the
Jordan-H\"{o}lder factors of $e^{p}M_{(n-2)p+r}$ can be computed in the
standard form by induction (if $n=1$ then $M_{(n-2)p+r}=-e^{1-p+r}M_{p-r-2}$
by $(\Delta _{g})$). Also, by induction we have an algorithm that allows us
to write $M_{(n-1)p+r}=\tsum\nolimits_{i\in I}J_{i}$ where $I$ is a finite
set and each $J_{i}$ is of the form $e^{m}\tprod%
\nolimits_{i=0}^{g-1}M_{k_{i}}^{[i]}$ for some integers $m,k_{0},...,k_{g-1}$
such that $0\leq m<q-1,$ $0\leq k_{0},...,k_{g-1}\leq p-1$. It is therefore
enough to show that we can compute the factors of $\left(
\tprod\nolimits_{i=0}^{g-1}M_{k_{i}}^{[i]}\right) M_{1}^{[1]}$ in standard
form, where $0\leq m<q-1,$ $0\leq k_{0},...,k_{g-1}\leq p-1$. The product
formula gives:%
\begin{equation*}
\left( \prod\nolimits_{i=0}^{g-1}M_{k_{i}}^{[i]}\right) M_{1}^{[1]}=\left(
\prod\nolimits_{\substack{ i=0  \\ i\neq 1}}^{g-1}M_{k_{i}}^{[i]}\right)
M_{k_{1}+1}^{[1]}+e^{p}\left( \prod\nolimits_{\substack{ i=0  \\ i\neq 1}}%
^{g-1}M_{k_{i}}^{[i]}\right) M_{k_{1}-1}^{[1]}.
\end{equation*}

If $k_{1}\neq p-1$, each of the two summands is either a Jordan-H\"{o}lder
factor in standard form, or it is zero. Otherwise we are left with the
determination of the constituents of the first summand. If $g=2$ the latter
equals $%
M_{k_{0}}M_{p}^{[1]}=M_{k_{0}}M_{1}+e^{p}M_{k_{0}}M_{p-2}^{[1]}=M_{k_{0}+1}+eM_{k_{0}-1}+e^{p}M_{k_{0}}M_{p-2}^{[1]} 
$ and this is not in standard form if and only if $k_{0}=p-1$, in which case
we can compute the constituents of $M_{k_{0}+1}=M_{p}$ in standard form by
using $(\Phi _{g})$: $M_{p}=M_{1}^{[1]}+eM_{p-2}$.

Assume now $g>2$ and $k_{1}=p-1$. We have, applying $(\Phi _{g})$:%
\begin{equation*}
\left( \prod\nolimits_{\substack{ i=0  \\ i\neq 1}}^{g-1}M_{k_{i}}^{[i]}%
\right) M_{p}^{[1]}=\left( \prod\nolimits_{\substack{ i=0  \\ i\neq 1,2}}%
^{g-1}M_{k_{i}}^{[i]}\right) \left( M_{k_{2}}M_{1}\right) ^{[2]}+e^{p}\left(
\prod\nolimits_{\substack{ i=0  \\ i\neq 1}}^{g-1}M_{k_{i}}^{[i]}\right)
M_{p-2}^{[1]}.
\end{equation*}

\noindent The second summand is already in standard form; for the first
summand we have:%
\begin{equation*}
\left( \prod\nolimits_{\substack{ i=0  \\ i\neq 1,2}}^{g-1}M_{k_{i}}^{[i]}%
\right) \left( M_{k_{2}}M_{1}\right) ^{[2]}=\left( \prod\nolimits_{\substack{
i=0  \\ i\neq 1,2}}^{g-1}M_{k_{i}}^{[i]}\right)
M_{k_{2}+1}^{[2]}+e^{p^{2}}\left( \prod\nolimits_{\substack{ i=0  \\ i\neq
1,2 }}^{g-1}M_{k_{i}}^{[i]}\right) M_{k_{2}-1}^{[2]}.
\end{equation*}

\noindent If $k_{2}\neq p-1$, each of the two summands is either a Jordan-H%
\"{o}lder factor in standard form, or it is zero. Otherwise we are left with
the determination of the constituents of the first summand. We proceed as
before, distinguishing the cases $g=3$ and $g>3$. It is easily seen by
induction that the algorithm produces the Jordan-H\"{o}lder factors of the
virtual representations appearing in each step as long as $k_{i}\neq p-1$
for some $1\leq i\leq g-1$. If $k_{1}=...=k_{g-1}=p-1$, we are left with the
determination the Jordan-H\"{o}lder factors of $%
M_{k_{0}}M_{p}^{[g-1]}=M_{k_{0}}(M_{1}+e^{p^{g-1}}M_{p-2}^{[g-1]})$. By the
product formula, we just need to find the constituents of $M_{k_{0}+1}$: if $%
k_{0}\neq p-1$ this is an irreducible representation; otherwise $%
M_{p}=M_{1}^{[1]}+eM_{p-2}$ and we are done. $\blacksquare $

\bigskip

\begin{corollary}
\label{intz}Let $g\geq 1$. Then $(\Delta _{g}),(\Phi _{g}),(\Pi _{g})$\
imply $(\Sigma _{g})$.
\end{corollary}

\textbf{Proof }By the previous lemma, we can compute the Jordan-H\"{o}lder
factors of each summand appearing in $(\Sigma _{g})$ (in standard form).
Since we know a priori that the Jordan-H\"{o}lder factors appearing in the
right and left hand sides of $(\Sigma _{g})$ have to appear with the same
multiplicities, $(\Sigma _{g})$ is a consequence of $(\Delta _{g}),(\Phi
_{g}),(\Pi _{g})$. $\blacksquare $

\bigskip

Notice that\ we were able to show that $(\Delta _{g}),(\Phi _{g}),(\Pi _{g})$%
\ imply $(\Sigma _{g})$ because we knew already that $(\Sigma _{g})$ was
true. It does not seem to be an easy task to directly deduce Serre's
relation from the set $(\Delta _{g}),(\Phi _{g}),(\Pi _{g})$. Serre's
relation will allow sometimes to bypass long computations involving
Frobenius twists - this will turn out to be useful in \cite{Re3} for $g>1$,
where we will give an explicit presentation of the ring $K_{0}(G)$.

We can finally prove:

\begin{theorem}
\label{final}Let $g\geq 1$. Using $(\Delta _{g}),(\Phi _{g}),(\Pi _{g})$ we
can compute the Jordan-H\"{o}lder factors in the standard form for any
algebraic sum of products of virtual representations of the form $%
e^{m}\tprod\nolimits_{i=0}^{g-1}M_{k_{i}}^{[i]}$ ($m,k_{0},...,k_{g-1}\in 
%TCIMACRO{\U{2124} }%
%BeginExpansion
\mathbb{Z}
%EndExpansion
$).
\end{theorem}

\textbf{Proof }If $g=1$, this is just Proposition \ref{3Id}. Assume $g\geq 2$%
; by applications of $(\Delta _{g})$ and of Lemma \ref{product-corollary},
it is enough to prove that we can compute the Jordan-H\"{o}lder factors in
the standard form for the \textit{representation} $M=\tbigotimes%
\nolimits_{i=0}^{g-1}M_{k_{i}}^{[i]}$ ($k_{0},...,k_{g-1}\geq 0$). We induct
on $\dim _{\mathbb{F}_{q}}M$. If $\dim _{\mathbb{F}_{q}}M=1$, we are done,
otherwise we distinguish two cases.

\underline{\textit{Case 1}}\textit{: There is some }$i$\textit{, }$0\leq
i\leq g-1$,\textit{\ such that }$M_{k_{i}}$\textit{\ is reducible.}

\noindent By applying an appropriate Frobenius twist, we can assume without
loss of generality that $M_{k_{0}}$ is reducible. By the previous lemma, we
can compute the Jordan-H\"{o}lder factors of $M_{k_{0}}$\ in the standard
form, say $M_{k_{0}}=\tsum\nolimits_{h\in I}J_{h}$ in $K_{0}(G)$, where $I$
is a finite set with at least \textit{two} elements and each $J_{h}$ is a
non-zero composition factor of $M_{k_{0}}$, written in standard form. It is
then enough to compute in standard form the constituents of $%
J_{h}\tprod\nolimits_{i=1}^{g-1}M_{k_{i}}^{[i]}$ for each $h\in I$. Fix an
element $h\in I$; up to twisting by a power of $e$ we can assume $%
J_{h}=\tprod\nolimits_{i=0}^{g-1}M_{r_{i}}^{[i]}$ where $0\leq
r_{0},...,r_{g-1}\leq p-1,$ so that an application of Lemma \ref{product}
gives:%
\begin{eqnarray*}
J_{h}\prod\nolimits_{i=1}^{g-1}M_{k_{i}}^{[i]}
&=&M_{r_{0}}\prod\nolimits_{i=1}^{g-1}\left( M_{r_{i}}M_{k_{i}}\right)
^{[i]}= \\
&=&M_{r_{0}}\prod\nolimits_{i=1}^{g-1}\left( \sum\nolimits_{j=0}^{\min
\{r_{i},k_{i}\}}e^{jp^{i}}M_{r_{i}+k_{i}-2j}^{[i]}\right) = \\
&=&\sum\nolimits_{_{\substack{ j_{i}=0  \\ (1\leq i\leq g-1)}}}^{\min
\{r_{i},k_{i}\}}e^{s\left( j_{1},...,j_{g-1}\right)
}M_{r_{0}}M_{r_{1}+k_{1}-2j_{1}}^{[1]}...M_{r_{g-1}+k_{g-1}-2j_{g-1}}^{[g-1]},
\end{eqnarray*}

\noindent where $s\left( j_{1},...,j_{g-1}\right) \in 
%TCIMACRO{\U{2124} }%
%BeginExpansion
\mathbb{Z}
%EndExpansion
$ and the last summation is over the $g-1$ indices $j_{1},...,j_{g-1}.$
Since 
\begin{equation*}
\dim _{\mathbb{F}_{q}}\left( M_{r_{0}}\otimes
M_{r_{1}+k_{1}-2j_{1}}^{[1]}\otimes ...\otimes
M_{r_{g-1}+k_{g-1}-2j_{g-1}}^{[g-1]}\right) <\dim _{\mathbb{F}_{q}}M
\end{equation*}
for any value of $j_{1},...,j_{g-1}$, by induction assumption we can compute
the Jordan-H\"{o}lder constituents of $J_{h}\tprod%
\nolimits_{i=1}^{g-1}M_{k_{i}}^{[i]}$ in the standard form.

\textit{\noindent }\noindent \noindent \noindent \noindent \underline{%
\textit{Case 2}}\textit{: For any }$i$\textit{, }$0\leq i\leq g-1$,\textit{\
the representation }$M_{k_{i}}$\textit{\ is irreducible.}

\noindent By the previous lemma, we can assume - up to twistings by powers
of $e$ - that we have written $M_{k_{i}}^{[i]}=\tprod%
\nolimits_{j=0}^{g-1}M_{r_{j}^{(i)}}^{[j]}$ for any $0\leq i\leq g-1$, where 
$0\leq r_{0}^{(i)},...,r_{g-1}^{(i)}\leq p-1$. Then:%
\begin{equation}
M=\prod\nolimits_{i=0}^{g-1}M_{k_{i}}^{[i]}=\prod\nolimits_{j=0}^{g-1}\left(
\prod\nolimits_{i=0}^{g-1}M_{r_{j}^{(i)}}\right) ^{[j]}.  \label{f2}
\end{equation}

\noindent Applying \textit{only} the product formula (cf. Corollary \ref%
{product-corollary}), we can write:%
\begin{equation}
\left( \prod\nolimits_{i=0}^{g-1}M_{r_{j}^{(i)}}\right)
^{[j]}=\sum\nolimits_{\alpha _{j}\in A_{j}}e^{s_{\alpha _{j}}}M_{r_{\alpha
_{j}}}^{[j]},  \label{f3}
\end{equation}

\noindent where, for any $0\leq j\leq g-1,$ $A_{j}$ is a non-empty finite
set and $s_{\alpha _{j}},r_{\alpha _{j}}\geq 0$ for $\alpha _{j}\in A_{j}$.
Combining (\ref{f2}) and (\ref{f3}) we obtain:%
\begin{equation}
M=\sum\nolimits_{\substack{ \alpha _{j}\in A_{j}  \\ (0\leq j\leq g-1)}}%
e^{s(\alpha _{0},...,\alpha _{g-1})}M_{r_{\alpha _{0}}}M_{r_{\alpha
_{1}}}^{[1]}...M_{r_{\alpha _{g-1}}}^{[g-1]},  \label{f1}
\end{equation}

\noindent where $s(\alpha _{0},...,\alpha _{g-1})\in 
%TCIMACRO{\U{2124} }%
%BeginExpansion
\mathbb{Z}
%EndExpansion
$ and the summation is over the $g$-tuples $\left( \alpha _{0},...,\alpha
_{g-1}\right) \in A_{0}\times ...\times A_{g-1}$. If each of the sets $%
A_{0},...,A_{g-1}$ contains exactly one element, then for any $0\leq j\leq
g-1$, at most one element in $\{r_{j}^{(0)},...,r_{j}^{(g-1)}\}$ is
positive. Indeed, if this were not the case, there would be some $j$ such
that $r_{j}^{(a)},r_{j}^{(b)}>0$ for some $a,b$ with $0\leq a<b\leq g-1$;
then by \textit{only} applying the product formula we would obtain: 
\begin{equation*}
\left( \prod\nolimits_{i=0}^{g-1}M_{r_{j}^{(i)}}\right) ^{[j]}=\left[ \left(
\prod\nolimits_{\substack{ i=0  \\ i\neq a,b}}^{g-1}M_{r_{j}^{(i)}}\right)
\left(
M_{r_{j}^{(a)}+r_{j}^{(b)}}+eM_{r_{j}^{(a)}-1}M_{r_{j}^{(b)}-1}\right) %
\right] ^{[j]}.
\end{equation*}

\noindent Since $r_{j}^{(a)}-1,r_{j}^{(b)}-1\geq 0$, the left hand side
above contains at least two non-zero summand, contradicting the fact that by 
\textit{only} applying the product formula we could also write $\left(
\prod\nolimits_{i=0}^{g-1}M_{r_{j}^{(i)}}\right)
^{[j]}=e^{s_{j}}M_{r_{j}}^{[j]}$ for some integers $s_{j},r_{j}$. \noindent
We conclude that if each of the sets $A_{0},...,A_{g-1}$ contains exactly
one element, then $M$ is irreducible and (\ref{f2}) is the standard Jordan-H%
\"{o}lder form of $M$.

\noindent If there is $0\leq j\leq g-1$ such that $A_{j}$ has at least two
elements, then in (\ref{f1}) at least tow non-zero terms appear, so that
each of the summand of (\ref{f1}) has dimension strictly less than $\dim _{%
\mathbb{F}_{q}}M$, and by induction we are done. $\blacksquare $

\bigskip

\subsection{Families of intertwining operators for $g>1\label{intertnow}$}

For $g=1$, one has available two intertwining operators acting on $\mathbb{F}%
_{p}[G]$-modules and shifting weights by $p\pm 1$, namely the Dickson
invariant $\Theta _{q}$ and the derivation map $D$ (cf. \ref{g=1_op}). For $%
g>1$, equation $(\Phi _{g})$ and the existence of partial Hasse invariants
and theta operators acting on spaces of mod $p$ Hilbert modular forms (cf. 
\cite{AG}) suggest that there should be other intertwining operators between
modular representations of $G$, generalizing $\Theta _{q}$ and $D$. In this
section we will construct such operators.

Unless otherwise specified, we will always assume $g>1$, and we will
consider all the tensor product over $\mathbb{F}_{q}$ ($q=p^{g}$).

\subsubsection{Generalized Dickson invariants\label{GDI}}

\begin{definition}
For any integer $\beta $ such that $1\leq \beta \leq g-1$, the (non-twisted) 
\textbf{generalized }$\beta $\textbf{th Dickson operator} is the element 
\begin{equation*}
\Theta _{\beta }=X\otimes Y^{p^{g-\beta }}-Y\otimes X^{p^{g-\beta }}
\end{equation*}

\noindent \noindent of the $G$-module $M_{1}\otimes M_{p^{g-\beta }}^{[\beta
]}$.

For integers $\alpha ,\beta $ such that $0\leq \alpha \leq g-1$ and $1\leq
\beta \leq g-1$, the $\alpha $\textbf{-twisted generalized }$\beta $\textbf{%
th Dickson operator} is the element 
\begin{equation*}
\Theta _{\beta }^{[\alpha ]}=X\otimes Y^{p^{g-\beta }}-Y\otimes
X^{p^{g-\beta }}
\end{equation*}

\noindent of the $G$-module $M_{1}^{[\alpha ]}\otimes M_{p^{g-\beta
}}^{[\alpha +\beta ]}$.
\end{definition}

\begin{lemma}
\label{conti}Let $k,h$ be two non-negative integers and let $\alpha ,\beta $
be two integers such that $0\leq \alpha \leq g-1$ and $1\leq \beta \leq g-1$%
. Multiplication by $\Theta _{\beta }^{[\alpha ]}$ in the $\mathbb{F}_{q}[G]$%
-algebra $\mathbb{F}_{q}[X,Y]^{[\alpha ]}\otimes \mathbb{F}%
_{q}[X,Y]^{[\alpha +\beta ]}$ induces an injective $G$-homomorphism:%
\begin{equation*}
\Theta _{\beta }^{[\alpha ]}:\det\nolimits^{p^{\alpha }}\otimes
M_{k}^{[\alpha ]}\otimes M_{h}^{[\alpha +\beta ]}\mathbb{\rightarrow }%
M_{k+1}^{[\alpha ]}\otimes M_{h+p^{g-\beta }}^{[\alpha +\beta ]}.
\end{equation*}
\end{lemma}

\textbf{Proof }We can assume $\alpha =0$. To prove $G$-equivariance of the
map $\Theta _{\beta }$,\ it is enough to show that $\gamma \Theta _{\beta
}=\det \gamma \cdot \Theta _{\beta }$ for all $\gamma \in G$. Let $\gamma
=\left( \QATOP{a}{c}\QATOP{b}{d}\right) \in G$; then $\gamma \Theta _{\beta
} $ equals:%
\begin{eqnarray*}
&&(aX+cY)\otimes (b^{p^{\beta }}X+d^{p^{\beta }}Y)^{p^{g-\beta
}}-(bX+dY)\otimes (a^{p^{\beta }}X+c^{p^{\beta }}Y)^{p^{g-\beta }} \\
&=&(aX+cY)\otimes (b^{q}X^{p^{g-\beta }}+d^{q}Y^{p^{g-\beta
}})-(bX+dY)\otimes (a^{q}X^{p^{g-\beta }}+c^{q}Y^{p^{g-\beta }}) \\
&=&(aX+cY)\otimes (bX^{p^{g-\beta }}+dY^{p^{g-\beta }})-(bX+dY)\otimes
(aX^{p^{g-\beta }}+cY^{p^{g-\beta }}) \\
&=&adX\otimes Y^{p^{g-\beta }}+bcY\otimes X^{p^{g-\beta }}-bcX\otimes
Y^{p^{g-\beta }}-adY\otimes X^{p^{g-\beta }} \\
&=&\det \gamma \cdot \left( X\otimes Y^{p^{g-\beta }}-Y\otimes X^{p^{g-\beta
}}\right) \\
&=&\det \gamma \cdot \Theta _{\beta }.\text{ }
\end{eqnarray*}

\noindent To show injectivity of $\Theta _{\beta }$, notice that there is an
isomorphism of $\mathbb{F}_{q}[G]$-algebras $\mathbb{F}_{q}[X,Y]\otimes 
\mathbb{F}_{q}[X,Y]^{[\beta ]}\simeq \mathbb{F}_{q}[Z,W,T^{p^{\beta
}},U^{p^{\beta }}]$ obtained by sending the ordered tuple $(X\otimes
1,Y\otimes 1,1\otimes X,1\otimes Y)$ into the ordered tuple $%
(Z,W,T^{p^{\beta }},U^{p^{\beta }})$, were we are letting $G$ acts on $%
Z,W,T,U$ as follows: for $\gamma =\left( \QATOP{a}{c}\QATOP{b}{d}\right) \in
G$,%
\begin{eqnarray*}
\gamma Z &=&aZ+cW, \\
\gamma W &=&bZ+dW, \\
\gamma T &=&aT+cU, \\
\gamma U &=&bT+dU.
\end{eqnarray*}%
Under the above identification, the map $\Theta _{\beta }$ corresponds to
multiplication by $ZU^{q}-WT^{q}$ on $\mathbb{F}_{q}[Z,W,T^{p^{\beta
}},U^{p^{\beta }}]$, and it is therefore injective. $\blacksquare $

\bigskip

In addition to the above operators, the classical Dickson invariant also
gives rise to an intertwining map:

\begin{proposition}
Let $k$ be a non-negative integer and let $\alpha $ be an integer such that $%
0\leq \alpha \leq g-1$. Let $\Theta ^{\lbrack \alpha ]}=XY^{q}-YX^{q}$ be
the classical Dickson invariant, viewed as an element of $M_{q+1}^{[\alpha
]} $. Multiplication by $\Theta ^{\lbrack \alpha ]}$ in the $\mathbb{F}%
_{q}[G]$-algebra $\mathbb{F}_{q}[X,Y]^{[\alpha ]}$ induces an injective $G$%
-homomorphism:%
\begin{equation*}
\Theta ^{\lbrack \alpha ]}:\det\nolimits^{p^{\alpha }}\otimes M_{k}^{[\alpha
]}\longrightarrow M_{k+(q+1)}^{[\alpha ]}.
\end{equation*}
\end{proposition}

\textbf{Proof }This follows from section \ref{Dickson}. $\blacksquare $

\bigskip

Notice that the operators $\Theta _{\beta }^{[\alpha ]}$ and $\Theta
^{\lbrack \alpha ]}$ ($0\leq \alpha \leq g-1$, $1\leq \beta \leq g-1$)
pairwise commute, as it follows by seeing them as multiplication by
polynomials in some polynomial algebra over $\mathbb{F}_{q}$ (cf. the end of
the proof of Lemma \ref{conti}).

\begin{remark}
Let us fix a convention that will make the notation easier in the sequel.
For non-negative integers $k_{0},...,k_{g-1}$, the $G$-module $%
M_{k_{0}}\otimes M_{k_{1}}^{[1]}\otimes ...\otimes M_{k_{g-1}}^{[g-1]}$ will
be identified with the $G$-module obtained by permuting in any possible way
the tensor factors. Also, for integers $\alpha ,\beta $ and any $G$-module $%
M $, the notation $M^{[\alpha +\beta ]}$ will denote the $\gamma $th
Frobenius twist of $M$, where $\gamma $ is the smallest non-negative integer
such that $\gamma \equiv \alpha +\beta (\func{mod}g)$.
\end{remark}

We can summarize the above results as follows:

\begin{theorem}
\label{TTh}Let us fix non-negative integers $k_{0},...,k_{g-1}$. For any
integers $\alpha ,\beta $ subject to the constraints $0\leq \alpha \leq g-1$
and $1\leq \beta \leq g-1$, there are pairwise commuting injective $G$%
-intertwining operators as follows:%
\begin{equation*}
\Theta _{\beta }^{[\alpha ]}:\det\nolimits^{p^{\alpha }}\otimes
\dbigotimes\nolimits_{i}M_{k_{i}}^{[i]}\longrightarrow \left(
\dbigotimes\nolimits_{i\neq \alpha ,\alpha +\beta }M_{k_{i}}^{[i]}\right)
\otimes M_{k_{\alpha }+1}^{[\alpha ]}\otimes M_{k_{\alpha +\beta
}+p^{g-\beta }}^{[\alpha +\beta ]};
\end{equation*}%
\begin{equation*}
\Theta ^{\lbrack \alpha ]}:\det\nolimits^{p^{\alpha }}\otimes
\dbigotimes\nolimits_{i}M_{k_{i}}^{[i]}\longrightarrow \left(
\dbigotimes\nolimits_{i\neq \alpha }M_{k_{i}}^{[i]}\right) \otimes
M_{k_{\alpha }+(q+1)}^{[\alpha ]},
\end{equation*}

where the tensor products indices run over the integers $i$ such that $0\leq
i\leq g-1$, unless otherwise specified.
\end{theorem}

\begin{remark}
\label{AG1}The operators $\Theta _{g-1}^{[\alpha ]}$ for $0\leq \alpha \leq
g-1$ give, under suitable assumptions, cohomological analogues of the theta
operators defined in \cite{AG} in the context of Hilbert modular forms. We
do not know of any geometric interpretation of the other generalized Dickson
operators.
\end{remark}

We can picture the weight shiftings allowed by the $g(g-1)+g=g^{2}$
generalized Dickson operators with the following self-explanatory tables:

\begin{equation*}
\fbox{$%
\begin{array}{cc}
\Theta _{1} & \left( 1,p^{g-1},0,0,...,0,0\right) \\ 
\Theta _{1}^{[1]} & \left( 0,1,p^{g-1},0,...,0,0\right) \\ 
\Theta _{1}^{[2]} & \left( 0,0,1,p^{g-1},...,0,0\right) \\ 
... & ... \\ 
\Theta _{1}^{[g-2]} & \left( 0,0,0,0,...,1,p^{g-1}\right) \\ 
\Theta _{1}^{[g-1]} & \left( p^{g-1},0,0,0,...,0,1\right)%
\end{array}%
$}\text{ \ \ }\fbox{$%
\begin{array}{cc}
\Theta _{2} & \left( 1,0,p^{g-2},0,...,0,0\right) \\ 
\Theta _{2}^{[1]} & \left( 0,1,0,p^{g-2},...,0,0\right) \\ 
\Theta _{2}^{[2]} & \left( 0,0,1,0,...,0,0\right) \\ 
... & ... \\ 
\Theta _{2}^{[g-2]} & \left( p^{g-2},0,0,0,...,1,0\right) \\ 
\Theta _{2}^{[g-1]} & \left( 0,p^{g-2},0,0,...,0,1\right)%
\end{array}%
$}
\end{equation*}

\begin{equation*}
...
\end{equation*}

\begin{equation*}
\fbox{$%
\begin{array}{cc}
\Theta _{g-1} & \left( 1,0,0,0,...,0,p\right) \\ 
\Theta _{g-1}^{[1]} & \left( p,1,0,0,...,0,0\right) \\ 
\Theta _{g-1}^{[2]} & \left( 0,p,1,0,...,0,0\right) \\ 
... & ... \\ 
\Theta _{g-1}^{[g-2]} & \left( 0,0,0,0,...,1,0\right) \\ 
\Theta _{g-1}^{[g-1]} & \left( 0,0,0,0,...,p,1\right)%
\end{array}%
$}\text{ \ \ }\fbox{$%
\begin{array}{cc}
\Theta & \left( q+1,0,0,...,0,0\right) \\ 
\Theta ^{\lbrack 1]} & \left( 0,q+1,0,...,0,0\right) \\ 
\Theta ^{\lbrack 2]} & \left( 0,0,q+1,...,0,0\right) \\ 
... & ... \\ 
\Theta ^{\lbrack g-2]} & \left( 0,0,0,...,q+1,0\right) \\ 
\Theta ^{\lbrack g-1]} & \left( 0,0,0,...,0,q+1\right)%
\end{array}%
$}
\end{equation*}

\bigskip

For example, if $g=2$ the generalized Dickson operators give all and only
the weight shiftings of the form:%
\begin{equation*}
a_{1}(1,p)+a_{2}(p,1)+a_{3}(0,p^{2}+1)+a_{4}(p^{2}+1,0),
\end{equation*}

\noindent for any non-negative integers $a_{1},a_{2},a_{3},a_{4}.$ For $g>2$
a new phenomenon occurs, as the operators do not allow weight shiftings of
the form:%
\begin{equation*}
(1,p,0,...,0,0),(0,1,p,...,0,0),...,(0,0,0,...,1,p),(p,0,0,...,0,1).
\end{equation*}

\noindent This happens not because of limitations intrinsic to our
intertwining maps, but because of the structure of $G$-modules:

\begin{proposition}
Assume $g>2$ and let $k,h$ be integers such that $0\leq k,h\leq p-1$. For
any integer $\alpha $ such that $0\leq \alpha \leq g-1$ and any integer $m$,
there are no $G$-module embeddings $\det^{m}\otimes M_{k}^{[\alpha ]}\otimes
M_{h}^{[\alpha +1]}\mathbb{\rightarrow }M_{k+1}^{[\alpha ]}\otimes
M_{h+p}^{[\alpha +1]}$.
\end{proposition}

\textbf{Proof }It is enough to prove the non existence of embeddings for $%
\alpha =0$. Using $(\Phi _{g})$ and $(\Delta _{g})$\ we have, in $K_{0}(G)$:%
\begin{equation*}
M_{h+p}^{[1]}=M_{h}^{[1]}M_{1}^{[2]}+e^{p(h+1)}M_{p-h-2}^{[1]}.
\end{equation*}

\noindent If $k\neq p-1$, as $g>2$, we deduce that the Jordan-H\"{o}lder
factors of $M_{k+1}\otimes M_{h+p}^{[1]}$ are $M_{k+1}\otimes
M_{h}^{[1]}\otimes M_{1}^{[2]}$ and $\det\nolimits^{p(h+1)}\otimes
M_{k+1}\otimes M_{p-h-2}^{[1]}$, unless $h=p-1$, in which case only the
first factor occurs. None of these factors coincides with $\det^{m}\otimes
M_{k}\otimes M_{h}^{[1]}$.

\noindent If $k=p-1$, write $M_{p}=M_{1}^{[1]}+eM_{p-2}$ in $K_{0}(G)$.
Applying $(\Pi _{g})$ we obtain:%
\begin{eqnarray*}
M_{p}M_{h+p}^{[1]} &=&\left( M_{1}^{[1]}+eM_{p-2}\right) \left(
M_{h}^{[1]}M_{1}^{[2]}+e^{p(h+1)}M_{p-h-2}^{[1]}\right) \\
&=&M_{h+1}^{[1]}M_{1}^{[2]}+e^{p}M_{h-1}^{[1]}M_{1}^{[2]}+e^{p(h+1)}M_{p-h-1}^{[1]}+e^{p(h+2)}M_{p-h-3}^{[1]}
\\
&&+eM_{p-2}M_{h}^{[1]}M_{1}^{[2]}+e^{p(h+1)+1}M_{p-2}M_{p-h-2}^{[1]}.
\end{eqnarray*}

\noindent If $h\neq p-1$, the above formula shows that none of the Jordan-H%
\"{o}lder factors of $M_{p}\otimes M_{h+p}^{[1]}$ equals $\det^{m}\otimes
M_{p-1}\otimes M_{h}^{[1]}$. If $h=p-1$, we have:

\begin{eqnarray*}
M_{p}M_{2p-1}^{[1]}
&=&M_{1}^{[2]}M_{1}^{[2]}+e^{p}M_{p-2}^{[1]}M_{1}^{[2]}+e^{p}M_{p-2}^{[1]}M_{1}^{[2]}+eM_{p-2}M_{p-1}^{[1]}M_{1}^{[2]}
\\
&=&M_{2}^{[2]}+e^{p^{2}}+2e^{p}M_{p-2}^{[1]}M_{1}^{[2]}+eM_{p-2}M_{p-1}^{[1]}M_{1}^{[2]},
\end{eqnarray*}

\noindent and $\det^{m}\otimes M_{p-1}\otimes M_{p-1}^{[1]}$ is not a
constituent of $M_{p}\otimes M_{2p-1}^{[1]}$ if $p\neq 2.$ If $p=2$,
decomposing $M_{2}^{[2]}$ we get to the same conclusion. $\blacksquare $

\bigskip

We conclude this section by noticing the following consequence of
Proposition \ref{Dickson prop}:

\begin{proposition}
Let us fix non-negative integers $k_{0},...,k_{g-1}$. For any integer $%
\alpha $ such that $0\leq \alpha \leq g-1$ consider the $G$-map $\Theta
^{\lbrack \alpha ]}:\det\nolimits^{p^{\alpha }}\otimes
\tbigotimes\nolimits_{i}M_{k_{i}}^{[i]}\mathbb{\rightarrow }\left(
\tbigotimes\nolimits_{i\neq \alpha }M_{k_{i}}^{[i]}\right) \otimes
M_{k_{\alpha }+(q+1)}^{[\alpha ]}$. We have:%
\begin{equation*}
\limfunc{coker}\Theta ^{\lbrack \alpha ]}\simeq \left(
\dbigotimes\nolimits_{i\neq \alpha }M_{k_{i}}^{[i]}\right) \otimes \left[ 
\limfunc{Ind}\nolimits_{B}^{G}\left( \eta ^{k_{\alpha }+2}\right) \right]
^{[\alpha ]},
\end{equation*}

where $B$ is the subgroup of $G$ consisting of upper triangular matrices,
and $\eta $ is the character of $B$ defined extending the character $%
\limfunc{diag}(a,b)\mapsto a$ of the standard maximal torus of $G.$
\end{proposition}

\begin{remark}
Even though the Jordan-H\"{o}lder constituents of $\limfunc{coker}\Theta
_{\beta }^{[\alpha ]}$ can be explicitly computed using the results of this
paper, We do not know of any interesting description of the cokernel of the
operators $\Theta _{\beta }^{[\alpha ]}$.
\end{remark}

\subsubsection{Generalized $D$-operators\label{GDO}}

Let us denote by $\partial _{X}$ (resp. $\partial _{Y}$) the operator of
partial derivation with respect to $X$ (resp. $Y$) acting on the polynomial
algebra $\mathbb{F}_{q}[X,Y]$; if $f\in \mathbb{F}_{q}[X,Y]$, denote by the
same symbol the $\mathbb{F}_{q}$-vector space endomorphism of $\mathbb{F}%
_{q}[X,Y]$ induced by multiplication by $f$. The operators $\partial
_{X}\otimes f,\partial _{Y}\otimes f,f$ $\otimes \partial _{X}$ and $%
f\otimes \partial _{Y}$ are therefore derivation of the $\mathbb{F}_{q}$%
-algebra $\mathbb{F}_{q}[X,Y]\otimes \mathbb{F}_{q}[X,Y].$

\begin{definition}
Let $k,h$ be two non-negative integers. For any integer $\beta $ such that $%
1\leq \beta \leq g-1$, the (non-twisted) \textbf{generalized }$\beta $%
\textbf{th }$D$\textbf{-operator} is the $\mathbb{F}_{q}$-vector space
homomorphism: 
\begin{equation*}
D_{\beta }=\partial _{X}\otimes X^{p^{g-\beta }}+\partial _{Y}\otimes
Y^{p^{g-\beta }}:M_{k}\otimes M_{h}^{[\beta ]}\longrightarrow M_{k-1}\otimes
M_{h+p^{g-\beta }}^{[\beta ]}.
\end{equation*}

For any integers $\alpha ,\beta $ such that $0\leq \alpha \leq g-1$ and $%
1\leq \beta \leq g-1$, the $\alpha $\textbf{-twisted generalized }$\beta $%
\textbf{th }$D$\textbf{-operator} is the $\mathbb{F}_{q}$-vector space
homomorphism: 
\begin{equation*}
D_{\beta }^{[\alpha ]}=\partial _{X}\otimes X^{p^{g-\beta }}+\partial
_{Y}\otimes Y^{p^{g-\beta }}:M_{k}^{[\alpha ]}\otimes M_{h}^{[\alpha +\beta
]}\longrightarrow M_{k-1}^{[\alpha ]}\otimes M_{h+p^{g-\beta }}^{[\alpha
+\beta ]}.
\end{equation*}
\end{definition}

\begin{lemma}
Let $k,h$ be two non-negative integers and let $\alpha ,\beta $ be integers
such that $0\leq \alpha \leq g-1$ and $1\leq \beta \leq g-1$. The operator $%
D_{\beta }^{[\alpha ]}:M_{k}^{[\alpha ]}\otimes M_{h}^{[\alpha +\beta ]}%
\mathbb{\rightarrow }M_{k-1}^{[\alpha ]}\otimes M_{h+p^{g-\beta }}^{[\alpha
+\beta ]}$ is a $G$-homomorphism; it is injective if $0<k\leq p-1$ and $%
0\leq h\leq p-1$.
\end{lemma}

\textbf{Proof }By twisting, we can assume that $\alpha =0$. Fix $f_{1}\in
M_{k}$, $f_{2}\in M_{h}^{[\beta ]}$ and let $\gamma =\left( \QATOP{a}{c}%
\QATOP{b}{d}\right) \in G$; denote by $\gamma ^{\sigma ^{\beta }}$ the
matrix $\left( \QATOP{a^{\sigma ^{\beta }}}{c^{\sigma ^{\beta }}}\QATOP{%
b^{\sigma ^{\beta }}}{d^{\sigma ^{\beta }}}\right) $, where $\sigma $
denotes the arithmetic Frobenius element of $\limfunc{Gal}(\mathbb{F}_{q}/%
\mathbb{F}_{p}).$ $D_{\beta }\left( \gamma (f_{1}\otimes f_{2})\right) $
equals: 
\begin{eqnarray*}
&&\left[ a\cdot \left( \partial _{X}f_{1}\right) (\gamma X,\gamma Y)+b\cdot
\left( \partial _{Y}f_{1}\right) (\gamma X,\gamma Y)\right] \otimes
X^{p^{g-\beta }}f_{2}(\gamma ^{\sigma ^{\beta }}X,\gamma ^{\sigma ^{\beta
}}Y) \\
&&+\left[ c\cdot \left( \partial _{X}f_{1}\right) (\gamma X,\gamma Y)+d\cdot
\left( \partial _{Y}f_{1}\right) (\gamma X,\gamma Y)\right] \otimes
Y^{p^{g-\beta }}f_{2}(\gamma ^{\sigma ^{\beta }}X,\gamma ^{\sigma ^{\beta
}}Y) \\
&=&\left( \partial _{X}f_{1}\right) (\gamma X,\gamma Y)\otimes
(aX^{p^{g-\beta }}+cY^{p^{g-\beta }})f_{2}(\gamma ^{\sigma ^{\beta
}}X,\gamma ^{\sigma ^{\beta }}Y) \\
&&+\left( \partial _{Y}f_{1}\right) (\gamma X,\gamma Y)\otimes
(bX^{p^{g-\beta }}+dY^{p^{g-\beta }})f_{2}(\gamma ^{\sigma ^{\beta
}}X,\gamma ^{\sigma ^{\beta }}Y) \\
&=&\gamma D_{\beta }\left( f_{1}\otimes f_{2}\right) .
\end{eqnarray*}

\noindent For the injectivity statement, notice that if $0<k\leq p-1$ and $%
0\leq h\leq p-1$, then $M_{k}\otimes M_{h}^{[\beta ]}$ is an irreducible $G$%
-module, so it is enough to show that $D_{\beta }$ is non-zero on $%
M_{k}\otimes M_{h}^{[\beta ]}$. We have $D_{\beta }(X^{k}\otimes
X^{h})=kX^{k-1}\otimes X^{h+p^{g-\beta }}$, and this is non-zero as $k$ is
prime with $p$. $\blacksquare $

\bigskip

In addition to the above operators, the $D$-map defined by Serre also gives
an intertwining map:

\begin{proposition}
Let $k$ be a non-negative integer and let $\alpha $ be an integer such that $%
0\leq \alpha \leq g-1$. Then the Frobenius twists of Serre's operator $%
D^{[\alpha ]}=X^{q}\partial _{X}+Y^{q}\partial _{Y}$ define $G$%
-homomorphisms:%
\begin{equation*}
D^{[\alpha ]}:M_{k}^{[\alpha ]}\longrightarrow M_{k+(q-1)}^{[\alpha ]}
\end{equation*}

\noindent which are injective if $1\leq k\leq p-1$.
\end{proposition}

\textbf{Proof }After twisting, we can assume $\alpha =0$. The result then
follows from section \ref{D-map} and the irreducibility of $M_{k}^{[\alpha
]} $ in the range $1\leq k\leq p-1$. $\blacksquare $

\bigskip

We can summarize the above results as follows:

\begin{theorem}
\label{DTh}Let us fix non-negative integers $k_{0},...,k_{g-1}$. For any
integers $\alpha ,\beta $ subject to the constraints $0\leq \alpha \leq g-1$
and $1\leq \beta \leq g-1$, there are $G$-intertwining operators as follows:%
\begin{equation*}
D_{\beta }^{[\alpha
]}:\dbigotimes\nolimits_{i}M_{k_{i}}^{[i]}\longrightarrow \left(
\dbigotimes\nolimits_{i\neq \alpha ,\alpha +\beta }M_{k_{i}}^{[i]}\right)
\otimes M_{k_{\alpha }-1}^{[\alpha ]}\otimes M_{k_{\alpha +\beta
}+p^{g-\beta }}^{[\alpha +\beta ]};
\end{equation*}%
\begin{equation*}
D^{[\alpha ]}:\dbigotimes\nolimits_{i}M_{k_{i}}^{[i]}\longrightarrow \left(
\dbigotimes\nolimits_{i\neq \alpha }M_{k_{i}}^{[i]}\right) \otimes
M_{k_{\alpha }+(q-1)}^{[\alpha ]},
\end{equation*}

\noindent \noindent where the tensor product indices run over the integers $%
i $ such that $0\leq i\leq g-1$, unless otherwise specified. If $0<k_{\alpha
}\leq p-1$, then $D^{[\alpha ]}$ is injective; if in addition $0\leq
k_{\alpha +\beta }\leq p-1$, then $D_{\beta }^{[\alpha ]}$ is injective.
\end{theorem}

\begin{remark}
\label{AG2}The operators $D_{g-1}^{[\alpha ]}$ for $0\leq \alpha \leq g-1$
give, under suitable assumptions, cohomological analogues of the partial
Hasse invariants defined in \cite{AG} in the context of mod $p$ Hilbert
modular forms. We do not know of any geometric interpretation of the other $%
D $-maps introduced above.
\end{remark}

We can picture the weight shiftings allowed by the $g(g-1)+g=g^{2}$
generalized $D$-maps as follows:

\begin{equation*}
\fbox{$%
\begin{array}{cc}
D_{1} & \left( -1,p^{g-1},0,0,...,0,0\right) \\ 
D_{1}^{[1]} & \left( 0,-1,p^{g-1},0,...,0,0\right) \\ 
D_{1}^{[2]} & \left( 0,0,-1,p^{g-1},...,0,0\right) \\ 
... & ... \\ 
D_{1}^{[g-2]} & \left( 0,0,0,0,...,-1,p^{g-1}\right) \\ 
D_{1}^{[g-1]} & \left( p^{g-1},0,0,0,...,0,-1\right)%
\end{array}%
$}\text{ \ \ }\fbox{$%
\begin{array}{cc}
D_{2} & \left( -1,0,p^{g-2},0,...,0,0\right) \\ 
D_{2}^{[1]} & \left( 0,-1,0,p^{g-2},...,0,0\right) \\ 
D_{2}^{[2]} & \left( 0,0,-1,0,...,0,0\right) \\ 
... & ... \\ 
D_{2}^{[g-2]} & \left( p^{g-2},0,0,0,...,-1,0\right) \\ 
D_{2}^{[g-1]} & \left( 0,p^{g-2},0,0,...,0,-1\right)%
\end{array}%
$}
\end{equation*}

\begin{equation*}
...
\end{equation*}

\begin{equation*}
\fbox{$%
\begin{array}{cc}
D_{g-1} & \left( -1,0,0,0,...,0,p\right) \\ 
D_{g-1}^{[1]} & \left( p,-1,0,0,...,0,0\right) \\ 
D_{g-1}^{[2]} & \left( 0,p,-1,0,...,0,0\right) \\ 
... & ... \\ 
D_{g-1}^{[g-2]} & \left( 0,0,0,0,...,-1,0\right) \\ 
D_{g-1}^{[g-1]} & \left( 0,0,0,0,...,p,-1\right)%
\end{array}%
$}\text{ \ \ }\fbox{$%
\begin{array}{cc}
D & \left( q-1,0,0,...,0,0\right) \\ 
D^{[1]} & \left( 0,q-1,0,...,0,0\right) \\ 
D^{[2]} & \left( 0,0,q-1,...,0,0\right) \\ 
... & ... \\ 
D^{[g-2]} & \left( 0,0,0,...,q-1,0\right) \\ 
D^{[g-1]} & \left( 0,0,0,...,0,q-1\right)%
\end{array}%
$}
\end{equation*}

\bigskip

Similarly to what happened for the generalized Dickson operators, the non
existence of shiftings of the form

\begin{equation*}
(-1,p,0,...,0,0),(0,-1,p,...,0,0),...,(0,0,0,...,-1,p),(p,0,0,...,0,-1)
\end{equation*}

\noindent when $g>2$ is a consequence of the structure of the irreducible $G$%
-modules:

\begin{proposition}
Assume $g>2$ and let $k,h$ be integers such that $0\leq k,h\leq p-1$. For
any integer $\alpha $ such that $0\leq \alpha \leq g-1$ and any integer $m$,
there are no $G$-module embeddings $\det^{m}\otimes M_{k}^{[\alpha ]}\otimes
M_{h}^{[\alpha +1]}\mathbb{\rightarrow }M_{k-1}^{[\alpha ]}\otimes
M_{h+p}^{[\alpha +1]}$.
\end{proposition}

\textbf{Proof }It is enough to consider the case $\alpha =0$; we can also
assume that $k\neq 0$. Using formulae $(\Phi _{g})$ and $(\Delta _{g})$\ we
have, in $K_{0}(G)$:%
\begin{equation*}
M_{h+p}^{[1]}=M_{h}^{[1]}M_{1}^{[2]}+e^{p(h+1)}M_{p-h-2}^{[1]}.
\end{equation*}

\noindent As $g>2$, the Jordan-H\"{o}lder factors of $M_{k-1}\otimes
M_{h+p}^{[1]}$ are $M_{k-1}\otimes M_{h}^{[1]}\otimes M_{1}^{[2]}$ and $%
\det\nolimits^{p(h+1)}\otimes M_{k-1}\otimes M_{p-h-2}^{[1]}$, unless $h=p-1$%
, in which case only the first factor occurs: none of these factors
coincides with $\det^{m}\otimes M_{k}\otimes M_{h}^{[1]}.\blacksquare $

\bigskip

We conclude by noticing the following consequence of Theorem \ref{cris}:

\begin{proposition}
Let us fix non-negative integers $k_{0},...,k_{g-1}$; let $\alpha $ be an
integer such that $0\leq \alpha \leq g-1$ and assume $2\leq k_{\alpha }\leq
p-1$, $k_{\alpha }\neq \frac{q+1}{2}$. Consider the injective $G$-map $%
D^{[\alpha ]}:D^{[\alpha ]}:\tbigotimes\nolimits_{i}M_{k_{i}}^{[i]}\mathbb{%
\rightarrow }\left( \tbigotimes\nolimits_{i\neq \alpha
}M_{k_{i}}^{[i]}\right) \otimes M_{k_{\alpha }+(q-1)}^{[\alpha ]}$. We have:%
\begin{equation*}
\limfunc{coker}D^{[\alpha ]}\simeq \left( \tbigotimes\nolimits_{i\neq \alpha
}M_{k_{i}}^{[i]}\right) \otimes \left[ \bar{\Xi}\left( \chi ^{k_{\alpha
}}\right) \right] ^{[\alpha ]},
\end{equation*}

where: $\bar{\Xi}\left( \chi ^{k_{\alpha }}\right) =H_{\limfunc{cris}}^{1}(%
\mathcal{C}_{/\mathbb{F}_{q}})_{-k_{\alpha }}\otimes _{W(\mathbb{F}_{q})}%
\mathbb{F}_{q}$, $\mathcal{C}$ is the Deligne-Lusztig variety of $SL_{2/%
\mathbb{F}_{q}}$\textit{\ and }the $(-k_{\alpha })$-eigenspace of $H_{%
\limfunc{cris}}^{1}(\mathcal{C}_{/\mathbb{F}_{q}})$ is computed with respect
to the natural action of $\ker (\limfunc{Nm}_{\mathbb{F}_{q^{2}}^{\times }/%
\mathbb{F}_{q}^{\times }})$ on $H_{\limfunc{cris}}^{1}(\mathcal{C}_{/\mathbb{%
F}_{q}})$\textit{. (Here }$W(\mathbb{F}_{q})$ denotes the ring of Witt
vectors of $\mathbb{F}_{q}$).
\end{proposition}

\begin{remark}
We do not know of any interesting description of the cokernel of the
operators $D_{\beta }^{[\alpha ]}$. The Jordan-H\"{o}lder constituents of $%
\limfunc{coker}D_{\beta }^{[\alpha ]}$ can be explicitly computed using the
results of this paper.
\end{remark}

\newpage

\part{Weight shiftings for automorphic forms}

We apply the results of the previous sections to obtain weight shiftings for
automorphic forms on definite quaternion algebras whose center is a totally
real field $F$ unramified at the prime $p>2$. In section \ref{sec4}\ we
mostly\ treat the case in which the tensor factors - corresponding to the
prime decomposition of $p$ in $F$ - of the weight that we want to shift are
all of dimension greater than one: this is what we call a weight not
containing a $(2,...,2)$-block. In section \ref{22222} we consider shiftings
for irreducible weights that contain a $(2,...,2)$-block.

\section{Shiftings for weights not containing $(2,...,2)$-blocks\label{sec4}}

Let us fix some notation that will be used throughout this section and the
next one. Let $F$ be a totally real number field of degree $g$ over $%
%TCIMACRO{\U{211a} }%
%BeginExpansion
\mathbb{Q}
%EndExpansion
$, and let $p>2$ be a rational prime which is unramified in $F/%
%TCIMACRO{\U{211a} }%
%BeginExpansion
\mathbb{Q}
%EndExpansion
$. Denote by $\mathcal{O}_{F}$ the ring of integers of $F$ and write:%
\begin{equation*}
p\mathcal{O}_{F}=\dprod\nolimits_{j=1}^{r}\mathfrak{P}_{j},
\end{equation*}

\noindent where the $\mathfrak{P}_{j}$'s are distinct maximal ideals of $%
\mathcal{O}_{F}$.

Fix an integer $j$ with $1\leq j\leq r$. Let $f_{j}$ be the residual degree
of $\mathfrak{P}_{j}$ over $p%
%TCIMACRO{\U{2124} }%
%BeginExpansion
\mathbb{Z}
%EndExpansion
$, so that $\mathbb{F}_{\mathfrak{P}_{j}}:=\mathcal{O}_{F}/\mathfrak{P}_{j}$
is an extension of $\mathbb{F}_{p}=%
%TCIMACRO{\U{2124} }%
%BeginExpansion
\mathbb{Z}
%EndExpansion
/p%
%TCIMACRO{\U{2124} }%
%BeginExpansion
\mathbb{Z}
%EndExpansion
$ of degree $f_{j}$. Let $F_{\mathfrak{P}_{j}}$ be the completion of $F$ at $%
\mathfrak{P}_{j}$, and denote by $\mathcal{O}_{F_{\mathfrak{P}_{j}}}$ its
ring of integers. Fix an algebraic closure $\overline{%
%TCIMACRO{\U{211a} }%
%BeginExpansion
\mathbb{Q}
%EndExpansion
}_{p}$ of $%
%TCIMACRO{\U{211a} }%
%BeginExpansion
\mathbb{Q}
%EndExpansion
_{p}$; let $n$ be the positive least common multiple of the integers $%
f_{1},...,f_{r}$ and let $E$ be the maximal unramified extension of $%
%TCIMACRO{\U{211a} }%
%BeginExpansion
\mathbb{Q}
%EndExpansion
_{p}$\ inside $\overline{%
%TCIMACRO{\U{211a} }%
%BeginExpansion
\mathbb{Q}
%EndExpansion
}_{p}$ having\ degree $n$ over $%
%TCIMACRO{\U{211a} }%
%BeginExpansion
\mathbb{Q}
%EndExpansion
_{p}$, so that $\limfunc{Hom}(F,\overline{%
%TCIMACRO{\U{211a} }%
%BeginExpansion
\mathbb{Q}
%EndExpansion
}_{p})=\limfunc{Hom}(F,E)$. Denote by $\mathcal{O}$ the ring of integers of $%
E$ and let $\mathbb{F}$ be its residue field. Let $\sigma $ be the
arithmetic Frobenius of the extension $E/%
%TCIMACRO{\U{211a} }%
%BeginExpansion
\mathbb{Q}
%EndExpansion
_{p}$. Set:

\begin{equation*}
\limfunc{Hom}(F_{\mathfrak{P}_{j}},E)=\left\{ \sigma _{i}^{(j)}:0\leq i\leq
f_{j}-1\right\} ,
\end{equation*}

\noindent where the labeling is chosen so that, for any $i$, we have: 
\begin{equation*}
\sigma \circ \sigma _{i}^{(j)}=\sigma _{i+1}^{(j)}.
\end{equation*}%
Here the subscripts read modulo $f_{j}$ and in the range $0\leq i\leq
f_{j}-1.$

Denote by a bar the analogous morphisms for the residue fields, so that $%
\overline{\sigma }$ is the arithmetic Frobenius of the extension $\mathbb{F}/%
\mathbb{F}_{p}$, and: 
\begin{equation*}
\limfunc{Hom}(\mathbb{F}_{\mathfrak{P}_{j}},\mathbb{F})=\{\overline{\sigma }%
_{i}^{(j)}:0\leq i\leq f_{j}-1\}
\end{equation*}%
are labeled so that: 
\begin{equation*}
\overline{\sigma }\circ \overline{\sigma }_{i}^{(j)}=\overline{\sigma }%
_{i+1}^{(j)},
\end{equation*}
where the subscripts read modulo $f_{j}$ and in the range $0\leq i\leq
f_{j}-1.$

We let $\mathbb{A}_{F}$ be the topological ring of ad\`{e}les of $F$, and we
denote by $\mathbb{A}_{F}^{\infty }$ the subring of finite ad\`{e}les. We
let $\mathfrak{M}_{F,f}$ (resp. $\mathfrak{M}_{F,\infty }$) be the set of
finite (resp. infinite)\ places of $F$ and we identify $\mathfrak{M}_{F,f}$
with the set of maximal ideals of $\mathcal{O}_{F}$.

\subsection{Some motivations: geometric Hilbert modular forms\label%
{motivations}}

Denote by $d_{F}$\ the discriminant of $F/%
%TCIMACRO{\U{211a} }%
%BeginExpansion
\mathbb{Q}
%EndExpansion
$ and fix a fractional ideal $\mathfrak{a}$ of $F$ with its natural positive
cone $\mathfrak{a}^{+}$, so that $(\mathfrak{a},\mathfrak{a}^{+})$
represents an element in the strict class group of $F$. Let $N\geq 4$ be an
integer and recall that, by previous assumptions, $p$ does not divide $d_{F}$%
. Let $S$ be a scheme over $\limfunc{Spec}(%
%TCIMACRO{\U{2124} }%
%BeginExpansion
\mathbb{Z}
%EndExpansion
\lbrack \frac{1}{d_{F}}])$.

There is an $S$-scheme $\mathcal{M}$ parametrizing isomorphism classes $%
[(A,\lambda ,\iota ,\varepsilon )/T/S]$ of $(\mathfrak{a},\mathfrak{a}^{+})$%
-polarized Hilbert-Blumenthal abelian $T$-schemes $(A,\lambda )$ of relative
dimension $g$ ($T$ is an $S$-scheme), endowed with real multiplication $%
\iota $ by $\mathcal{O}_{F}$, $\mu _{N}$-level structure $\varepsilon $, and
satisfying the Deligne-Pappas condition (or, equivalently since $d_{F}$ is
invertible in $S$, satisfying the Rapoport condition). $\mathcal{M}$ has
relative dimension $g$ over $S$ and is geometrically irreducible; see \cite%
{DP} and \cite{AG} for more details.

Let $\mathbb{G}=\limfunc{Res}\nolimits_{\mathcal{O}_{F}/%
%TCIMACRO{\U{2124} }%
%BeginExpansion
\mathbb{Z}
%EndExpansion
}(\mathbb{G}_{m,\mathcal{O}_{F}})$ be the Weil restriction to $%
%TCIMACRO{\U{2124} }%
%BeginExpansion
\mathbb{Z}
%EndExpansion
$ of the algebraic $\mathcal{O}_{F}$-group $\mathbb{G}_{m,\mathcal{O}_{F}}$%
.\ \ For any scheme $T$, denote by $\mathbb{X}_{T}=\limfunc{Hom}(\mathbb{G}%
_{T},\mathbb{G}_{m,T})$ the group of characters of the base change $\mathbb{G%
}_{T}$ of $\mathbb{G}$ to $T$. If $S$ is the scheme over $\limfunc{Spec}(%
%TCIMACRO{\U{2124} }%
%BeginExpansion
\mathbb{Z}
%EndExpansion
\lbrack \frac{1}{d_{F}}])$ fixed above, a geometric $(\mathfrak{a},\mathfrak{%
a}^{+})$-polarized Hilbert modular form $f$ over $S$ having weight $\chi \in 
\mathbb{X}_{S}$ and level $\mu _{N}$ is a rule that assigns to any affine
scheme $\limfunc{Spec}(R)\rightarrow S$, any $R$-point $\left[ (A,\lambda
,\iota ,\varepsilon )/R/S\right] $ of $\mathcal{M}$, and any generator $%
\omega $ of the $R\otimes _{%
%TCIMACRO{\U{2124} }%
%BeginExpansion
\mathbb{Z}
%EndExpansion
}\mathcal{O}_{F}$-module $\Omega _{A/R}^{1}$, an element $f(A,\lambda ,\iota
,\varepsilon ,\omega )\in R$ such that:%
\begin{equation*}
f(A,\lambda ,\iota ,\varepsilon ,\alpha ^{-1}\omega )=\chi (\alpha )\cdot
f(A,\lambda ,\iota ,\varepsilon ,\omega )
\end{equation*}%
for $\alpha \in \mathbb{G}(R)$, and such that some compatibility conditions
are satisfied (cf. \cite{AG}, 5). We denote by $M_{\chi }(\mu _{N},S)$ the $%
\Gamma (S,\mathcal{O}_{S})$-module of such functions.

We remark that the formation of spaces of geometric Hilbert modular forms 
\textit{does not} commute with base change: for example, if $g>1$ and $1\leq
j\leq r,$ $0\leq i\leq f_{j}-1$, the $(j,i)$th partial Hasse invariant that
we will consider below is a non-zero, non-cuspidal modular forms over $%
\limfunc{Spec}\mathbb{F}_{\mathfrak{P}_{j}}$ that cannot be lifted to a
modular forms over $\limfunc{Spec}\mathcal{O}_{F}$: the natural reduction
morphism $M_{\chi }(\mu _{N},\mathcal{O}_{F})\mathbb{\rightarrow }M_{\chi
}(\mu _{N},\mathbb{F}_{\mathfrak{P}_{j}})$ is in general not surjective.

Assume $g>1$ for the rest of this paragraph. We consider modular forms over $%
S=\limfunc{Spec}(\mathbb{F})$. The labeling of the embeddings $\overline{%
\sigma }_{i}^{(j)}$ for $1\leq j\leq r$ and $0\leq i\leq f_{j}-1$\ induces a
canonical splitting:%
\begin{eqnarray*}
\mathbb{G}_{\mathbb{F}} &=&\dbigoplus\nolimits_{j=1}^{r}\left( \limfunc{Res}%
\nolimits_{\mathbb{F}_{\mathfrak{P}_{j}}/\mathbb{F}_{p}}(\mathbb{G}_{m,%
\mathbb{F}_{\mathfrak{P}_{j}}})\times _{\limfunc{Spec}\mathbb{F}_{p}}%
\limfunc{Spec}\mathbb{F}\right) \\
&=&\dbigoplus\nolimits_{j=1}^{r}\dbigoplus\nolimits_{\overline{\sigma }%
_{i}^{(j)}:\mathbb{F}_{\mathfrak{P}_{j}}\hookrightarrow \mathbb{F}}\mathbb{G}%
_{m,\mathbb{F}},
\end{eqnarray*}%
such that the projection $\chi _{(j,i)}$ of $\mathbb{G}_{\mathbb{F}}$ onto
the $(j,i)$th factor is induced by $\overline{\sigma }_{i}^{(j)}$. The
character group $\mathbb{X}_{\mathbb{F}}$ of $\mathbb{G}_{\mathbb{F}}$ is
the free $%
%TCIMACRO{\U{2124} }%
%BeginExpansion
\mathbb{Z}
%EndExpansion
$-module or rank $g$ generated by these projections. A geometric Hilbert
modular form over $\limfunc{Spec}(\mathbb{F})$ whose weight is $%
\tprod\nolimits_{j=1}^{r}\tprod\nolimits_{i=0}^{f_{j}-1}\mathbb{\chi }%
_{(j,i)}^{a_{i}^{(j)}}$ for some $a_{i}^{(j)}\in 
%TCIMACRO{\U{2124} }%
%BeginExpansion
\mathbb{Z}
%EndExpansion
$ is also said to have weight vector $\vec{a}=(\vec{a}^{(1)},...,\vec{a}%
^{(r)})$ where $\vec{a}^{(j)}=(a_{0}^{(j)},...,a_{f_{j}-1}^{(j)})$ for $%
1\leq j\leq r$.

Theorem 2.1 of \cite{Go} shows that, for any $1\leq j\leq r$ and $0\leq
i\leq f_{j}-1$, there is an $(\mathfrak{a},\mathfrak{a}^{+})$-polarized
Hilbert modular form $h_{(j,i)}$ over $\limfunc{Spec}(\mathbb{F})$ having
weight $\mathbb{\chi }_{(j,i-1)}^{p}\mathbb{\chi }_{(j,i)}^{-1}$ and level $%
1 $, whose $q$-expansion at every $(\mathfrak{a},\mathfrak{a}^{+})$%
-polarized unramified $\mathbb{F}_{p}$-rational cusp is one. $h_{(j,i)}$ is
called the $(j,i)$th \textit{partial Hasse invariant}. As mentioned earlier,
the forms $h_{(j,i)}$ are not liftable to characteristic zero; even the
total Hasse invariant, i.e., the form $h=\tprod\nolimits_{(j,i)}h_{(j,i)}$,
having parallel weight $(p-1,p-1,...,p-1)$, is not always liftable to
characteristic zero (cf. Proposition 3.1 in \cite{Go}).

As a consequence of the existence of the partial Hasse invariants, one can
produce (geometric) weight shiftings. More precisely, fix an integer $j$
such that $1\leq j\leq r$ and assume $\chi \in \mathbb{X}_{\mathbb{F}}$ is
such that $M_{\chi }(\mu _{N},\mathbb{F})\neq 0$; denote the weight vector
associated to $\chi $ by $\vec{a}=(\vec{a}^{(1)},...,\vec{a}^{(r)})$.
Multiplication by $h_{(j,i)}$ for an integer $i$ such that $0\leq i\leq
f_{j}-1$ induces a Hecke injection (\textit{weight shifting}) of $M_{\chi
}(\mu _{N},\mathbb{F})$ into $M_{\chi ^{\prime }}(\mu _{N},\mathbb{F})$,
where the weight vector associated to $\chi ^{\prime }$ is $\vec{a}+\vec{t}$
and $\vec{t}=(\vec{t}^{(1)},...,\vec{t}^{(r)})$ is such that $\vec{t}^{(r)}=%
\vec{0}$ if $r\neq j$, while $\vec{t}^{(j)}$ is one of the following $f_{j}$%
-tuples: 
\begin{eqnarray*}
&&(-1,0,0,...,0,p)\text{ \ \ if }i=0, \\
&&(p,-1,0,...,0,0)\text{ \ \ if }i=1, \\
&&(0,p,-1,...,0,0)\text{ \ \ if }i=2, \\
&&... \\
&&(0,0,0,...,p,-1)\text{ \ \ if }i=f_{j}-1.
\end{eqnarray*}

\noindent In this case, we will say that $h_{(j,i)}$ induces a weight
shifting by $\vec{t}$.

In \cite{Ka} 2.5. and \cite{AG} 12, generalized theta operators acting on
spaces of geometric Hilbert modular forms over $\limfunc{Spec}(\mathbb{F})$
are defined, allowing additional weight shiftings. For example, if $p$ is
inert in $F/%
%TCIMACRO{\U{211a} }%
%BeginExpansion
\mathbb{Q}
%EndExpansion
$, these operators induce shiftings by the vectors:%
\begin{eqnarray*}
&&(1,0,0,...,0,p), \\
&&(p,1,0,...,0,0), \\
&&(0,p,1,...,0,0), \\
&&... \\
&&(0,0,0,...,p,1).
\end{eqnarray*}

The reader will notice that the two sets of weight shifting vectors
described above are contained in the sets of weight shifting vectors
produced in \ref{GDI} and \ref{GDO} for $\mathbb{\bar{F}}_{p}$%
-representation of $GL_{2}(\mathbb{F})$. Exploiting the adelic definition of
Hilbert modular forms, we will see that all the geometric weight shiftings
can be obtained as cohomological weight shiftings via the operators
considered in Section 3. The purely cohomological picture will be reacher,
as more shiftings will be allowed. The formation of spaces of adelic
automorphic forms on definite quaternion algebra will have the big advantage
of being compatible with base changes, under suitable assumptions
(Proposition \ref{onto}). Finally, our cohomological weight shiftings
translate into weight shiftings for ($\func{mod}p$) Galois representations
arising from automorphic forms on $GL_{2}(\mathbb{A}_{F})$.

\subsection{Automorphic forms on definite quaternion algebras\label{adelic
HMF}}

We recall the definition and some properties of automorphic forms on
definite quaternion algebras over totally real number fields. The exposition
follows \cite{Taylor} and \cite{Ki}; cf. also \cite{Ta}.

Fix a finite set $\Sigma \subset \mathfrak{M}_{F,f}$ that is disjoint from
the set of places of $F$ lying above $p$ and such that $\#\Sigma +[F:%
%TCIMACRO{\U{211a} }%
%BeginExpansion
\mathbb{Q}
%EndExpansion
]\equiv 0(\func{mod}2)$. Let $D$ be a quaternion algebra over $F$ whose
ramification set is $\mathfrak{M}_{F,\infty }\cup \Sigma $. Let $\mathcal{O}%
_{D}$ be a fixed maximal order of $D$ and for any $v\in \mathfrak{M}%
_{F,f}-\Sigma $ fix ring isomorphisms $\left( \mathcal{O}_{D}\right)
_{v}\simeq M_{2}(\mathcal{O}_{F_{v}})$.

Let $U$ be a compact open subgroup of $\left( D\otimes _{F}\mathbb{A}%
_{F}^{\infty }\right) ^{\times }$ such that:

\begin{enumerate}
\item $U=\tprod\nolimits_{v\in \mathfrak{M}_{F,f}}U_{v}$, where $U_{v}$ is a
subgroup of $\left( \mathcal{O}_{D}\right) _{v}^{\times }$;

\item $U_{v}=\left( \mathcal{O}_{D}\right) _{v}^{\times }$ if $v\in \Sigma $;

\item if $v|p$, then $U_{v}=GL_{2}(\mathcal{O}_{F_{v}}).$
\end{enumerate}

Let $A$ be a topological $%
%TCIMACRO{\U{2124} }%
%BeginExpansion
\mathbb{Z}
%EndExpansion
_{p}$-algebra. Let $v$ be a place of $F$ above $p$, say $v=v_{j}:=\mathfrak{P%
}_{j}$ for some integer $j$ such that $1\leq j\leq r$; let $W_{\tau _{j}}$
be a free $A$-module of finite rank and fix a continuous homomorphism 
\begin{equation*}
\tau _{j}:U_{v_{j}}=GL_{2}(\mathcal{O}_{F_{\mathfrak{P}_{j}}})%
\longrightarrow \limfunc{Aut}(W_{\tau _{j}}),
\end{equation*}

\noindent where $\limfunc{Aut}(W_{\tau _{j}})$ is the group of continuous $A$%
-linear automorphisms of $W_{\tau _{j}}$. \noindent Let $W_{\tau
}=\tbigotimes\nolimits_{j=1}^{r}W_{\tau _{j}},$ where the tensor products
are over $A$, and denote by $\tau $ the corresponding group homomorphism $%
\tau :\tprod\nolimits_{j=1}^{r}U_{v_{j}}\mathbb{\rightarrow }\limfunc{Aut}%
(W_{\tau })$. If no confusion arises, we also denote by $\tau $ the action
of $U$ on $W_{\tau }$ induced by precomposing the latter morphism with the
natural projection $U\mathbb{\rightarrow }\tprod\nolimits_{j=1}^{r}U_{v_{j}}$%
.

For $A$ as above, let $\psi :\left( \mathbb{A}_{F}^{\infty }\right) ^{\times
}/F^{\times }\mathbb{\rightarrow }A^{\times }$ be a continuous character
such that, for any $v\in \mathfrak{M}_{F,f}$:%
\begin{equation*}
\tau _{|U_{v}\cap \mathcal{O}_{F_{v}}^{\times }}(u)=\psi ^{-1}(u)\cdot
Id_{W_{\tau }}\text{, \ \ \ for all }u\in U_{v}\cap \mathcal{O}%
_{F_{v}}^{\times }\text{.}
\end{equation*}

\noindent We say that such a Hecke character $\psi $ is \textit{compatible}
with $\tau $.

\begin{definition}
\label{defhmf}For $D,U,A,\tau ,W_{\tau }$ and $\psi $ as above, the space $%
S_{\tau ,\psi }(U,A)$ of automorphic forms on $D$ having level $U$, weight $%
\tau $, character $\psi $ and coefficients in $A$ is the $A$-module
consisting of all the functions:%
\begin{equation*}
f:D^{\times }\backslash \left( D\otimes _{F}\mathbb{A}_{F}^{\infty }\right)
^{\times }\longrightarrow W_{\tau }
\end{equation*}

satisfying:

\begin{description}
\item[(a)] $f(gu)=\tau (u)^{-1}f(g)$ for all $g\in \left( D\otimes _{F}%
\mathbb{A}_{F}^{\infty }\right) ^{\times }$ and all $u\in U;$

\item[(b)] $f(gz)=\psi (z)f(g)$ for all $g\in \left( D\otimes _{F}\mathbb{A}%
_{F}^{\infty }\right) ^{\times }$ and all $z\in \left( \mathbb{A}%
_{F}^{\infty }\right) ^{\times }$.
\end{description}
\end{definition}

As in \cite{Ki}, we will always assume, unless otherwise stated, that for
all $t\in \left( D\otimes _{F}\mathbb{A}_{F}^{\infty }\right) ^{\times }$,
the finite group $(U\cdot \left( \mathbb{A}_{F}^{\infty }\right) ^{\times
}\cap t^{-1}D^{\times }t)/F^{\times }$ has order prime to $p$. This
assumption is automatically satisfied if $U$ is sufficiently small, as Lemma
1.1. of \cite{Taylor} implies that in this case $(U\cdot \left( \mathbb{A}%
_{F}^{\infty }\right) ^{\times }\cap t^{-1}D^{\times }t)/F^{\times }$ is a $%
2 $-group. We obtain as a consequence (cf. \cite{Taylor}, Corollary 1.2):

\begin{proposition}
\label{onto}Let $B$ a topological $A$-algebra. Then the natural morphism $%
S_{\tau ,\psi }(U,A)\otimes _{A}B\mathbb{\rightarrow }S_{\tau \otimes
_{A}B,\psi \otimes _{A}B}(U,B)$ is an isomorphism of $B$-modules.
\end{proposition}

Define a left action of $\left( D\otimes _{F}\mathbb{A}_{F}^{\infty }\right)
^{\times }$ on the set of functions $D^{\times }\backslash \left( D\otimes
_{F}\mathbb{A}_{F}^{\infty }\right) ^{\times }\rightarrow W_{\tau }$ by
setting $(gf)(x):=f(xg)$ for all $g,x\in \left( D\otimes _{F}\mathbb{A}%
_{F}^{\infty }\right) ^{\times }$. Let $S$ be a set of primes\ of $F$\
containing the ramification set of $D$, the primes above $p$ and the primes $%
v$ for which $U_{v}$ is not a maximal compact subgroup of $D_{v}^{\times }$.
Let $\mathbb{T}_{S,A}^{univ}=A[T_{v},S_{v}:v\notin S]$ be the commutative
polynomial $A$-algebra in the indicated indeterminates. For each finite
place $v\notin S$, let $\varpi _{v}$ be a fixed uniformizer for $F_{v}$. $%
S_{\tau ,\psi }(U,A)$ has a natural action of $\mathbb{T}_{S,A}^{univ}$,
with $S_{v}$ acting via the double coset $U\left( \QTATOP{\varpi _{v}}{{}}%
\QTATOP{{}}{\varpi _{v}}\right) U$ and $T_{v}$ via $U\left( \QTATOP{\varpi
_{v}}{{}}\QTATOP{{}}{1}\right) U$ (cf. \cite{Taylor}, 1); this action does
not depend upon the choices of uniformizers that we made. The image of $%
\mathbb{T}_{S,A}^{univ}$ in the ring of $A$-module endomorphisms of $S_{\tau
,\psi }(U,A)$ is the Hecke algebra $\mathbb{T}_{S,A}$ acting on $S_{\tau
,\psi }(U,A)$. The isomorphism of Proposition \ref{onto} is Hecke
equivariant.

\subsection{Behavior of Hecke eigensystems under reduction modulo $\mathfrak{%
M}_{R}$}

For a discrete valuation ring $R$, we will denote by $\mathfrak{M}_{R}$ its
maximal ideal. If the residual characteristic of $R$ is $p>0$ and no
confusion arises, we will also improperly refer to reduction modulo $%
\mathfrak{M}_{R}$ as reduction modulo $p$. If $\mathcal{T}$ is a commutative
algebra, a system of eigenvalues of $\mathcal{T}$ with values in $R$ is a
set theoretic map $\Omega :\mathcal{T}\rightarrow R$; the reduction of $%
\Omega $ modulo $p$, denoted $\bar{\Omega}$, is the function obtained by
composing $\Omega $ with the reduction morphism $R\rightarrow \frac{R}{%
\mathfrak{M}_{R}}$. Let $R\mathcal{T}=R\otimes _{%
%TCIMACRO{\U{2124} }%
%BeginExpansion
\mathbb{Z}
%EndExpansion
}\mathcal{T}$; if $M$ is an $R\mathcal{T}$-module, we say that a system of
eigenvalues $\Omega :\mathcal{T}\mathbb{\rightarrow }R$ occurs in $M$ if
there is a non-zero element $m\in M$ such that $Tm=\Omega (T)m$ for all $%
T\in \mathcal{T}$. Such a non-zero $m$ is called an $\Omega $-eigenvector.

Fixing $R$ and $\mathcal{T}$ as above. We have:

\begin{lemma}
\label{AS1}Let $M$ be an $R\mathcal{T}$-module which is finitely generated
over $R$. If $\Omega :\mathcal{T}\mathbb{\rightarrow }R$ is a system of
eigenvalues of $\mathcal{T}$ occurring in $M$, then $\bar{\Omega}:\mathcal{T}%
\mathbb{\rightarrow }\frac{R}{\mathfrak{M}_{R}}$ is a system of eigenvalues
of $\mathcal{T}$ occurring in $\bar{M}:=M\otimes _{R}\frac{R}{\mathfrak{M}%
_{R}}.$
\end{lemma}

\textbf{Proof }Cf. \cite{AS1}, Proposition 1.2.3. $\blacksquare $

\bigskip

\begin{lemma}
\label{AS2}Let $M$ be an $R\mathcal{T}$-module which is finite and free over 
$R$. Let $\bar{\Omega}:\mathcal{T}\mathbb{\rightarrow }\frac{R}{\mathfrak{M}%
_{R}}$ be a system of eigenvalues of $\mathcal{T}$ occurring in $\bar{M}%
=M\otimes _{R}\frac{R}{\mathfrak{M}_{R}}$. There exists a finite extension
of discrete valuation rings $R^{\prime }/R$ such that $\mathfrak{M}%
_{R^{\prime }}\cap R=\mathfrak{M}_{R}$\ and a system of eigenvalues $\Omega
^{\prime }:\mathcal{T}\rightarrow R^{\prime }$ of $\mathcal{T}$ occurring in 
$M\otimes _{R}R^{\prime }$ such that, for all $T\in \mathcal{T}$, $\Omega
^{\prime }(T)(\func{mod}\mathfrak{M}_{R^{\prime }})=\bar{\Omega}(T)$ in $%
\frac{R^{\prime }}{\mathfrak{M}_{R^{\prime }}}.$ (Here we view $\frac{R}{%
\mathfrak{M}_{R}}\subseteq \frac{R^{\prime }}{\mathfrak{M}_{R^{\prime }}}$
by the given embedding $R\subseteq R^{\prime }$).
\end{lemma}

\textbf{Proof }Cf. \cite{DS}, Lemme 6.11. A generalization of the result is
given in \cite{AS1}, Proposition 1.2.2. $\blacksquare $

\bigskip

Let $D$, $U$, $\tau $, $W_{\tau }$ and $\psi $ be as in \ref{adelic HMF},
and set $A=\mathcal{O}$. In particular, we assume that $\psi $ is compatible
with $\left( \tau ,W_{\tau }\right) $, $U$ is small enough and $p$ is odd.
Denote by a bar the operation of tensoring over $\mathcal{O}$ with $\mathbb{F%
}$. From now on, unless otherwise stated, we assume fixed a set $S$ of
primes\ of $F$\ containing the ramification set of $D$, the primes above $p$
and the primes $v$ for which $U_{v}$ is not a maximal compact subgroup of $%
D_{v}^{\times }$. The Hecke eigensystems considered below will always be
with respect to the Hecke algebra $\mathbb{T}_{S,A^{\prime }}^{univ}$ for
some topological $%
%TCIMACRO{\U{2124} }%
%BeginExpansion
\mathbb{Z}
%EndExpansion
_{p}$-algebra $A^{\prime }$.

\begin{proposition}
\label{shift}Fix an $\mathcal{O}$-valued weight $\left( \tau ^{\prime
},W_{\tau ^{\prime }}\right) $ together with a compatible Hecke character $%
\psi ^{\prime }:\left( \mathbb{A}_{F}^{\infty }\right) ^{\times }/F^{\times }%
\mathbb{\rightarrow }\mathcal{O}^{\times }$ such that $\bar{\psi}^{\prime }=%
\bar{\psi}$. Let $\varphi :\left( \bar{\tau},W_{\bar{\tau}}\right) \mathbb{%
\rightarrow }\left( \bar{\tau}^{\prime },W_{\bar{\tau}^{\prime }}\right) $
be a non-zero intertwining operator for $\mathbb{F}$-representations of $U$. 
$\varphi $ induces a Hecke equivariant map $\varphi _{\ast }:S_{\bar{\tau},%
\bar{\psi}}(U,\mathbb{F})\mathbb{\rightarrow }S_{\bar{\tau}^{\prime },\bar{%
\psi}}(U,\mathbb{F})$.

Assume $\varphi $ is injective: then if $\Omega $ is a Hecke eigensystem
occurring in $S_{\tau ,\psi }(U,\mathcal{O})$, there is a finite extension
of $E$, with ring of integer $\mathcal{O}^{\prime }$ such that $\mathfrak{M}%
_{\mathcal{O}^{\prime }}\cap \mathcal{O}=\mathfrak{M}_{\mathcal{O}}$,\ and
there is a Hecke eigensystem $\Omega ^{\prime }$ occurring in $S_{\tau
^{\prime },\psi ^{\prime }}(U,\mathcal{O}^{\prime })$ such that:%
\begin{equation*}
\Omega ^{\prime }(\func{mod}\mathfrak{M}_{\mathcal{O}^{\prime }})=\Omega (%
\func{mod}\mathfrak{M}_{\mathcal{O}})\text{ \ \ in }\frac{\mathcal{O}%
^{\prime }}{\mathfrak{M}_{\mathcal{O}^{\prime }}}.
\end{equation*}
\end{proposition}

\textbf{Proof }For $f\in S_{\bar{\tau},\bar{\psi}}(U,\mathbb{F})$ set $%
\varphi _{\ast }(f):=\varphi \circ f$. If $g\in \left( D\otimes _{F}\mathbb{A%
}_{F}^{\infty }\right) ^{\times },u\in U$ and $z\in \left( \mathbb{A}%
_{F}^{\infty }\right) ^{\times }$ we have:

\begin{eqnarray*}
\varphi _{\ast }(f)(gu) &=&\varphi (f(gu))=\varphi \left( \bar{\tau}%
(u^{-1})f(g)\right) =\bar{\tau}^{\prime }(u)^{-1}\varphi \left( f(g)\right) ,
\\
\varphi _{\ast }(f)(gz) &=&\varphi (f(gz))=\varphi \left( \bar{\psi}%
(z)f(g)\right) =\bar{\psi}(z)\varphi \left( f(g)\right) .
\end{eqnarray*}

\noindent Since $\bar{\psi}^{\prime }=\bar{\psi}$, we have that $\bar{\tau}%
^{\prime }$ and $\bar{\psi}$ are compatible and we conclude that $\varphi
_{\ast }(f)\in S_{\bar{\tau}^{\prime },\bar{\psi}}(U,\mathbb{F})$. If $%
g,x\in \left( D\otimes _{F}\mathbb{A}_{F}^{\infty }\right) ^{\times }$, we
have:%
\begin{eqnarray*}
\left( g\cdot \varphi _{\ast }(f)\right) \left( x\right) &=&\left( \varphi
\circ f\right) \left( xg\right) \\
&=&\left( \varphi \circ (g\cdot f)\right) \left( x\right) \\
&=&\left( \varphi _{\ast }(g\cdot f)\right) \left( x\right) ,
\end{eqnarray*}

\noindent so that $\varphi $ is Hecke-equivariant. Assume now that $\varphi $
is injective and notice that this implies the injectivity of $\varphi _{\ast
}$. Let $\Omega $ be a Hecke eigensystem occurring in the finite $\mathcal{O}
$-module with Hecke action $S_{\tau ,\psi }(U,\mathcal{O})$; by Proposition %
\ref{onto}, reduction modulo $p$ induces a Hecke equivariant surjection:%
\begin{equation*}
\pi :S_{\tau ,\psi }(U,\mathcal{O})\longrightarrow S_{\bar{\tau},\bar{\psi}%
}(U,\mathbb{F}).
\end{equation*}

By Lemma \ref{AS1}, the Hecke eigensystem $\bar{\Omega}:=\Omega (\func{mod}%
\mathfrak{M}_{\mathcal{O}})$ occurs in $S_{\bar{\tau},\bar{\psi}}(U,\mathbb{F%
})$, and hence in $S_{\bar{\tau}^{\prime },\bar{\psi}}(U,\mathbb{F})$ as $%
\varphi _{\ast }$ is Hecke equivariant and injective. Now, applying Lemma %
\ref{AS2} to the Hecke equivariant surjection $S_{\tau ^{\prime },\psi
^{\prime }}(U,\mathcal{O})\mathbb{\rightarrow }S_{\bar{\tau}^{\prime },\bar{%
\psi}}(U,\mathbb{F})$, we deduce the existence of a finite extension of
discrete valuation rings $\mathcal{O}^{\prime }/\mathcal{O}$ such that $%
\mathfrak{M}_{\mathcal{O}^{\prime }}\cap \mathcal{O}=\mathfrak{M}_{\mathcal{O%
}}$, and of a Hecke eigensystem $\Omega ^{\prime }:\mathbb{T}_{S,\mathcal{O}%
^{\prime }}\rightarrow \mathcal{O}^{\prime }$ occurring in $S_{\tau ^{\prime
},\psi ^{\prime }}(U,\mathcal{O})\otimes _{\mathcal{O}}\mathcal{O}^{\prime }$
whose reduction modulo $\mathfrak{M}_{\mathcal{O}^{\prime }}$ has value in $%
\mathbb{F}\subset \frac{\mathcal{O}^{\prime }}{\mathfrak{M}_{\mathcal{O}%
^{\prime }}}$ and coincide with $\bar{\Omega}$. By Proposition \ref{onto}, $%
S_{\tau ^{\prime },\psi ^{\prime }}(U,\mathcal{O})\otimes _{\mathcal{O}}%
\mathcal{O}^{\prime }\simeq S_{\tau ^{\prime },\psi ^{\prime }}(U,\mathcal{O}%
^{\prime })$ as Hecke modules, and we are done. $\blacksquare $

\subsection{Holomorphic weights\label{conv_on_emb}}

For any integer $j$ such that $1\leq j\leq r$ let us fix two tuples $\vec{k}%
^{(j)}=(k_{0}^{(j)},...,k_{f_{j}-1}^{(j)})\in 
%TCIMACRO{\U{2124} }%
%BeginExpansion
\mathbb{Z}
%EndExpansion
_{\geq 2}^{f_{j}}$ and $\vec{w}^{(j)}=(w_{0}^{(j)},...,w_{f_{j}-1}^{(j)})\in 
%TCIMACRO{\U{2124} }%
%BeginExpansion
\mathbb{Z}
%EndExpansion
^{f_{j}}$. Define the finite free $\mathcal{O}$-module with $GL_{2}(\mathcal{%
O})$-action: 
\begin{equation*}
W_{(\vec{k}^{(j)},\vec{w}^{(j)})}:=\dbigotimes\nolimits_{i=0}^{f_{j}-1}%
\limfunc{Sym}\nolimits^{k_{i}^{(j)}-2}\mathcal{O}^{2}\otimes
\det\nolimits^{w_{i}^{(j)}}
\end{equation*}%
where the tensor products are over $\mathcal{O}$.

If we let the group $GL_{2}(\mathcal{O}_{F_{\mathfrak{P}_{j}}})$\ act on the
tensor factor $\limfunc{Sym}\nolimits^{k_{i}^{(j)}-2}\mathcal{O}^{2}\otimes
\det\nolimits^{w_{i}^{(j)}}$($0\leq i\leq f_{j}-1$) via the embedding $%
GL_{2}(\mathcal{O}_{F_{\mathfrak{P}_{j}}})\mathbb{\rightarrow }GL_{2}(%
\mathcal{O})$ induced by $\sigma _{i}^{(j)}=\sigma ^{i}\circ \sigma
_{0}^{(j)}$, $W_{(\vec{k}^{(j)},\vec{w}^{(j)})}$ can be seen as a
representation of $GL_{2}(\mathcal{O}_{F_{\mathfrak{P}_{j}}})$. We convene
of viewing $GL_{2}(\mathcal{O}_{F_{\mathfrak{P}_{j}}})$ as a subgroup of $%
GL_{2}(\mathcal{O})$ via the embedding $\sigma _{0}^{(j)}$, and we write the 
$GL_{2}(\mathcal{O}_{F_{\mathfrak{P}_{j}}})$-representation $W_{(\vec{k}%
^{(j)},\vec{w}^{(j)})}$ as:%
\begin{equation*}
W_{(\vec{k}^{(j)},\vec{w}^{(j)})}=\dbigotimes\nolimits_{i=0}^{f_{j}-1}\left( 
\limfunc{Sym}\nolimits^{k_{i}^{(j)}-2}\mathcal{O}^{2}\otimes
\det\nolimits^{w_{i}^{(j)}}\right) ^{[i]},
\end{equation*}%
where the superscript $[i]$ indicates twisting by the $i$th power of the
Frobenius element $\sigma $. In the sequel, unless otherwise stated, we
always view $GL_{2}(\mathcal{O}_{F_{\mathfrak{P}_{j}}})\subseteq GL_{2}(%
\mathcal{O})$ via $\sigma _{0}^{(j)}$.

Denote by $\tau _{(\vec{k}^{(j)},\vec{w}^{(j)})}$ the continuous action of $%
GL_{2}(\mathcal{O}_{F_{\mathfrak{P}_{j}}})$ on $W_{(\vec{k}^{(j)},\vec{w}%
^{(j)})}$ and let $\tau _{(\vec{k},\vec{w})}=\tbigotimes\nolimits_{j=1}^{r}%
\tau _{(\vec{k}^{(j)},\vec{w}^{(j)})}$, where the tensor products are over $%
\mathcal{O}$ and $\vec{k}=(\vec{k}^{(1)},...,\vec{k}^{(r)})$. We have: 
\begin{equation*}
\tau _{(\vec{k},\vec{w})}:\dprod\nolimits_{j=1}^{r}GL_{2}(\mathcal{O}_{F_{%
\mathfrak{P}_{j}}})\longrightarrow \limfunc{Aut}W_{(\vec{k},\vec{w})},
\end{equation*}

\noindent with $W_{(\vec{k},\vec{w})}=\tbigotimes\nolimits_{j=1}^{r}W_{(\vec{%
k}^{(j)},\vec{w}^{(j)})}$ (tensor product over $\mathcal{O}$).

If there is some integer $j$ such that $\vec{k}^{(j)}=(2,...,2)$, we say
that the weight $\tau _{(\vec{k},\vec{w})}$ contains a $(2,...,2)$\textit{%
-block relative to the prime }$\mathfrak{P}_{j}$. This terminology is not
standard but it is used throughout the paper.

We say that $\tau _{(\vec{k},\vec{w})}$ is a \textit{holomorphic weight} if
there exists an integer $w$ such that:%
\begin{equation}
k_{i}^{(j)}+2w_{i}^{(j)}-1=w  \tag{*}  \label{star}
\end{equation}%
for all $1\leq j\leq r$ and all $0\leq i\leq f_{j}-1$ (cf. \cite{H}).

The pair $(\vec{k},\vec{w})\in 
%TCIMACRO{\U{2124} }%
%BeginExpansion
\mathbb{Z}
%EndExpansion
_{\geq 2}^{g}\times 
%TCIMACRO{\U{2124} }%
%BeginExpansion
\mathbb{Z}
%EndExpansion
^{g}$ is called the parameter pair for $\tau _{(\vec{k},\vec{w})}$. If $\tau
_{(\vec{k},\vec{w})}$ is a holomorphic weight, it is also determined by the
parameter pair $(\vec{k},w)\in 
%TCIMACRO{\U{2124} }%
%BeginExpansion
\mathbb{Z}
%EndExpansion
_{\geq 2}^{g}\times 
%TCIMACRO{\U{2124} }%
%BeginExpansion
\mathbb{Z}
%EndExpansion
$, with $w$ as in (\ref{star}).

\subsubsection{Some results on holomorphic weight shiftings}

\begin{lemma}
\label{compat}Let us view the holomorphic weight $\tau _{(\vec{k},w)}$ as an 
$\mathcal{O}$-representation of the fixed level $U\subset \left( D\otimes
_{F}\mathbb{A}_{F}^{\infty }\right) ^{\times }$. A Hecke character $\psi
:\left( \mathbb{A}_{F}^{\infty }\right) ^{\times }/F^{\times }\mathbb{%
\rightarrow }\mathcal{O}^{\times }$ is compatible with $\tau _{(\vec{k},w)}$
if and only if the following two conditions are satisfied:

\begin{description}
\item[(a)] $\psi (u)=1$ for all $u\in U_{v}\cap \mathcal{O}_{F_{v}}^{\times
} $, where $v\in \mathfrak{M}_{F,f}$ and $v\NEG{|}p;$

\item[(b)] $\psi (u)=\left( \limfunc{Nm}_{F_{\mathfrak{P}_{j}}/%
%TCIMACRO{\U{211a} }%
%BeginExpansion
\mathbb{Q}
%EndExpansion
_{p}}(u)\right) ^{1-w}$ for all $u\in \mathcal{O}_{F_{\mathfrak{P}%
_{j}}}^{\times }$, where $1\leq j\leq r$.
\end{description}
\end{lemma}

\textbf{Proof }The reason for condition (a) is clear, as the representation $%
\tau _{(\vec{k},w)}$ factors through $\tprod\nolimits_{j=1}^{r}GL_{2}(%
\mathcal{O}_{F_{\mathfrak{P}_{j}}})$. Let $j$ be such that $1\leq j\leq r$
and fix $u\in \mathcal{O}_{F_{\mathfrak{P}_{j}}}^{\times }$; recall that we
embed $\mathcal{O}_{F_{\mathfrak{P}_{j}}}$ in $\mathcal{O}$\ via $\sigma
_{0}^{(j)}$. The matrix $\left( \QTATOP{u}{{}}\QTATOP{{}}{u}\right) \in
GL_{2}(\mathcal{O}_{F_{\mathfrak{P}_{j}}})$ acts on $W_{(\vec{k}^{(j)},w)}$
as the automorphism:%
\begin{eqnarray*}
\dbigotimes\nolimits_{i=0}^{f_{j}-1}\left( \sigma
^{i}(u)^{k_{i}^{(j)}-2+2w_{i}^{(j)}}\cdot Id_{i}\right)
&=&\dbigotimes\nolimits_{i=0}^{f_{j}-1}\sigma ^{i}(u)^{w-1}\cdot Id_{i} \\
&=&\left( \limfunc{Nm}\nolimits_{F_{\mathfrak{P}_{j}}/%
%TCIMACRO{\U{211a} }%
%BeginExpansion
\mathbb{Q}
%EndExpansion
_{p}}(u)\right) ^{w-1}\cdot \dbigotimes\nolimits_{i=0}^{f_{j}-1}Id_{i} \\
&=&\left( \limfunc{Nm}\nolimits_{F_{\mathfrak{P}_{j}}/%
%TCIMACRO{\U{211a} }%
%BeginExpansion
\mathbb{Q}
%EndExpansion
_{p}}(u)\right) ^{w-1}\cdot Id_{W_{(\vec{k}^{(j)},w)}},
\end{eqnarray*}

\noindent where $Id_{i}$ denotes the identity map of the $\mathcal{O}$%
-vector space: 
\begin{equation*}
\left( \limfunc{Sym}\nolimits^{k_{i}^{(j)}-2}\mathcal{O}^{2}\otimes
\det\nolimits^{w_{i}^{(j)}}\right) ^{[i]},
\end{equation*}

\noindent and we used the assumption that the local extension $F_{\mathfrak{P%
}_{j}}/%
%TCIMACRO{\U{211a} }%
%BeginExpansion
\mathbb{Q}
%EndExpansion
_{p}$ is unramified with Galois group generated by the restriction of $%
\sigma $ to $F_{\mathfrak{P}_{j}}$. The result now follows, as we need to
have $\tau _{|\mathcal{O}_{F_{\mathfrak{P}_{j}}}^{\times }}(u)=\psi
^{-1}(u)\cdot Id_{W_{(\vec{k},w)}}$. $\blacksquare $

\bigskip

\begin{lemma}
\label{exist}Let $w$ be an even integer. Then there exists a continuous
character $\psi :\left( \mathbb{A}_{F}^{\infty }\right) ^{\times }/F^{\times
}\mathbb{\rightarrow }\mathcal{%
%TCIMACRO{\U{2124} }%
%BeginExpansion
\mathbb{Z}
%EndExpansion
}_{p}^{\times }$ such that:

\begin{description}
\item[(a)] $\psi (u)=1$ for all $u\in \mathcal{O}_{F_{v}}^{\times }$, where $%
v\in \mathfrak{M}_{F,f}$ and $v\NEG{|}p;$

\item[(b)] $\psi (u)=\left( \limfunc{Nm}_{F_{\mathfrak{P}_{j}}/%
%TCIMACRO{\U{211a} }%
%BeginExpansion
\mathbb{Q}
%EndExpansion
_{p}}(u)\right) ^{w}$ for all $u\in \mathcal{O}_{F_{\mathfrak{P}%
_{j}}}^{\times }$, where $1\leq j\leq r$.
\end{description}
\end{lemma}

\textbf{Proof }The ad\`{e}les norm map $\left( \mathbb{A}_{F}^{\infty
}\right) ^{\times }\mathbb{\rightarrow }\left( \mathbb{A}_{%
%TCIMACRO{\U{211a} }%
%BeginExpansion
\mathbb{Q}
%EndExpansion
}^{\infty }\right) ^{\times }$ induces a continuous homomorphism $\limfunc{Nm%
}:\left( \mathbb{A}_{F}^{\infty }\right) ^{\times }/F^{\times }\mathbb{%
\rightarrow }\left( \mathbb{A}_{%
%TCIMACRO{\U{211a} }%
%BeginExpansion
\mathbb{Q}
%EndExpansion
}^{\infty }\right) ^{\times }/%
%TCIMACRO{\U{211a} }%
%BeginExpansion
\mathbb{Q}
%EndExpansion
^{\times }$. The group-theoretic decomposition $\left( \mathbb{A}_{%
%TCIMACRO{\U{211a} }%
%BeginExpansion
\mathbb{Q}
%EndExpansion
}^{\infty }\right) ^{\times }=%
%TCIMACRO{\U{211a} }%
%BeginExpansion
\mathbb{Q}
%EndExpansion
^{\times }\cdot \hat{%
%TCIMACRO{\U{2124}}%
%BeginExpansion
\mathbb{Z}%
%EndExpansion
}^{\times }$ induces a continuous isomorphism $\beta :\left( \mathbb{A}_{%
%TCIMACRO{\U{211a} }%
%BeginExpansion
\mathbb{Q}
%EndExpansion
}^{\infty }\right) ^{\times }/%
%TCIMACRO{\U{211a} }%
%BeginExpansion
\mathbb{Q}
%EndExpansion
^{\times }\mathbb{\rightarrow }\hat{%
%TCIMACRO{\U{2124}}%
%BeginExpansion
\mathbb{Z}%
%EndExpansion
}^{\times }/\left\langle -1\right\rangle $. Finally, the map $%
\tprod\nolimits_{l}%
%TCIMACRO{\U{2124} }%
%BeginExpansion
\mathbb{Z}
%EndExpansion
_{l}^{\times }\mathbb{\rightarrow }%
%TCIMACRO{\U{2124} }%
%BeginExpansion
\mathbb{Z}
%EndExpansion
_{p}^{\times }$ defined by sending the tuple $\left( a_{l}\right) _{l}\in
\tprod\nolimits_{l}%
%TCIMACRO{\U{2124} }%
%BeginExpansion
\mathbb{Z}
%EndExpansion
_{l}^{\times }$ into $a_{p}^{w}\in 
%TCIMACRO{\U{2124} }%
%BeginExpansion
\mathbb{Z}
%EndExpansion
_{p}^{\times }$ defines a continuous homomorphism $\alpha :\hat{%
%TCIMACRO{\U{2124}}%
%BeginExpansion
\mathbb{Z}%
%EndExpansion
}^{\times }/\left\langle -1\right\rangle \mathbb{\rightarrow }%
%TCIMACRO{\U{2124} }%
%BeginExpansion
\mathbb{Z}
%EndExpansion
_{p}^{\times }$ since $w$ is even. We check that the composition $\psi
:=\alpha \circ \beta \circ \limfunc{Nm}$ is a Hecke character with the
desired properties.

Assume $v=\mathfrak{P}_{j}|p$ and view a fixed $u\in $ $\mathcal{O}_{F_{%
\mathfrak{P}_{j}}}^{\times }$ as an element of $\left( \mathbb{A}%
_{F}^{\infty }\right) ^{\times }$ whose $v$-component is $u$ and whose $%
v^{\prime }$-component is $1$ for all finite places $v^{\prime }\neq v$ of $%
F $. Then $\limfunc{Nm}(u\cdot F^{\times })=\limfunc{Nm}_{F_{\mathfrak{P}%
_{j}}/%
%TCIMACRO{\U{211a} }%
%BeginExpansion
\mathbb{Q}
%EndExpansion
_{p}}(u)\cdot 
%TCIMACRO{\U{211a} }%
%BeginExpansion
\mathbb{Q}
%EndExpansion
^{\times }$, where we identify $\limfunc{Nm}_{F_{\mathfrak{P}_{j}}/%
%TCIMACRO{\U{211a} }%
%BeginExpansion
\mathbb{Q}
%EndExpansion
_{p}}(u)$ with the ad\`{e}le of $%
%TCIMACRO{\U{211a} }%
%BeginExpansion
\mathbb{Q}
%EndExpansion
$ whose $p$-component is the $p$-adic unit $\limfunc{Nm}_{F_{\mathfrak{P}%
_{j}}/%
%TCIMACRO{\U{211a} }%
%BeginExpansion
\mathbb{Q}
%EndExpansion
_{p}}(u)\in 
%TCIMACRO{\U{2124} }%
%BeginExpansion
\mathbb{Z}
%EndExpansion
_{p}^{\times }$ and whose other components are equal to $1$. Then $\left(
\alpha \circ \beta \right) \left( \limfunc{Nm}_{F_{\mathfrak{P}_{j}}/%
%TCIMACRO{\U{211a} }%
%BeginExpansion
\mathbb{Q}
%EndExpansion
_{p}}(u)\cdot 
%TCIMACRO{\U{211a} }%
%BeginExpansion
\mathbb{Q}
%EndExpansion
^{\times }\right) =\left( \limfunc{Nm}_{F_{\mathfrak{P}_{j}}/%
%TCIMACRO{\U{211a} }%
%BeginExpansion
\mathbb{Q}
%EndExpansion
_{p}}(u)\right) ^{w}\in \mathcal{%
%TCIMACRO{\U{2124} }%
%BeginExpansion
\mathbb{Z}
%EndExpansion
}_{p}^{\times }$.

Assume $v$ is a finite place of $F$ lying above some rational prime $l\neq p$
and let $u\in \mathcal{O}_{F_{v}}^{\times }$ viewed as an element of $\left( 
\mathbb{A}_{F}^{\infty }\right) ^{\times }$ in the usual way. Write $%
\limfunc{Nm}(u\cdot F^{\times })=\limfunc{Nm}_{F_{v}/%
%TCIMACRO{\U{211a} }%
%BeginExpansion
\mathbb{Q}
%EndExpansion
_{l}}(u)\cdot 
%TCIMACRO{\U{211a} }%
%BeginExpansion
\mathbb{Q}
%EndExpansion
^{\times }$; since the $p$-component of $\limfunc{Nm}_{F_{v}/%
%TCIMACRO{\U{211a} }%
%BeginExpansion
\mathbb{Q}
%EndExpansion
_{l}}(u)\in \hat{%
%TCIMACRO{\U{2124}}%
%BeginExpansion
\mathbb{Z}%
%EndExpansion
}^{\times }$ is trivial, $\psi (u)=1$.\ $\blacksquare $

\bigskip

Set $A=\mathcal{O}$ and let $D$, $U$, $\left( \tau ,W_{\tau }\right) $ and $%
\psi $ be as in \ref{adelic HMF}.

\begin{proposition}
\label{main}Assume $\tau =\tau _{(\vec{k},w)}$ and $\tau ^{\prime }=\tau _{(%
\vec{k}^{\prime },w^{\prime })}$ are \textit{holomorphic }$\mathcal{O}$%
\textit{-linear weights} for automorphic forms on $D$, with $w\equiv
w^{\prime }(\func{mod}p-1)$ and $w$ odd. Assume that $\tau _{(\vec{k},w)}$
and $\psi $ are compatible and that $\bar{\tau}_{(\vec{k},w)}$ is isomorphic
to an $\mathbb{F}$-linear $U$-subrepresentation of $\bar{\tau}_{(\vec{k}%
^{\prime },w^{\prime })}$. Then:

\begin{description}
\item[(a)] There is a Hecke character $\psi ^{\prime }:\left( \mathbb{A}%
_{F}^{\infty }\right) ^{\times }/F^{\times }\mathbb{\rightarrow }\mathcal{O}%
^{\times }$ which is compatible with $\tau _{(\vec{k}^{\prime },w^{\prime
})} $ and such that $\bar{\psi}^{\prime }=\bar{\psi}$;

\item[(b)] For any Hecke eigensystem $\Omega $\ occurring in $S_{\tau ,\psi
}(U,\mathcal{O})$ there is a finite extension of discrete valuation rings $%
\mathcal{O}^{\prime }/\mathcal{O}$ with $\mathfrak{M}_{\mathcal{O}^{\prime
}}\cap \mathcal{O}=\mathfrak{M}_{\mathcal{O}}$ and a Hecke eigensystem $%
\Omega ^{\prime }$ occurring in $S_{\tau ^{\prime },\psi ^{\prime }}(U,%
\mathcal{O}^{\prime })$ such that $\Omega ^{\prime }(\func{mod}\mathfrak{M}_{%
\mathcal{O}^{\prime }})=\Omega (\func{mod}\mathfrak{M}_{\mathcal{O}}).$
\end{description}
\end{proposition}

\textbf{Proof }Since $p>2$, the integer $1-w^{\prime }$ is even. By Lemma %
\ref{exist}, there exists a Hecke character $\psi ^{\prime \prime }:\left( 
\mathbb{A}_{F}^{\infty }\right) ^{\times }/F^{\times }\mathbb{\rightarrow }%
\mathcal{%
%TCIMACRO{\U{2124} }%
%BeginExpansion
\mathbb{Z}
%EndExpansion
}_{p}^{\times }\subset \mathcal{O}^{\times }$ such that $\psi ^{\prime
\prime }(u)=1$ for all $v\in \mathfrak{M}_{F,f}$ not lying above $p$ and all 
$u\in \mathcal{O}_{F_{v}}^{\times }$, and $\psi ^{\prime \prime }(u)=\left( 
\limfunc{Nm}_{F_{\mathfrak{P}_{j}}/%
%TCIMACRO{\U{211a} }%
%BeginExpansion
\mathbb{Q}
%EndExpansion
_{p}}(u)\right) ^{1-w^{\prime }}$ for $u\in \mathcal{O}_{F_{\mathfrak{P}%
_{j}}}^{\times }$ ($1\leq j\leq r$). By Lemma \ref{compat}, $\psi ^{\prime
\prime }$ is compatible with $\tau _{(\vec{k}^{\prime },w^{\prime })}$.

Let $\alpha $ denote the reduction modulo $\mathfrak{M}_{\mathcal{O}}$ of
the Hecke character $\psi ^{-1}\psi ^{\prime \prime }$. Since $w\equiv
w^{\prime }(\func{mod}p-1)$, by the compatibility of $\psi $ with $\tau _{(%
\vec{k},w)}$ and by the construction of $\psi ^{\prime \prime }$, the
continuous character $\alpha $ is trivial on the open subgroup 
\begin{equation*}
\tprod\nolimits_{v\NEG{|}p}\left( U_{v}\cap \mathcal{O}_{F_{v}}^{\times
}\right) \times \tprod\nolimits_{j=1}^{r}\mathcal{O}_{F_{\mathfrak{P}%
_{j}}}^{\times }
\end{equation*}%
of $\left( \mathcal{O}_{F}\otimes _{%
%TCIMACRO{\U{2124} }%
%BeginExpansion
\mathbb{Z}
%EndExpansion
}\hat{%
%TCIMACRO{\U{2124}}%
%BeginExpansion
\mathbb{Z}%
%EndExpansion
}\right) ^{\times }$. Therefore $\alpha $ factors through a finite discrete
quotient of $\left( \mathbb{A}_{F}^{\infty }\right) ^{\times }$. In
particular, the Teichm\"{u}ller lift $\tilde{\alpha}$ of $\alpha $ is a
continuous character $\left( \mathbb{A}_{F}^{\infty }\right) ^{\times
}/F^{\times }\mathbb{\rightarrow }\mathcal{O}^{\times }$. The $\mathcal{O}%
^{\times }$-valued Hecke character $\psi ^{\prime }:=\psi ^{\prime \prime }%
\tilde{\alpha}^{-1}$ is compatible with $\tau _{(\vec{k}^{\prime },w^{\prime
})}$ and satisfies $\bar{\psi}^{\prime }=\bar{\psi}$, so that (a) is proved.

Part (b) follows by applying Proposition \ref{shift} with $\psi ^{\prime }$
chosen as in (a). $\blacksquare $

\subsubsection{Link with classical automorphic forms on $D^{\times }$}

To conclude this paragraph, we make explicit the link between adelic
automorphic forms for a definite quaternion algebra $D$ having holomorphic
weights, and classical automorphic forms for the algebraic $%
%TCIMACRO{\U{211a} }%
%BeginExpansion
\mathbb{Q}
%EndExpansion
$-group $\mathbb{D}$ associated to $D^{\times }$.

Set $A=E$ and let $\tau :\tprod\nolimits_{j=1}^{r}GL_{2}(\mathcal{O}_{F_{%
\mathfrak{P}_{j}}})\mathbb{\rightarrow }\limfunc{Aut}(W_{\tau })$ be a
weight for adelic automorphic forms on $D$ as considered in \ref{adelic HMF}%
; suppose $W_{\tau }=W_{\tau ^{\limfunc{alg}}}\otimes _{E}W_{\tau ^{\limfunc{%
sm}}}$, where $W_{\tau ^{\limfunc{sm}}}$ is a smooth irreducible $E$%
-representation of $\tprod\nolimits_{j=1}^{r}GL_{2}(\mathcal{O}_{F_{%
\mathfrak{P}_{j}}})$, and $W_{\tau ^{\limfunc{alg}}}=\dbigotimes%
\nolimits_{j=1}^{r}\dbigotimes\nolimits_{i=0}^{f_{j}-1}\left( \limfunc{Sym}%
\nolimits^{k_{i}^{(j)}-2}E^{2}\otimes \det\nolimits^{w_{i}^{(j)}}\right)
^{[i]}$ is an irreducible algebraic representation of $\mathbb{D}(%
%TCIMACRO{\U{211a} }%
%BeginExpansion
\mathbb{Q}
%EndExpansion
_{p})=(D\otimes _{%
%TCIMACRO{\U{211a} }%
%BeginExpansion
\mathbb{Q}
%EndExpansion
}%
%TCIMACRO{\U{211a} }%
%BeginExpansion
\mathbb{Q}
%EndExpansion
_{p})^{\times }=\tprod\nolimits_{j=1}^{r}GL_{2}(F_{\mathfrak{P}_{j}})$. We
assume that $k_{i}^{(j)}+2w_{i}^{(j)}-1$ equals some fixed integer $w$ for
all $1\leq j\leq r$ and all $0\leq i\leq f_{j}-1$. Recall that, as usual, we
see $F_{\mathfrak{P}_{j}}$ embedded in $E$ via $\sigma _{0}^{(j)}$ for $%
1\leq j\leq r$; we can also write $W_{\tau ^{\limfunc{alg}%
}}=\dbigotimes\nolimits_{\sigma :F\hookrightarrow E}\left( \limfunc{Sym}%
\nolimits^{k_{\sigma }-2}E^{2}\otimes \det\nolimits^{w_{\sigma }}\right) $.
Let $\psi :\left( \mathbb{A}_{F}^{\infty }\right) ^{\times }/F^{\times }%
\mathbb{\rightarrow }E^{\times }$ be a Hecke character compatible with $\tau 
$.

Fix an isomorphism $\bar{%
%TCIMACRO{\U{211a}}%
%BeginExpansion
\mathbb{Q}%
%EndExpansion
}_{p}\simeq 
%TCIMACRO{\U{2102} }%
%BeginExpansion
\mathbb{C}
%EndExpansion
$, inducing an embedding $E\hookrightarrow 
%TCIMACRO{\U{2102} }%
%BeginExpansion
\mathbb{C}
%EndExpansion
$.\ View $W_{\tau _{%
%TCIMACRO{\U{2102} }%
%BeginExpansion
\mathbb{C}
%EndExpansion
}^{\limfunc{alg}}}:=W_{\tau ^{\limfunc{alg}}}\otimes _{E}%
%TCIMACRO{\U{2102} }%
%BeginExpansion
\mathbb{C}
%EndExpansion
$ (resp. $W_{\tau _{%
%TCIMACRO{\U{2102} }%
%BeginExpansion
\mathbb{C}
%EndExpansion
}^{\limfunc{sm}}}:=W_{\tau ^{\limfunc{sm}}}\otimes _{E}%
%TCIMACRO{\U{2102} }%
%BeginExpansion
\mathbb{C}
%EndExpansion
$) as a complex representation of $\mathbb{D}(%
%TCIMACRO{\U{211d} }%
%BeginExpansion
\mathbb{R}
%EndExpansion
):=\left( D\otimes _{%
%TCIMACRO{\U{211a} }%
%BeginExpansion
\mathbb{Q}
%EndExpansion
}%
%TCIMACRO{\U{211d} }%
%BeginExpansion
\mathbb{R}
%EndExpansion
\right) ^{\times }\subset \mathbb{D}(%
%TCIMACRO{\U{2102} }%
%BeginExpansion
\mathbb{C}
%EndExpansion
)\simeq \mathbb{D}(\bar{%
%TCIMACRO{\U{211a}}%
%BeginExpansion
\mathbb{Q}%
%EndExpansion
}_{p})$ (resp. of $\tprod\nolimits_{j=1}^{r}GL_{2}(\mathcal{O}_{F_{\mathfrak{%
P}_{j}}})$). Let $W_{\tau _{%
%TCIMACRO{\U{2102} }%
%BeginExpansion
\mathbb{C}
%EndExpansion
}}:=W_{\tau }\otimes _{E}%
%TCIMACRO{\U{2102} }%
%BeginExpansion
\mathbb{C}
%EndExpansion
$ be the corresponding complex representation of $\tprod%
\nolimits_{j=1}^{r}GL_{2}(\mathcal{O}_{F_{\mathfrak{P}_{j}}})\times
\tprod\nolimits_{v|\infty }(\mathcal{O}_{D})_{v}^{\times }$.

Let $U^{\prime }$ be a compact open subgroup of $\left( D\otimes _{F}\mathbb{%
A}_{F}^{\infty }\right) ^{\times }$ such that $U^{\prime
}=\tprod\nolimits_{v\in \mathfrak{M}_{F,f}}U_{v}^{\prime }$, where $%
U_{v}^{\prime }=U_{v}$ if $v\NEG{|}p$ and, for $v_{j}|p$, $U_{v_{j}}^{\prime
}\subseteq GL_{2}(\mathcal{O}_{F_{\mathfrak{P}_{j}}})$ acts trivially on $%
W_{\tau ^{\limfunc{sm}}}$. Denote by $\mathcal{C}^{\infty }(D^{\times
}\backslash \left( D\otimes _{F}\mathbb{A}_{F}\right) ^{\times }/U^{\prime
}) $ the complex vector space of smooth functions $f:D^{\times }\backslash
\left( D\otimes _{F}\mathbb{A}_{F}\right) ^{\times }\mathbb{\rightarrow }%
%TCIMACRO{\U{2102} }%
%BeginExpansion
\mathbb{C}
%EndExpansion
$ which are invariant by the action of $U^{\prime }$. Let $W_{\tau _{%
%TCIMACRO{\U{2102} }%
%BeginExpansion
\mathbb{C}
%EndExpansion
}}^{\ast }$ be the $%
%TCIMACRO{\U{2102} }%
%BeginExpansion
\mathbb{C}
%EndExpansion
$-linear dual of $W_{\tau _{%
%TCIMACRO{\U{2102} }%
%BeginExpansion
\mathbb{C}
%EndExpansion
}}$.

Define a map:%
\begin{equation*}
\alpha :S_{\tau ,\psi }(U,E)\longrightarrow \limfunc{Hom}\nolimits_{\left(
D\otimes _{%
%TCIMACRO{\U{211a} }%
%BeginExpansion
\mathbb{Q}
%EndExpansion
}%
%TCIMACRO{\U{211d} }%
%BeginExpansion
\mathbb{R}
%EndExpansion
\right) ^{\times }}\left( W_{\tau _{%
%TCIMACRO{\U{2102} }%
%BeginExpansion
\mathbb{C}
%EndExpansion
}}^{\ast },\mathcal{C}^{\infty }(D^{\times }\backslash \left( D\otimes _{F}%
\mathbb{A}_{F}\right) ^{\times }/U^{\prime })\right)
\end{equation*}

\noindent by sending $f\in S_{\tau ,\psi }(U,E)$ to the assignment:%
\begin{equation*}
w^{\ast }\longmapsto (g\longmapsto w^{\ast }(\tau _{%
%TCIMACRO{\U{2102} }%
%BeginExpansion
\mathbb{C}
%EndExpansion
}^{\limfunc{alg}}(g_{\infty }^{-1})\tau ^{\limfunc{alg}}(g_{p})f(g^{\infty
})),
\end{equation*}

\noindent where $w^{\ast }\in W_{\tau _{%
%TCIMACRO{\U{2102} }%
%BeginExpansion
\mathbb{C}
%EndExpansion
}}^{\ast }$ and $g\in \left( D\otimes _{F}\mathbb{A}_{F}\right) ^{\times }$.
We have the following (cf. \cite{Ki}, 3.1.14):

\begin{proposition}
The map $\alpha $ identifies $S_{\tau ,\psi }(U,E)\otimes _{E}%
%TCIMACRO{\U{2102} }%
%BeginExpansion
\mathbb{C}
%EndExpansion
$ with a space of automorphic forms for the group $D^{\times }$ having
central character $\psi _{%
%TCIMACRO{\U{2102} }%
%BeginExpansion
\mathbb{C}
%EndExpansion
}$ given by $\psi _{%
%TCIMACRO{\U{2102} }%
%BeginExpansion
\mathbb{C}
%EndExpansion
}(g)=\limfunc{Nm}_{F/%
%TCIMACRO{\U{211a} }%
%BeginExpansion
\mathbb{Q}
%EndExpansion
}(g_{\infty })^{1-w}\limfunc{Nm}_{F/%
%TCIMACRO{\U{211a} }%
%BeginExpansion
\mathbb{Q}
%EndExpansion
}(g_{p})^{w-1}\psi (g^{\infty })$ for $g\in \left( D\otimes _{F}\mathbb{A}%
_{F}\right) ^{\times }$.

If $\pi =\tbigotimes\nolimits_{v}\pi _{v}$ is an irreducible automorphic
representation for the group $D^{\times }$, then $\pi $ is generated by an
element in $\alpha (f)(W_{\tau _{%
%TCIMACRO{\U{2102} }%
%BeginExpansion
\mathbb{C}
%EndExpansion
}}^{\ast })$ for some $f\in S_{\tau ,\psi }(U,E^{\prime })$, some $U$ small
enough and some $E^{\prime }\supseteq E$ big enough, if and only if $\pi
_{\infty }\simeq W_{\tau _{%
%TCIMACRO{\U{2102} }%
%BeginExpansion
\mathbb{C}
%EndExpansion
}^{\limfunc{alg}}}^{\ast }$ and $\tbigotimes\nolimits_{v|p}\pi _{v}$
contains $W_{\tau _{%
%TCIMACRO{\U{2102} }%
%BeginExpansion
\mathbb{C}
%EndExpansion
}^{\limfunc{sm}}}^{\ast }$ as a representation of $\tprod%
\nolimits_{j=1}^{r}GL_{2}(\mathcal{O}_{F_{\mathfrak{P}_{j}}})$.
\end{proposition}

Assume furthermore that $F/%
%TCIMACRO{\U{211a} }%
%BeginExpansion
\mathbb{Q}
%EndExpansion
$ has even degree and that we chose $\Sigma $ to be the empty set. Let $\tau 
$ be a holomorphic weight with parameters $(\vec{k},w)\in 
%TCIMACRO{\U{2124} }%
%BeginExpansion
\mathbb{Z}
%EndExpansion
_{\geq 2}^{g}\times 
%TCIMACRO{\U{2124} }%
%BeginExpansion
\mathbb{Z}
%EndExpansion
$ and let $\psi :\mathbb{A}_{F}^{\times }/F^{\times }\mathbb{\rightarrow }%
\bar{%
%TCIMACRO{\U{211a}}%
%BeginExpansion
\mathbb{Q}%
%EndExpansion
}_{p}^{\times }$ be a continuous character such that $\psi (a)=\left( 
\limfunc{Nm}a\right) ^{1-w}$ for all $a$ contained inside an open subgroup
of $(F\otimes _{%
%TCIMACRO{\U{211a} }%
%BeginExpansion
\mathbb{Q}
%EndExpansion
}%
%TCIMACRO{\U{211a} }%
%BeginExpansion
\mathbb{Q}
%EndExpansion
_{p})^{\times }$. Fix an isomorphism $\bar{%
%TCIMACRO{\U{211a}}%
%BeginExpansion
\mathbb{Q}%
%EndExpansion
}_{p}\simeq 
%TCIMACRO{\U{2102} }%
%BeginExpansion
\mathbb{C}
%EndExpansion
$ as before.

\noindent As a consequence of the classical Jacquet-Langlands theorem, we
can identify the complexification of the space $S_{\tau ,\psi }(U,\bar{%
%TCIMACRO{\U{211a}}%
%BeginExpansion
\mathbb{Q}%
%EndExpansion
}_{p})$ ($\vec{k}\neq \vec{2}$) with a space of regular algebraic cuspidal
automorphic representations $\pi $ of $GL_{2}(\mathbb{A}_{F})$ such that $%
\pi _{\infty }$ has weight $(\vec{k},w)$ and $\pi $ has central character $%
\psi _{\infty }$. If $\vec{k}=\vec{2}$ the identification works if we
consider, instead of $S_{\tau ,\psi }(U,\bar{%
%TCIMACRO{\U{211a}}%
%BeginExpansion
\mathbb{Q}%
%EndExpansion
}_{p})$, the quotient of $S_{\tau ,\psi }(U,\bar{%
%TCIMACRO{\U{211a}}%
%BeginExpansion
\mathbb{Q}%
%EndExpansion
}_{p})$ by the subspace of functions factoring through the reduced norm. For
a detailed formulation of these last facts, cf. Theorem 2.1 of \cite{H}\ and
Lemma 1.3 of \cite{Taylor}.

\subsection{Holomorphic weight shiftings via generalized Dickson invariants
and $D$-operators}

Let $q$ be a power of $p$. The intertwining operators between $\mathbb{F}%
_{q} $-representations of $GL_{2}(\mathbb{F}_{q})$ studied in Section \ref%
{sec3} allow us to produce weight shiftings between spaces of automorphic
forms having holomorphic weights.

\subsubsection{Main theorem}

Let us set $A=\mathcal{O}$ and let $D$, $U$, $\left( \tau ,W_{\tau }\right) $
and $\psi $ be as in \ref{adelic HMF}. Recall in particular that $U$ is
small enough, and that $\psi $ is compatible with $\tau $. For simplicity,
if $\tau $ is a holomorphic weight with parameters $(\vec{k},w)\in 
%TCIMACRO{\U{2124} }%
%BeginExpansion
\mathbb{Z}
%EndExpansion
_{\geq 2}^{g}\times 
%TCIMACRO{\U{2124} }%
%BeginExpansion
\mathbb{Z}
%EndExpansion
$ and $f\in S_{\tau ,\psi }(U,\mathcal{O})$, we also say that $f$ has weight 
$(\vec{k},w)$ or, sometimes, that $f$ has weight $\vec{k}$. Recall that we
write $\vec{k}=(\vec{k}^{(1)},...,\vec{k}^{(r)})$ with $\vec{k}%
^{(j)}=(k_{0}^{(j)},...,k_{f_{j}-1}^{(j)})\in 
%TCIMACRO{\U{2124} }%
%BeginExpansion
\mathbb{Z}
%EndExpansion
_{\geq 2}^{f_{j}}$ for $1\leq j\leq r$, and that we define the vector $\vec{w%
}^{(j)}=(w_{0}^{(j)},...,w_{f_{j}-1}^{(j)})\in 
%TCIMACRO{\U{2124} }%
%BeginExpansion
\mathbb{Z}
%EndExpansion
^{f_{j}}$ by the relations $k_{i}^{(j)}+2w_{i}^{(j)}-1=w$, for all $0\leq
i\leq f_{j}-1$.

\begin{theorem}
\label{general}Assume $\tau $ is a holomorphic $\mathcal{O}$-linear weight
with parameters $(\vec{k},w)\in 
%TCIMACRO{\U{2124} }%
%BeginExpansion
\mathbb{Z}
%EndExpansion
_{\geq 2}^{g}\times 
%TCIMACRO{\U{2124} }%
%BeginExpansion
\mathbb{Z}
%EndExpansion
$ with $w$ odd. Let $f=\min \{f_{1},...,f_{r}\}$ and fix an integer $\beta $
such that $1\leq \beta \leq f$. For any integers $i,j$ with $1\leq j\leq r$
and $0\leq i\leq f_{j-1}\ $choose:%
\begin{equation*}
a_{i}^{(j)}\in \{p^{\beta }-1,p^{\beta }+1\}.\noindent
\end{equation*}%
Set $\vec{a}=(\vec{a}^{(1)},...,\vec{a}^{(r)})$ with $\vec{a}%
^{(j)}=(a_{0}^{(j)},...,a_{f_{j}-1}^{(j)})$, and let $w^{\prime
}=w+(p^{\beta }-1)$. Assume at least one of the following conditions is
satisfied:

\begin{description}
\item[(*)] Let $j$ be any integer such that $1\leq j\leq r$ and $\beta
<f_{j} $. Then for any $i$ with $0\leq i\leq f_{j}-1$ and $%
a_{i}^{(j)}=p^{\beta }-1$, we have that $2<k_{i}^{(j)}\leq p+1$, $2\leq
k_{i+f_{j}-\beta }^{(j)}\leq p+1$ and if $i^{\prime }\neq i$ is another
integer such that $0\leq i^{\prime }\leq f_{j}-1$ and $a_{i^{\prime
}}^{(j)}=p^{\beta }-1$, we also have $i\not\equiv i^{\prime }-\beta (\func{%
mod}f_{j})$.

\item Let $j$ be any integer such that $1\leq j\leq r$ and $\beta =f_{j}$.
Then for any $i$ with $0\leq i\leq f_{j}-1$ and $a_{i}^{(j)}=p^{\beta }-1$,
we have that $2<k_{i}^{(j)}\leq p+1.$

\item[(**)] The weight $(\vec{k},w)$ is $p$-small and generic, i.e., $%
2<k_{i}^{(j)}\leq p+1$ for all $i,j$.
\end{description}

Let $\psi :\left( \mathbb{A}_{F}^{\infty }\right) ^{\times }/F^{\times }%
\mathbb{\rightarrow }\mathcal{O}^{\times }$ be a Hecke character compatible
with $\tau $. Then, if $\Omega $ is a Hecke eigensystem occurring in the
space $S_{\tau ,\psi }(U,\mathcal{O})$, there is a finite local extension of
discrete valuation rings $\mathcal{O}^{\prime }/\mathcal{O}$ and an $%
\mathcal{O}^{\prime }$-valued Hecke eigensystem $\Omega ^{\prime }$
occurring in holomorphic weight $(\vec{k}+\vec{a},w^{\prime })$ and with
associated Hecke character $\psi ^{\prime }$ such that: 
\begin{equation*}
\Omega ^{\prime }(\func{mod}\mathfrak{M}_{\mathcal{O}^{\prime }})=\Omega (%
\func{mod}\mathfrak{M}_{\mathcal{O}}).
\end{equation*}%
\noindent\ The character $\psi ^{\prime }$ is compatible with the weight $(%
\vec{k}+\vec{a},w^{\prime })$ and it can be chosen so that $\bar{\psi}%
^{\prime }=\bar{\psi}$.
\end{theorem}

\textbf{Proof }Recall that $\tau $ is the $\mathcal{O}$-linear
representation $\tau :\tprod\nolimits_{j=1}^{r}GL_{2}(\mathcal{O}_{F_{%
\mathfrak{P}_{j}}})\mathbb{\rightarrow }\limfunc{Aut}W$, where $%
W=\tbigotimes\nolimits_{j=1}^{r}W_{j}$, $W_{j}=\tbigotimes%
\nolimits_{i=0}^{f_{j}-1}\left( \limfunc{Sym}\nolimits^{k_{i}^{(j)}-2}%
\mathcal{O}^{2}\otimes \det\nolimits^{w_{i}^{(j)}}\right) ^{[i]}$, $%
k_{i}^{(j)}+2w_{i}^{(j)}-1=w$. The group $GL_{2}(\mathcal{O}_{F_{\mathfrak{P}%
_{j}}})$ acts on $W$ via the action on $W_{j}$ induced by the embedding $%
\sigma _{0}^{(j)}:GL_{2}(\mathcal{O}_{F_{\mathfrak{P}_{j}}})\hookrightarrow
GL_{2}(\mathcal{O})$. The superscript $[i]$ indicates twisting by the $i$th
power of the arithmetic Frobenius element of $\limfunc{Gal}(E/%
%TCIMACRO{\U{211a} }%
%BeginExpansion
\mathbb{Q}
%EndExpansion
_{p})$.

The $\mathbb{F}$-linear representation $\bar{W}_{j}:=W_{j}\otimes _{\mathcal{%
O}}\mathbb{F}$ of $GL_{2}(\mathcal{O}_{F_{\mathfrak{P}_{j}}})$\ factors
through the reduction map $GL_{2}(\mathcal{O}_{F_{\mathfrak{P}%
_{j}}})\rightarrow GL_{2}(\mathbb{F}_{\mathfrak{P}_{j}})$; using the
notation introduced in Section \ref{sec2}, we can identify $\bar{W}_{j}$
with the $\mathbb{F[}GL_{2}(\mathbb{F}_{\mathfrak{P}_{j}})]$-module 
\begin{equation*}
\bar{W}_{j}=\dbigotimes\nolimits_{i=0}^{f_{j}-1}\left(
M_{k_{i}^{(j)}-2}\otimes \det\nolimits^{w_{i}^{(j)}}\right) ^{[i]},
\end{equation*}

\noindent where we see $GL_{2}(\mathbb{F}_{\mathfrak{P}_{j}})\hookrightarrow
GL_{2}(\mathbb{F})$ via $\bar{\sigma}_{0}^{(j)}$, and the superscript $[i]$
indicates twisting by the $i$th power of the arithmetic Frobenius element of 
$\limfunc{Gal}(\mathbb{F}/\mathbb{F}_{p})$.

For any fixed integer $j$, $1\leq j\leq r$, let $\mathcal{T}_{j}=\mathcal{\{}%
i:a_{i}^{(j)}=p^{\beta }+1\mathcal{\}}$ and $\mathcal{D}_{j}=\mathcal{\{}%
i:a_{i}^{(j)}=p^{\beta }-1\mathcal{\}}$. For $i\in \mathcal{T}_{j}$ set $%
\vartheta _{i}^{(j)}:=\Theta _{f_{j}-\beta }^{[i]}$ if $\beta <f_{j}$ and $%
\vartheta _{i}^{(j)}:=\Theta ^{\lbrack i]}$ if $\beta =f_{j}$, where $\Theta
_{f_{j}-\beta }^{[i]}$ and $\Theta ^{\lbrack i]}$ are the generalized\
Dickson invariants for the group $GL_{2}(\mathbb{F}_{\mathfrak{P}%
_{j}})\simeq GL_{2}(\mathbb{F}_{p^{f_{j}}})$ as defined in \ref{GDI}. For $%
i\in \mathcal{D}_{j}$ set $\delta _{i}^{(j)}:=D_{f_{j}-\beta }^{[i]}$ if $%
\beta <f_{j}$ and $\delta _{i}^{(j)}:=D^{[i]}$ if $\beta =f_{j}$, where $%
D_{f_{j}-\beta }^{[i]}$ and $D^{[i]}$ are the generalized\ $D$-operators for 
$GL_{2}(\mathbb{F}_{\mathfrak{P}_{j}})$ defined in \ref{GDO}. Set:%
\begin{equation*}
\Lambda _{j}=\left( \dbigodot\nolimits_{i\in \mathcal{T}_{j}}\vartheta
_{i}^{(j)}\right) \circ \left( \dbigodot\nolimits_{i\in \mathcal{D}%
_{j}}\delta _{i}^{(j)}\right) ,
\end{equation*}

\noindent where the symbol $\tbigodot $ denotes composition of functions,
and each of the two composition factors above is computed by ordering $%
\mathcal{T}_{j}$ and $\mathcal{D}_{j}$ in the natural way. As seen in
section \ref{sec3}, the operators $\vartheta _{i}^{(j)}$ and $\delta
_{i}^{(j)}$ give rise to morphisms of $\mathbb{F}_{\mathfrak{P}_{j}}\mathbb{[%
}GL_{2}(\mathbb{F}_{\mathfrak{P}_{j}})]$-modules, and hence to morphisms of $%
\mathbb{F[}GL_{2}(\mathbb{F}_{\mathfrak{P}_{j}})]$-modules via the scalar
extension $\bar{\sigma}_{0}^{(j)}:\mathbb{F}_{\mathfrak{P}%
_{j}}\hookrightarrow \mathbb{F}.$ We deduce that $\Lambda _{j}$ induces a $%
GL_{2}(\mathbb{F}_{\mathfrak{P}_{j}})$-equivariant and $\mathbb{F}$-linear
morphism: 
\begin{equation*}
\Lambda _{j}:\bar{W}_{j}\mathbb{\rightarrow }\bar{W}_{j}^{\prime },
\end{equation*}

\noindent where $\bar{W}_{j}^{\prime }$ is the $\mathbb{F[}GL_{2}(\mathbb{F}%
_{\mathfrak{P}_{j}})]$-module: 
\begin{eqnarray*}
\bar{W}_{j}^{\prime } &:&=\dbigotimes\nolimits_{i\in \mathcal{T}_{j}}\left(
M_{k_{i}^{(j)}+(p^{\beta }+1)-2}\otimes \det\nolimits^{w_{i}^{(j)}-1}\right)
^{[i]} \\
&&\otimes \dbigotimes\nolimits_{i\in \mathcal{D}_{j}}\left(
M_{k_{i}^{(j)}+(p^{\beta }-1)-2}\otimes \det\nolimits^{w_{i}^{(j)}}\right)
^{[i]}.
\end{eqnarray*}

\noindent Indeed, by Theorem \ref{TTh}, $\Theta _{f_{j}-\beta }^{[i]}$
increases $k_{i}^{(j)}$ by $1$, $k_{i+f_{j}-\beta }^{(j)}$ by $p^{\beta }$, $%
w_{i}^{(j)}$ by $-1$, and does not change $k_{s}^{(j)}$ for $s\neq
i,i+f_{j}-\beta $ or $w_{s}^{(j)}$ for $s\neq i$; $\Theta ^{\lbrack i]}$
increases $k_{i}^{(j)}$ by $p^{f_{j}}+1$, $w_{i}^{(j)}$ by $-1$, and does
not change $k_{s}^{(j)}$ or $w_{s}^{(j)}$ for $s\neq i$. On the other hand,
by Theorem \ref{DTh}, the operator $D_{f_{j}-\beta }^{[i]}$ increases $%
k_{i}^{(j)}$ by $-1$, $k_{i+f_{j}-\beta }^{(j)}$ by $p^{\beta }$, and does
not change $k_{s}^{(j)}$ for $s\neq i,i+f_{j}-\beta $ or $w_{s}^{(j)}$ for
any $s$; $D^{[i]}$ increases $k_{i}^{(j)}$ by $p^{f_{j}}-1$, and does not
change $k_{s}^{(j)}$ for $s\neq i$ or $w_{s}^{(j)}$ for any $s$.

By Theorem \ref{TTh}, $\tbigodot\nolimits_{i\in \mathcal{T}_{j}}\vartheta
_{i}^{(j)}$ is injective. If $(\ast )$ is satisfied, the injectivity
statement of Theorem \ref{DTh} implies that $\tbigodot\nolimits_{i\in 
\mathcal{D}_{j}}\delta _{i}^{(j)}$ is injective on $\bar{W}_{j}$. The image
of $\tbigotimes\nolimits_{i=0}^{f_{j}-1}\left( X^{k_{i}^{(j)}-2}\otimes
1\right) ^{[i]}\in \bar{W}_{j}$\ under $\tbigodot\nolimits_{i\in \mathcal{D}%
_{j}}\delta _{i}^{(j)}$ is easily seen to be of the form $%
\tprod\nolimits_{i\in \mathcal{D}_{j}}(k_{i}^{(j)}-2)\cdot u$ for some
non-zero $u\in \bar{W}_{j}$. If $(\ast \ast )$ holds, $\tprod\nolimits_{i\in 
\mathcal{D}_{j}}(k_{i}^{(j)}-2)$ is non-zero in $\mathbb{F}$ and, being $%
\bar{W}_{j}$ an irreducible representation of $GL_{2}(\mathbb{F}_{\mathfrak{P%
}_{j}})$,$\ $we deduce that $\tbigodot\nolimits_{i\in \mathcal{D}_{j}}\delta
_{i}^{(j)}$ is injective on $\bar{W}_{j}$. \noindent We conclude that under
assumptions $(\ast )$ or $(\ast \ast )$, all the maps $\Lambda _{j}$ for $%
1\leq j\leq r$ are injective.

Let $b_{i}^{(j)}=-1$ if $i\in \mathcal{T}_{j}$ and $b_{i}^{(j)}=0$ if $i\in 
\mathcal{D}_{j}$. Define the $\mathcal{O}[GL_{2}(\mathcal{O}_{F_{\mathfrak{P}%
_{j}}})]$-module:%
\begin{equation*}
W_{j}^{\prime }=\dbigotimes\nolimits_{i=0}^{f_{j}-1}\left( \limfunc{Sym}%
\nolimits^{k_{i}^{(j)}+a_{i}^{(j)}-2}\mathcal{O}^{2}\otimes
\det\nolimits^{w_{i}^{(j)}+b_{i}^{(j)}}\right) ^{[i]},
\end{equation*}

\noindent so that $W_{j}^{\prime }\otimes _{\mathcal{O}}\mathbb{F}=\bar{W}%
_{j}^{\prime }$ as $\mathbb{F}$-representations of $GL_{2}(\mathcal{O}_{F_{%
\mathfrak{P}_{j}}})$ or, equivalently, of $GL_{2}(\mathbb{F}_{\mathfrak{P}%
_{j}})$. Set $W^{\prime }=\tbigotimes\nolimits_{j=1}^{r}W_{j}^{\prime }$ and
denote by $\tau ^{\prime }$ the action of $U$ on $W^{\prime }$ induced by
the projection $U\rightarrow \tprod\nolimits_{j=1}^{r}GL_{2}(\mathcal{O}_{F_{%
\mathfrak{P}_{j}}})$. Let $w^{\prime }=w+(p^{\beta }-1)$; for all the values
of $i$ and $j$ for which the following integers are defined, we have $%
k_{i}^{(j)}+a_{i}^{(j)}\geq k_{i}^{(j)}\geq 2$ and:%
\begin{eqnarray*}
&&\left( k_{i}^{(j)}+a_{i}^{(j)}\right) +2\left(
w_{i}^{(j)}+b_{i}^{(j)}\right) -1 \\
&=&\left( k_{i}^{(j)}+2w_{i}^{(j)}-1\right) +p^{\beta }-1 \\
&=&w+(p^{\beta }-1).
\end{eqnarray*}

\noindent Therefore $\tau ^{\prime }$ is a holomorphic weight for
automorphic forms on $D$ with parameters $(\vec{k}+\vec{a},w^{\prime })\in 
%TCIMACRO{\U{2124} }%
%BeginExpansion
\mathbb{Z}
%EndExpansion
_{\geq 2}^{g}\times 
%TCIMACRO{\U{2124} }%
%BeginExpansion
\mathbb{Z}
%EndExpansion
.$

The injections $\Lambda _{j}$ ($1\leq j\leq r$) constructed above allow us
to see $\bar{W}=\tbigotimes\nolimits_{j=1}^{r}\bar{W}_{j}$ as an $\mathbb{F}$%
-linear $U$-subrepresentation of $\bar{W}^{\prime
}=\tbigotimes\nolimits_{j=1}^{r}\bar{W}_{j}^{\prime }$. Since $w$ is odd and 
$w\equiv w^{\prime }(\func{mod}p-1)$, we can apply Proposition \ref{main}.
We conclude that there exists a Hecke character $\psi ^{\prime }:\left( 
\mathbb{A}_{F}^{\infty }\right) ^{\times }/F^{\times }\mathbb{\rightarrow }%
\mathcal{O}^{\times }$ compatible with $\tau ^{\prime }$ and such that $\bar{%
\psi}^{\prime }=\bar{\psi}$; furthermore, for any Hecke eigensystem $\Omega $%
\ occurring in $S_{\tau ,\psi }(U,\mathcal{O})$ there is a finite extension
of discrete valuation rings $\mathcal{O}^{\prime }/\mathcal{O}$ with $%
\mathfrak{M}_{\mathcal{O}^{\prime }}\cap \mathcal{O}=\mathfrak{M}_{\mathcal{O%
}}$ and a Hecke eigensystem $\Omega ^{\prime }$ occurring in $S_{\tau
^{\prime },\psi ^{\prime }}(U,\mathcal{O}^{\prime })$ such that $\Omega
^{\prime }(\func{mod}\mathfrak{M}_{\mathcal{O}^{\prime }})=\Omega (\func{mod}%
\mathfrak{M}_{\mathcal{O}}).$ $\blacksquare $

\bigskip

\begin{corollary}
\noindent Under the same notation and assumptions of Theorem \ref{general},
any $\mathbb{\bar{F}}_{p}$-linear continuous Galois representation arising
from a Hecke eigenform in $S_{\tau ,\psi }(U,\mathcal{O})$, where $\tau $ is
a holomorphic weight of parameter $\vec{k}$,\ also arises from an eigenform
in $S_{\tau ^{\prime },\psi ^{\prime }}(U,\mathcal{\bar{%
%TCIMACRO{\U{2124}}%
%BeginExpansion
\mathbb{Z}%
%EndExpansion
}}_{p})$, where $\tau ^{\prime }$ is a holomorphic weight of parameters $%
\vec{k}+\vec{a}$ and $\psi ^{\prime }$ is some $\mathcal{O}^{\times }$%
-valued Hecke character compatible with $\tau ^{\prime }$ and such that $%
\bar{\psi}^{\prime }=\bar{\psi}$.
\end{corollary}

\begin{remark}
\label{more}We remark what follows:

\begin{enumerate}
\item Condition $(\ast )$ of Theorem \ref{general} is true if, for example,
for any $j$ with $1\leq j\leq r$, there is at most one $i$, $0\leq i\leq
f_{j}-1$, such that $a_{i}^{(j)}=p^{\beta }-1$, and for these values of $i$
and $j$ we have $2<k_{i}^{(j)}\leq p+1$ and $2\leq k_{i+f_{j}-\beta
}^{(j)}\leq p+1$.

\item The reason for which in the above result we limit $a_{i}^{(j)}$ to be
in the set $\{p^{\beta }-1,p^{\beta }+1\}$ for all $i,j$ is that we want to
preserve the holomorphicity of the weights of the automorphic forms
involved. More weight shiftings are possible using the generalized Dickson
and $D$-operators if we do not impose the holomorphicity condition. On the
other side, we will see in \ref{more3} that when $g>1$ our operators allow
more holomorphic weight shiftings than the ones described in the Theorem \ref%
{general}.

\item As a consequence of Remark \ref{AG1} and Remark \ref{AG2}, the above
result gives rise in general to more holomorphic weight shiftings than the
ones obtained by the theory of generalized theta operators and Hasse
invariants for geometric $(\func{mod}p)$ Hilbert modular forms (cf. \ref%
{motivations}).
\end{enumerate}
\end{remark}

\subsubsection{Analysis of the case $f_{j}<3\label{more3}$}

We determine additional holomorphic weight shiftings using the generalized
Dickson and $D$-operators. Since the combinatorics involved in the
computations becomes very complicated as $\max \{f_{1},...,f_{r}\}$ grows,
we assume that $f_{j}<3$ for all $j.$ A procedure similar to the one
described below could be applied in greater generality.

The most interesting cases for us arise when some of the residue degrees $%
f_{j}$ equal two, so that we assume without loss of generality $g=2$, $p%
\mathcal{O}_{F}=\mathfrak{P}$ and $[\mathbb{F}_{\mathfrak{P}}:\mathbb{F}%
_{p}]=2$. We maintain the notation introduced at the beginning of the
section, but since $r=1$ we drop the index $j$ wherever it appeared before.
Notice that we have $E=F_{\mathfrak{P}}$, $\mathbb{F}=\mathbb{F}_{\mathfrak{P%
}}$ and we can assume that $\sigma _{0}$ (resp. $\bar{\sigma}_{0}$) is the
identity automorphism of $F_{\mathfrak{P}}$ (resp. $\mathbb{F}_{\mathfrak{P}%
} $).

Fix $\vec{k}=(k_{0},k_{1})\in 
%TCIMACRO{\U{2124} }%
%BeginExpansion
\mathbb{Z}
%EndExpansion
_{\geq 2}^{2}$\ and $\vec{w}=(w_{0},w_{1})\in 
%TCIMACRO{\U{2124} }%
%BeginExpansion
\mathbb{Z}
%EndExpansion
^{2}$\ such that $k_{i}+2w_{i}-1=w$ ($i=0,1$) with $w$ odd. Let $(\tau ,W)$
be the holomorphic $\mathcal{O}$-linear weight with parameters $(\vec{k},w)$%
, so that the reduction modulo $\mathfrak{P}$ of $W$ is the $\mathbb{F}_{%
\mathfrak{P}}\mathbb{[}GL_{2}(\mathbb{F}_{\mathfrak{P}})]$-module:%
\begin{equation*}
\bar{W}=\left( M_{k_{0}-2}\otimes \det\nolimits^{w_{0}}\right) \otimes
\left( M_{k_{1}-2}\otimes \det\nolimits^{w_{1}}\right) ^{[1]}.
\end{equation*}

\noindent Fix non-negative integers $n,m,r,s,t,u,v,z$ and let:%
\begin{equation*}
\Lambda =\Theta ^{\lbrack 1],u}\circ \Theta ^{\lbrack 0],t}\circ \Theta
_{1}^{[1],m}\circ \Theta _{1}^{[0],n}\circ D^{[1],z}\circ D^{[0],v}\circ
D_{1}^{[1],s}\circ D_{1}^{[0],r},
\end{equation*}

\noindent where the above operators are defined as in \ref{GDI} and \ref{GDO}%
. $\Lambda $ defines a $\mathbb{F}_{\mathfrak{P}}\mathbb{[}GL_{2}(\mathbb{F}%
_{\mathfrak{P}})]$-homomorphism having source $\bar{W}$ as long as $r+2\leq
k_{0}$ and $s+2\leq k_{1}+pr$; we assume therefore:%
\begin{equation}
\left\{ 
\begin{array}{c}
r+2\leq k_{0} \\ 
s+2\leq k_{1}.%
\end{array}%
\right.  \tag{$\clubsuit ^{\prime }$}
\end{equation}%
If $(\clubsuit ^{\prime })$ holds, we have $\Lambda :\bar{W}\mathbb{%
\rightarrow }\bar{W}^{\prime }$, where:%
\begin{equation*}
\bar{W}^{\prime }=\left( M_{k_{0}^{\prime }-2}\otimes
\det\nolimits^{w_{0}^{\prime }}\right) \otimes \left( M_{k_{1}^{\prime
}-2}\otimes \det\nolimits^{w_{1}^{\prime }}\right) ^{[1]},
\end{equation*}

\noindent with:

\begin{eqnarray*}
k_{0}^{\prime } &=&k_{0}+n+pm-r+ps+(p^{2}+1)t+(p^{2}-1)v \\
k_{1}^{\prime } &=&k_{1}+pn+m+pr-s+(p^{2}+1)u+(p^{2}-1)z \\
w_{0}^{\prime } &=&w_{0}-n-t+\alpha (p^{2}-1) \\
w_{1}^{\prime } &=&w_{1}-m-u+\beta (p^{2}-1).
\end{eqnarray*}

\noindent Here $\alpha ,\beta $ can be chosen to be any integers, as $%
\det^{p^{2}-1}=1$ on $GL_{2}(\mathbb{F}_{\mathfrak{P}})$.

Assume that the following are satisfied:

\begin{description}
\item[(A)] relations $(\clubsuit ^{\prime })$ hold and $\Lambda $ is
injective;

\item[(B)] $k_{0}^{\prime },k_{1}^{\prime }\geq 2$;

\item[(C)] $k_{0}^{\prime }+2w_{0}^{\prime }-1=k_{1}^{\prime
}+2w_{1}^{\prime }-1=:w^{\prime };$

\item[(D)] $w\equiv w^{\prime }(\func{mod}p-1)$.
\end{description}

\noindent Then we can apply Proposition \ref{main} to obtain holomorphic
weight shiftings for Hecke eigensystems associated to automorphic forms on $%
D $. We therefore want to translate the above four conditions into relations
between the integral parameters $k_{i},w_{i},\alpha ,\beta ,n,m,r,s,t,u,v,z$.

We easily see that:%
\begin{eqnarray*}
k_{0}^{\prime }+2w_{0}^{\prime }-1 &=&w-n+pm-r+ps+(p^{2}-1)(t+v+2\alpha ) \\
k_{1}^{\prime }+2w_{1}^{\prime }-1 &=&w+pn-m+pr-s+(p^{2}-1)(u+z+2\beta ),
\end{eqnarray*}

\noindent so that condition $(\mathbf{C})$ is equivalent to:%
\begin{equation}
(m-n)+(s-r)+(p-1)\left( (t-u)+(v-z)+2(\alpha -\beta )\right) =0. 
\tag{$\QTR{bf}{\spadesuit }$}
\end{equation}

\noindent Computing $m$ and $r$ from $(\mathbf{\spadesuit })$ we obtain:%
\begin{eqnarray*}
k_{0}^{\prime } &=&k_{0}+(p+1)(n+t)+(p-1)(r+v)+p(p-1)(u+z+2(\beta -\alpha ))
\\
k_{1}^{\prime } &=&k_{1}+(p+1)(m+u)+(p-1)(s+z)+p(p-1)(t+v+2(\alpha -\beta ))
\\
w^{\prime } &=&w+(p-1)(n+t+r+v+2\alpha +p(u+z+2\beta )).
\end{eqnarray*}

\noindent Condition $(\mathbf{D})$ is then automatically satisfied, and $%
\mathbf{(B)}$ holds if $\alpha =\beta $. If $r=s=v=z=0$, condition $\mathbf{%
(A)}$ is satisfied, as the generalized Dickson invariants induce injective
morphisms of $\mathbb{F}_{\mathfrak{P}}\mathbb{[}GL_{2}(\mathbb{F}_{%
\mathfrak{P}})]$-modules, and $(\clubsuit ^{\prime })$ is then a consequence
of $k_{0},k_{1}\geq 2$.

We claim that if $2<k_{0,}k_{1}\leq p+1$ and the non-negative integers $%
r,s,v,z$ satisfy:%
\begin{equation}
\left\{ 
\begin{array}{c}
r+v+2\leq k_{0} \\ 
s+z+2\leq k_{1},%
\end{array}%
\right.  \tag{$\clubsuit $}
\end{equation}

\noindent then $\mathbf{(A)}$ holds. (Notice that $(\clubsuit )$ implies $%
(\clubsuit ^{\prime })$). To prove this, first observe that if $%
2<k_{0,}k_{1}\leq p+1$, then $\bar{W}$ is irreducible for the action of $%
GL_{2}(\mathbb{F}_{\mathfrak{P}})$, so we only need to show that under the
above assumptions $\Lambda \neq 0$. Write $a=k_{0}-2$ and $b=k_{1}-2$ and
denote for simplicity the element $\left( X^{a}\otimes 1\right) \otimes
\left( X^{b}\otimes 1\right) ^{[1]}\ $of $\bar{W}$ by $X^{a}\otimes X^{b}$.
We have:

\begin{eqnarray}
&&\left( D^{[1],z}\circ D^{[0],v}\circ D_{1}^{[1],s}\circ
D_{1}^{[0],r}\right) \left( X^{a}\otimes X^{b}\right) 
\TCItag{$\blacklozenge $} \\
&=&c\cdot X^{a-r+ps+(p^{2}-1)v}\otimes X^{b+pr-s+(p^{2}-1)z},  \notag
\end{eqnarray}

\noindent where:%
\begin{equation*}
c=\frac{a!}{(a-r)!}\cdot \frac{(b+pr)!}{(b+pr-s)!}\cdot \frac{(a-r+ps)!}{%
(a-r+ps-v)!}\cdot \frac{(b+pr-s)!}{(b+pr-s-z)!}(\func{mod}p).
\end{equation*}

\noindent The exponents in the right hand side of $(\blacklozenge )$ and the
integers in the above formula for $c$ are non-negative under the assumption $%
(\clubsuit )$. Since $0<a,b\leq p-1$, $(\clubsuit )$ also implies that $p$
does not divide the integer $\frac{a!}{(a-r)!}\cdot \frac{(b+pr)!}{(b+pr-s)!}
$. Assume $r+v\leq a$ and $v>0$; if $p$ divided the integer: 
\begin{equation*}
\frac{(a-r+ps)!}{(a-r+ps-v)!}=\dprod\nolimits_{j=0}^{v-1}(a-r+ps-j),
\end{equation*}

\noindent then $p$ would divide $a-r-j$ for some $0\leq j\leq v-1$, which is
impossible as $1\leq a-r-j\leq p-1$. Similarly, if $s+z\leq b$ we see that $%
p $ does not divide $\frac{(b-s+pr)!}{(b-s+pr-z)!}$.

\noindent We conclude that $c\neq 0$ and hence $D^{[1],z}\circ
D^{[0],v}\circ D_{1}^{[1],s}\circ D_{1}^{[0],r}\neq 0$. The injectivity of
the generalized Dickson invariants implies then the claim.

\bigskip

Let us set $A=\mathcal{O=O}_{F_{\mathfrak{P}}}$ and let $D$, $U$ be as in %
\ref{adelic HMF}. The above considerations and Proposition \ref{main} prove
the following:

\begin{theorem}
\label{main2}Assume $g=f=2$ and let $\tau $ be an $\mathcal{O}$-linear
holomorphic weight for automorphic forms on $D$ of parameters $%
(k_{0},k_{1};w)\in 
%TCIMACRO{\U{2124} }%
%BeginExpansion
\mathbb{Z}
%EndExpansion
_{\geq 2}^{2}\times 
%TCIMACRO{\U{2124} }%
%BeginExpansion
\mathbb{Z}
%EndExpansion
$ with $w$ odd; let $\psi $ be a Hecke character compatible with $\tau $.
Fix $\alpha \in 
%TCIMACRO{\U{2124} }%
%BeginExpansion
\mathbb{Z}
%EndExpansion
$ and non-negative integers $n,m,r,s,t,u,v$ and $z$. Assume at least one of
the following two conditions is satisfied:

\begin{description}
\item[(*)] $r=s=v=z=0;$

\item[(**)] $2<k_{0,}k_{1}\leq p+1$, and $r+v\leq k_{0}-2$, $s+z\leq
k_{1}-2. $
\end{description}

\noindent Assume furthermore that the relation:%
\begin{equation}
(m-n)+(s-r)=(p-1)\cdot \left( (u-t)+(z-v)\right)  \tag{$\spadesuit $}
\end{equation}

\noindent \noindent holds. Define:%
\begin{equation*}
\left\{ 
\begin{array}{l}
k_{0}^{\prime }=k_{0}+(p+1)(n+t)+(p-1)(r+v)+p(p-1)(u+z) \\ 
k_{1}^{\prime }=k_{1}+(p+1)(m+u)+(p-1)(s+z)+p(p-1)(t+v) \\ 
w^{\prime }=w+(p-1)(n+t+r+v+2\alpha +p(u+z+2\alpha )).%
\end{array}%
\right.
\end{equation*}

Then if $\Omega $ is a Hecke eigensystem occurring in $S_{\tau ,\psi }(U,%
\mathcal{O})$, there is a finite local extension of discrete valuation rings 
$\mathcal{O}^{\prime }/\mathcal{O}$ and an $\mathcal{O}^{\prime }$-valued
Hecke eigensystem $\Omega ^{\prime }$ occurring in holomorphic weight $%
(k_{0}^{\prime },k_{1}^{\prime };w^{\prime })$ and with associated Hecke
character $\psi ^{\prime }$\ such that: 
\begin{equation*}
\Omega ^{\prime }(\func{mod}\mathfrak{M}_{\mathcal{O}^{\prime }})=\Omega (%
\func{mod}\mathfrak{M}_{\mathcal{O}}).
\end{equation*}
The character $\psi ^{\prime }$ is compatible with the weight $%
(k_{0}^{\prime },k_{1}^{\prime };w^{\prime })$ and it can be chosen so that $%
\bar{\psi}^{\prime }=\bar{\psi}$.
\end{theorem}

\begin{remark}
Many of the weight shiftings produced by Theorem \ref{main2} do not arise
from Theorem \ref{general} or from the operators described in \ref%
{motivations}.
\end{remark}

\section{Shiftings for weights containing $(2,...,2)$-blocks\label{22222}}

While the generalized Dickson invariants induce injective maps on the
trivial $\mathbb{F}$-representation of $GL_{2}(\mathbb{F}_{\mathfrak{P}%
_{j}}) $, the $D$-operators are identically zero on this module. Starting
with automorphic forms whose weight contains a $(2,...,2)$-block (cf.
definition in \ref{conv_on_emb}), we can then produce weight shiftings
through the operators $\Theta _{\alpha }^{[\beta ]}$ but we cannot always
successfully use the operators $D_{\alpha }^{[\beta ]}$. On the other side,
the study of weight shiftings "by $\overrightarrow{p-1}$" for automorphic
forms whose weight contains a $(2,...,2)$-block is motivated by the weight
part of Serre's modularity conjecture for totally real fields (cf. Remark %
\ref{SerreF}\ below).

In this section we present a result of Edixhoven and Khare (cf. \cite{EK})
to produce weight shiftings "by $\overrightarrow{p-1}$" starting from forms
whose weight is not necessarily parallel but contains $(2,...,2)$-blocks
relative to some primes of $F$ above $p$. We always assume that $p>2$ is
unramified in the totally real number field $F$.

\bigskip

We keep the notation introduced in \ref{sec4}, and we furthermore assume
that $F$ has even degree over $%
%TCIMACRO{\U{211a} }%
%BeginExpansion
\mathbb{Q}
%EndExpansion
$ and that the quaternion $F$-algebra $D$ is ramified at all and only the
infinite places of $F$, i.e., $\Sigma =\varnothing $. We fix an isomorphism $%
\left( D\otimes _{F}\mathbb{A}_{F}^{\infty }\right) ^{\times }\simeq GL_{2}(%
\mathbb{A}_{F}^{\infty })$.

The symbols $\mathbb{F}$, $U$, $\left( \tau ,W_{\tau }\right) $, $\psi $, $S$
and $\mathbb{T}_{S,\mathbb{F}}^{univ}$\ will have the same meaning as in \ref%
{adelic HMF}. We assume that $\tau $ is a (non necessarily holomorphic) $%
\mathbb{F}$-linear weight with parameters $(\vec{k},\vec{w})\in 
%TCIMACRO{\U{2124} }%
%BeginExpansion
\mathbb{Z}
%EndExpansion
_{\geq 2}^{g}\times 
%TCIMACRO{\U{2124} }%
%BeginExpansion
\mathbb{Z}
%EndExpansion
^{g}$, where $\vec{k}=(\vec{k}^{(1)},...,\vec{k}^{(r)})$ and $\vec{k}%
^{(j)}=(k_{0}^{(j)},...,k_{f_{j}-1}^{(j)})\in 
%TCIMACRO{\U{2124} }%
%BeginExpansion
\mathbb{Z}
%EndExpansion
_{\geq 2}^{f_{j}}$; $\vec{w}=(\vec{w}^{(1)},...,\vec{w}^{(r)})$ and $\vec{w}%
^{(j)}=(w_{0}^{(j)},...,w_{f_{j}-1}^{(j)})\in 
%TCIMACRO{\U{2124} }%
%BeginExpansion
\mathbb{Z}
%EndExpansion
^{f_{j}}$, for $1\leq j\leq r$.

We write $W_{\tau }=\tbigotimes\nolimits_{j=1}^{r}W_{\tau _{j}}$ where $%
W_{\tau _{j}}$ is the $\mathbb{F}$-representation of $GL_{2}(\mathcal{O}_{F_{%
\mathfrak{P}_{j}}})$ defined by:%
\begin{equation*}
W_{\tau _{j}}=\dbigotimes\nolimits_{i=0}^{f_{j}-1}\left( \limfunc{Sym}%
\nolimits^{k_{i}^{(j)}-2}\mathbb{F}^{2}\otimes
\det\nolimits^{w_{i}^{(j)}}\right) ^{[i]}.
\end{equation*}

\noindent If the weight $\tau $ is holomorphic, it is also determined by the
pair $(\vec{k},w)\in 
%TCIMACRO{\U{2124} }%
%BeginExpansion
\mathbb{Z}
%EndExpansion
_{\geq 2}^{g}\times 
%TCIMACRO{\U{2124} }%
%BeginExpansion
\mathbb{Z}
%EndExpansion
$ where $k_{i}^{(j)}+2w_{i}^{(j)}-1=w$, for all $i$ and $j$.

Choose a prime $\mathfrak{P}$ of $F$ above $p$ and let $\varpi $ be a fixed
choice of uniformizer for the ring of integers of the completion of $F$ at $%
\mathfrak{P}$. We can assume, up to relabeling, that $\mathfrak{P=P}_{1}$.
Define the matrix of $GL_{2}(F_{\mathfrak{P}_{1}})$:%
\begin{equation*}
\Pi =\left( 
\begin{array}{cc}
1 & 0 \\ 
0 & \varpi%
\end{array}%
\right) ,
\end{equation*}

\noindent and view it as an element of $GL_{2}(\mathbb{A}_{F}^{\infty })$
whose components away from $\mathfrak{P}_{1}$ are trivial.

If $g\ $is an element of $GL_{2}(\mathbb{A}_{F}^{\infty })$ and $Q$ is a
finite set of finite places of $F$, we denote by $g^{Q}$ the element of $%
GL_{2}(\mathbb{A}_{F}^{\infty })$ whose components at each place of $Q$ are
trivial, and whose components away from $Q$ coincide with those of $g$. We
let $g_{Q}=g/g^{Q}$. A similar convention is used for subgroups of $GL_{2}(%
\mathbb{A}_{F}^{\infty })$ which are products of subgroups of $GL_{2}(F_{v})$
for $v$ varying over the finite places of $F$. In particular, by assumption
we have $U_{p}=GL_{2}(\mathcal{O}_{F}\otimes _{%
%TCIMACRO{\U{2124} }%
%BeginExpansion
\mathbb{Z}
%EndExpansion
}%
%TCIMACRO{\U{2124} }%
%BeginExpansion
\mathbb{Z}
%EndExpansion
_{p})$.

\noindent We denote the action by right translation of $GL_{2}(\mathbb{A}%
_{F}^{\infty })$ on $S_{\tau ,\psi }(U,\mathbb{F})$ by a dot.

Set:%
\begin{equation*}
U_{0}=\left\{ u\in U:u_{\mathfrak{P}_{1}}\equiv \left( 
\begin{array}{cc}
\ast & \ast \\ 
0 & \ast%
\end{array}%
\right) (\func{mod}\varpi )\right\} .
\end{equation*}

\noindent By restricting $\tau $ to $U_{0}$, we define $S_{\tau ,\psi
}(U_{0},\mathbb{F})$ as in Definition \ref{defhmf}; notice that the level of
the automorphic forms belonging to this space is \textit{not} prime-to-$p$.

We have the following result, which is a not-prime-to-$p$ version of Lemma
3.1 of \cite{Taylor}:

\begin{lemma}
\label{dege}Assume that $\tau $ is an irreducible (non necessarily
holomorphic) $\mathbb{F}$-linear weight with parameters $(\vec{k},\vec{w}%
)\in 
%TCIMACRO{\U{2124} }%
%BeginExpansion
\mathbb{Z}
%EndExpansion
_{\geq 2}^{g}\times 
%TCIMACRO{\U{2124} }%
%BeginExpansion
\mathbb{Z}
%EndExpansion
^{g}$\ such that $\vec{k}^{(1)}=\vec{2}$. Then the map:%
\begin{equation*}
\alpha :S_{\tau ,\psi }(U,\mathbb{F})\oplus S_{\tau ,\psi }(U,\mathbb{F}%
)\longrightarrow S_{\tau ,\psi }(U_{0},\mathbb{F})
\end{equation*}

\noindent defined by:%
\begin{equation*}
\left( f_{1},f_{2}\right) \longmapsto f_{1}+\Pi \cdot f_{2}
\end{equation*}

\noindent is a Hecke-equivariant $\mathbb{F}$-morphism whose kernel is
Eisenstein, i.e., the localization $\left( \ker \alpha \right) _{\mathfrak{M}%
}$ vanishes for all maximal ideals $\mathfrak{M}\ $of $\mathbb{T}_{S,\mathbb{%
F}}^{univ}$ which are non-Eisenstein.
\end{lemma}

\textbf{Proof }It is straightforward to check that $\alpha $ is well
defined, using the fact that $GL_{2}(\mathcal{O}_{F_{\mathfrak{P}_{1}}})$
acts on $W_{\tau _{1}}$ via an integral power of the $(\func{mod}\varpi )$
determinant character. Also, $\alpha $ is equivariant for the action of the
algebra $\mathbb{T}_{S,\mathbb{F}}^{univ}$.

Write $\Pi U\Pi ^{-1}=U^{\mathfrak{P}_{1}}\times \Pi GL_{2}(\mathcal{O}_{F_{%
\mathfrak{P}_{1}}})\Pi ^{-1}$. Define an $\mathbb{F}$-linear action of the
subgroup $\Pi U\Pi ^{-1}$ of $GL_{2}(\mathbb{A}_{F}^{\infty })$ on $W_{\tau
} $ by letting $U^{\mathfrak{P}_{1}}$ act on $\tbigotimes%
\nolimits_{j=2}^{r}W_{\tau _{j}}$ via the restriction of $\tau $ to $U^{%
\mathfrak{P}_{1}}$, and by letting $\Pi GL_{2}(\mathcal{O}_{F_{\mathfrak{P}%
_{1}}})\Pi ^{-1}$ act on $W_{\tau _{1}}$ via the reduction modulo $\varpi $
of the determinant character raised to the power of $\tsum%
\nolimits_{i=0}^{f_{1}-1}w_{i}^{(1)}p^{i}$. Observe that this action is
compatible with the given action $\tau $ of $U$ on $W_{\tau }$.

If $\left( f_{1},f_{2}\right) \in \ker \alpha $, we see that $f_{1}\left(
gu\right) =u^{-1}f_{1}(g)$ for all $u\ $in $U$ and all $u\ $in $\Pi U\Pi
^{-1}$, so that $f_{1}\left( gu\right) =u^{-1}f_{1}(g)$ for every $u\ $in $%
SL_{2}(F_{\mathfrak{P}_{1}})U\subset GL_{2}(\mathbb{A}_{F}^{\infty })$. Here 
$SL_{2}(F_{\mathfrak{P}_{1}})$ acts on $W_{\tau _{1}}$ trivially.

Assume that $W_{\tau }^{U}\neq \{0\}$, i.e., that $W_{\tau }=\mathbb{F}$ is
the trivial representation of $U$. If $\left( f_{1},f_{2}\right) \in \ker
\alpha $, then $f_{1}$ is invariant under right translations by elements of $%
D^{\times }U$; strong approximation for $SL_{2}$ then implies that $f_{1}$
is invariant under right translations by any element of $SL_{2}(\mathbb{A}%
_{F}^{\infty })$, and hence it factors through the reduced norm map $%
D^{\times }\backslash \left( D\otimes _{F}\mathbb{A}_{F}^{\infty }\right)
^{\times }\mathbb{\rightarrow }F^{\times }\backslash \left( \mathbb{A}%
_{F}^{\infty }\right) ^{\times }$. Since any maximal ideal of $\mathbb{T}_{S,%
\mathbb{F}}^{univ}$ in the support of the space of functions $D^{\times
}\backslash \left( D\otimes _{F}\mathbb{A}_{F}^{\infty }\right) ^{\times }%
\mathbb{\rightarrow }W_{\tau }$ factoring through the reduced norm is
Eisenstein, we obtain the desired result.

Assume now that $W_{\tau }^{U}=\{0\}$ and let $\left( f_{1},f_{2}\right) \in
\ker \alpha $. Using strong approximation, we see that for any $g\in GL_{2}(%
\mathbb{A}_{F}^{\infty })$ and $u\in \tprod\nolimits_{j=1}^{r}GL_{2}(%
\mathcal{O}_{F_{\mathfrak{P}_{j}}})$ we can find an element $\delta \in
D^{\times }\cap gSL_{2}(F_{\mathfrak{P}_{1}})Ug^{-1}$ such that for all $%
j=1,...,r:$%
\begin{equation*}
g_{\mathfrak{P}_{j}}^{-1}\delta g_{\mathfrak{P}_{j}}\in u_{\mathfrak{P}%
_{j}}+M_{2}(\mathfrak{P}_{j}).
\end{equation*}

\noindent In particular, we obtain:%
\begin{equation*}
f_{1}(g)=f_{1}(\delta ^{-1}g)=f_{1}(g(g^{-1}\delta ^{-1}g))
\end{equation*}%
and, since $g^{-1}\delta ^{-1}g\in SL_{2}(F_{\mathfrak{P}_{1}})U$:%
\begin{equation*}
f_{1}(g)=\left( g^{-1}\delta g\right) f_{1}(g)=uf_{1}(g)\text{.}
\end{equation*}

\noindent Since $u$ is arbitrary, we conclude that $f_{1}(g)\in W_{\tau
}^{U} $ for any $g\in GL_{2}(\mathbb{A}_{F}^{\infty })$, so that $f_{1}=0$, $%
f_{2}=0$ and $\alpha $ is injective. $\blacksquare $

\bigskip

Let $\mathcal{F}_{\tau }$ denote the space consisting of all the functions $%
f:D^{\times }\backslash \left( D\otimes _{F}\mathbb{A}_{F}^{\infty }\right)
^{\times }\mathbb{\rightarrow }W_{\tau }$, and define a left $\mathbb{F}$%
-linear\ action of $U$ on $\mathcal{F}_{\tau }$ by:%
\begin{equation*}
(uf)(g)=\tau (u)f(gu)
\end{equation*}

\noindent for all $u\in U$, $g\in \left( D\otimes _{F}\mathbb{A}_{F}^{\infty
}\right) ^{\times }$ and $f\in \mathcal{F}_{\tau }$. Set:%
\begin{equation*}
S_{\tau }(U,\mathbb{F})=H^{0}(U,\mathcal{F}_{\tau }).
\end{equation*}

In what follows, we work for simplicity with the spaces $S_{\tau }(U,\mathbb{%
F})$, forgetting about the action of the center of $\left( D\otimes _{F}%
\mathbb{A}_{F}^{\infty }\right) ^{\times }$ on $\mathcal{F}_{\tau }$.

Following the proof of Proposition 1 at page 48 of \cite{EK}, and using
Lemma \ref{dege}, we obtain the following result:

\begin{theorem}
\label{holomorphic 2}Assume that $\tau $ is an irreducible (non necessarily
holomorphic) $\mathbb{F}$-linear weight with parameters $(\vec{k},\vec{w}%
)\in 
%TCIMACRO{\U{2124} }%
%BeginExpansion
\mathbb{Z}
%EndExpansion
_{\geq 2}^{g}\times 
%TCIMACRO{\U{2124} }%
%BeginExpansion
\mathbb{Z}
%EndExpansion
^{g}$ such that $\vec{k}^{(j)}=\vec{2}$ for some $1\leq j\leq r$. Let $\tau
^{\prime }$ be the $\mathbb{F}$-linear weight associated to the parameters $%
\vec{k}^{\prime }=(\vec{k}^{(1)},...,\vec{k}^{(j)}+\overrightarrow{p-1},...,%
\vec{k}^{(r)})$ and $\vec{w}^{\prime }=\vec{w}$. For any non-Eisenstein
maximal ideal $\mathfrak{M}$ of $\mathbb{T}_{S,\mathbb{F}}^{univ}$, there is
an injective Hecke-equivariant $\mathbb{F}$-morphism:%
\begin{equation*}
S_{\tau }(U,\mathbb{F})_{\mathfrak{M}}\hookrightarrow S_{\tau ^{\prime }}(U,%
\mathbb{F})_{\mathfrak{M}}.
\end{equation*}
\end{theorem}

\textbf{Proof }Assume without loss of generality that $j=1$. Via the
surjection $U\rightarrow GL_{2}(\mathbb{F}_{\mathfrak{P}_{1}})$,\ the group $%
U$ acts on the $\mathbb{F}_{\mathfrak{P}_{1}}$-points $\mathbb{P}^{1}(%
\mathbb{F}_{\mathfrak{P}_{1}})$ of the projective $\mathbb{F}_{p}$-line, and
we can identify the coset space $U/U_{0}$ with $\mathbb{P}^{1}(\mathbb{F}_{%
\mathfrak{P}_{1}})$. Recall that we are viewing $\mathbb{F}_{\mathfrak{P}%
_{1}}$ as a subfield of $\mathbb{F}$ via the fixed embedding $\bar{\sigma}%
_{0}^{(1)}$.

By Shapiro's lemma applied to the pair $(U,U_{0})$ and the left $\mathbb{F}%
[U]$-module $\mathcal{F}_{\tau }$, we obtain an isomorphism:%
\begin{equation}
H^{0}(U_{0},\mathcal{F}_{\tau })\overset{\sim }{\longrightarrow }H^{0}(U,%
\mathcal{F}_{\tau }\mathcal{\otimes }_{\mathbb{F}}\mathbb{F}[\mathbb{P}^{1}(%
\mathbb{F}_{\mathfrak{P}_{1}})]).  \tag{1}  \label{shapiroo}
\end{equation}

\noindent Here $U$ acts on $\mathbb{F}[\mathbb{P}^{1}(\mathbb{F}_{\mathfrak{P%
}_{1}})]=\{\varphi :\mathbb{P}^{1}(\mathbb{F}_{\mathfrak{P}_{1}})\mathbb{%
\rightarrow }\mathbb{F}\}$ via its quotient $GL_{2}(\mathbb{F}_{\mathfrak{P}%
_{1}})$\ and by the rule $(u\varphi )(P)=\varphi (u^{-1}P)$ for $u\in GL_{2}(%
\mathbb{F}_{\mathfrak{P}_{1}})$ and $P\in \mathbb{P}^{1}(\mathbb{F}_{%
\mathfrak{P}_{1}})$. Furthermore $U$ acts diagonally on $\mathcal{F}_{\tau }%
\mathcal{\otimes }\mathbb{F}[\mathbb{P}^{1}(\mathbb{F}_{\mathfrak{P}_{1}})]$%
. By Lemma 1.1.4 of \cite{AS1}, the isomorphism (\ref{shapiroo}) preserves
the Hecke action on both sides.

By Lemma 2.6 of \cite{Re}, there is an isomorphism of $\mathbb{F}[GL_{2}(%
\mathbb{F}_{\mathfrak{P}_{1}})]$-modules:%
\begin{equation*}
\mathbb{F}[\mathbb{P}^{1}(\mathbb{F}_{\mathfrak{P}_{1}})]\simeq \mathbb{%
F\oplus }\limfunc{Sym}\nolimits^{p^{f_{1}}-1}\left( \mathbb{F}^{2}\right)
=M_{0}\oplus M_{p^{f_{1}}-1}\text{,}
\end{equation*}

\noindent inducing a surjection:%
\begin{equation}
H^{0}(U,\mathcal{F}_{\tau }\mathcal{\otimes }\mathbb{F}[\mathbb{P}^{1}(%
\mathbb{F}_{\mathfrak{P}_{1}})])\longrightarrow H^{0}(U,\mathcal{F}_{\tau }%
\mathcal{\otimes }M_{p^{f_{1}}-1}).  \tag{2}  \label{oneone}
\end{equation}

Observe that the composition of the restriction map $H^{0}(U,\mathcal{F}%
_{\tau })\rightarrow H^{0}(U_{0},\mathcal{F}_{\tau })$ with the surjection:%
\begin{equation*}
H^{0}(U_{0},\mathcal{F}_{\tau })\simeq H^{0}(U,\mathcal{F}_{\tau }\mathcal{%
\otimes }\mathbb{F}[\mathbb{P}^{1}(\mathbb{F}_{\mathfrak{P}%
_{1}})])\rightarrow H^{0}(U,\mathcal{F}_{\tau })
\end{equation*}

\noindent is given by $f\mapsto \frac{1}{[U:U_{0}]}\tsum\nolimits_{u\in
U/U_{0}}1\otimes uf=1\otimes f$. This implies that the first summand of $%
H^{0}(U,\mathcal{F}_{\tau })^{\oplus 2}$ is identified via the map $\alpha $%
\ of Lemma \ref{dege} and the Shapiro isomorphism with the direct summand $%
H^{0}(U,\mathcal{F}_{\tau })$ of $H^{0}(U,\mathcal{F}_{\tau }\mathcal{%
\otimes }\mathbb{F}[\mathbb{P}^{1}(\mathbb{F}_{\mathfrak{P}_{1}})])$.

Using the map $\alpha $, the Shapiro isomorphism, the projection (\ref%
{oneone}), and the isomorphism of $\mathbb{F}[GL_{2}(\mathbb{F}%
_{p^{f_{1}}})] $-modules $M_{p^{f_{1}}-1}\simeq
\tbigotimes\nolimits_{i=0}^{f_{1}-1}M_{p-1}^{[i]}$, we obtain a Hecke
equivariant morphism:%
\begin{equation*}
\beta :H^{0}(U,\mathcal{F}_{\tau })^{\oplus 2}\longrightarrow H^{0}\left( U,%
\mathcal{F}_{\tau }\mathcal{\otimes }\tbigotimes%
\nolimits_{i=0}^{f_{1}-1}M_{p-1}^{[i]}\right) .
\end{equation*}

By Lemma \ref{dege}, precomposing $\beta $ with the injection $H^{0}(U,%
\mathcal{F}_{\tau })\hookrightarrow H^{0}(U,\mathcal{F}_{\tau })^{\oplus 2}$
given by $f\mapsto (0,f)$ we obtain a Hecke equivariant injective morphism:%
\begin{equation}
H^{0}(U,\mathcal{F}_{\tau })_{\mathfrak{M}}\hookrightarrow H^{0}\left( U,%
\mathcal{F}_{\tau }\mathcal{\otimes }\tbigotimes%
\nolimits_{i=0}^{f_{1}-1}M_{p-1}^{[i]}\right) _{\mathfrak{M}}  \tag{3}
\label{one}
\end{equation}

\noindent for any non-Eisenstein maximal ideal $\mathfrak{M}$ of $\mathbb{T}%
_{S,\mathbb{F}}^{univ}$.

Let $\vec{k}^{\prime }=(\overrightarrow{p+1},\vec{k}^{(2)},...,\vec{k}%
^{(r)}) $ and set $\vec{w}^{\prime }=\vec{w}$. Observe that if $\tau
^{\prime }$ is the representation of $U$ associated to the parameters $(\vec{%
k}^{\prime },\vec{w}^{\prime })$ then $W_{\tau ^{\prime }}\simeq W_{\tau }%
\mathcal{\otimes }_{\mathbb{F}}\tbigotimes%
\nolimits_{i=0}^{f_{1}-1}M_{p-1}^{[i]}$. The $U$-equivariant map $\mathcal{F}%
_{\tau }\mathcal{\otimes }\tbigotimes\nolimits_{i=0}^{f_{1}-1}M_{p-1}^{[i]}%
\rightarrow \mathcal{F}_{\tau ^{\prime }}$ induced by the assignment:%
\begin{equation*}
f\otimes m\longmapsto \lbrack g\mapsto f(g)\otimes m]
\end{equation*}

\noindent for $g\in D^{\times }\backslash \left( D\otimes _{F}\mathbb{A}%
_{F}^{\infty }\right) ^{\times }$ is injective. We deduce that for any
non-Eisenstein maximal ideal $\mathfrak{M}$ of $\mathbb{T}_{S,\mathbb{F}%
}^{univ}$, there is a Hecke equivariant monomorphism: 
\begin{equation*}
H^{0}\left( U,\mathcal{F}_{\tau }\mathcal{\otimes }\tbigotimes%
\nolimits_{i=0}^{f_{1}-1}M_{p-1}^{[i]}\right) _{\mathfrak{M}}\hookrightarrow
H^{0}(U,\mathcal{F}_{\tau ^{\prime }})_{\mathfrak{M}}.
\end{equation*}%
Combining this with (\ref{one}), we are done. $\blacksquare $

\bigskip

\begin{remark}
Under the assumptions of the above theorem, $\tau ^{\prime }$ is an
irreducible representation of $U$. This implies that, if the number of
indices $j$ such that $\vec{k}^{(j)}=\vec{2}$ is larger than one, Theorem %
\ref{holomorphic 2} can be further applied to obtain weight shiftings "in
blocks" by $\overrightarrow{p-1}$.
\end{remark}

\begin{remark}
\label{SerreF}The content of Theorem \ref{holomorphic 2} generalizes Lemma
4.6.8 of \cite{Gee22}, which is proved in loc. cit. via Lemma 1.5.5 of \cite%
{Kisin3}.
\end{remark}

The weight shifting produced by Theorem \ref{holomorphic 2}\ is not in
general of holomorphic type: for example, if $r>1$ and $\tau $ is
holomorphic, then $\tau ^{\prime }$ is never holomorphic. Nevertheless we
have:

\begin{corollary}
\label{HMF case}Assume that $\tau $ is the irreducible holomorphic $\mathbb{F%
}$-linear weight with parameters $(\vec{2},w)\in 
%TCIMACRO{\U{2124} }%
%BeginExpansion
\mathbb{Z}
%EndExpansion
_{\geq 2}^{g}\times \left( 2%
%TCIMACRO{\U{2124} }%
%BeginExpansion
\mathbb{Z}
%EndExpansion
+1\right) $. Let $\tau ^{\prime }$ be the holomorphic weight associated to
the parameters $(\overrightarrow{p+1},w+(p-1))\in 
%TCIMACRO{\U{2124} }%
%BeginExpansion
\mathbb{Z}
%EndExpansion
_{\geq 2}^{g}\times \left( 2%
%TCIMACRO{\U{2124} }%
%BeginExpansion
\mathbb{Z}
%EndExpansion
+1\right) $. For any non-Eisenstein maximal ideal $\mathfrak{M}$ of $\mathbb{%
T}_{S,\mathbb{F}}^{univ}$, there is an injective Hecke-equivariant $\mathbb{F%
}$-morphism:%
\begin{equation*}
S_{\tau }(U,\mathbb{F})_{\mathfrak{M}}\hookrightarrow S_{\tau ^{\prime }}(U,%
\mathbb{F})_{\mathfrak{M}}.
\end{equation*}
\end{corollary}

\textbf{Proof }Fix a non-Eisenstein maximal ideal $\mathfrak{M}$ of $\mathbb{%
T}_{S,\mathbb{F}}^{univ}$. Applying Theorem \ref{holomorphic 2} $r$ times we
obtain a Hecke equivariant injection $S_{\tau }(U,\mathbb{F})_{\mathfrak{M}%
}\hookrightarrow S_{\tau ^{\prime }}(U,\mathbb{F})_{\mathfrak{M}}$, where $%
\tau ^{\prime }$ is the irreducible $\mathbb{F}$-linear weight with
parameters $(\overrightarrow{p+1},\vec{v})\in 
%TCIMACRO{\U{2124} }%
%BeginExpansion
\mathbb{Z}
%EndExpansion
_{\geq 2}^{g}\times 
%TCIMACRO{\U{2124} }%
%BeginExpansion
\mathbb{Z}
%EndExpansion
^{g}$ and each component of $\vec{v}$ equals the integer $\frac{w-1}{2}$.
This weight is holomorphic with parameters $(\overrightarrow{p+1},w+(p-1))$. 
$\blacksquare $

\bigskip

In terms of Galois representations we obtain:

\begin{corollary}
Assume that $\tau $ is an irreducible (non necessarily holomorphic) $\mathbb{%
F}$-linear weight with parameters $(\vec{k},\vec{w})\in 
%TCIMACRO{\U{2124} }%
%BeginExpansion
\mathbb{Z}
%EndExpansion
_{\geq 2}^{g}\times 
%TCIMACRO{\U{2124} }%
%BeginExpansion
\mathbb{Z}
%EndExpansion
^{g}$ such that $\vec{k}^{(j)}=\vec{2}$ for some $1\leq j\leq r$. Then an
irreducible continuous representation $\rho :\limfunc{Gal}(\bar{F}/F)\mathbb{%
\rightarrow }GL_{2}(\mathbb{\bar{F}}_{p})$ arising from an automorphic
eigenform on $S_{\tau }(U,\mathbb{F})$ also arises from an automorphic
eigenform on $S_{\tau ^{\prime }}(U,\mathbb{F})$, where $\tau ^{\prime }$ is
the irreducible weight associated to the parameters $\vec{k}^{\prime }=(\vec{%
k}^{(1)},...,\vec{k}^{(j)}+\overrightarrow{p-1},...,\vec{k}^{(r)})$ and $%
\vec{w}^{\prime }=\vec{w}$.
\end{corollary}

The Jacquet-Langlands correspondence and Corollary \ref{HMF case} imply the
following (cf. \cite{EK}):

\begin{corollary}
An irreducible continuous representation $\rho :\limfunc{Gal}(\bar{F}/F)%
\mathbb{\rightarrow }GL_{2}(\mathbb{\bar{F}}_{p})$ arising from a
holomorphic Hilbert modular form of level $U\subset GL_{2}(\mathbb{A}%
_{F}^{\infty })$ and parallel weight $\vec{2}$ also arises from a
holomorphic Hilbert modular form of level $U$ and parallel weight $%
\overrightarrow{p+1}$.
\end{corollary}

\bibliographystyle{amsplain}
\bibliography{davide2}
\bigskip

\end{document}